\documentclass{ruthesis} 

\usepackage{epsfig}  		
\usepackage{epic,eepic}       

\usepackage{amsmath,amsthm,amssymb,latexsym,url,graphicx,cite}
\usepackage{hyperref}
\usepackage{enumerate}
\usepackage{mathbbol}
\usepackage{mathrsfs}
\usepackage{longtable}
\usepackage[small]{diagrams}
\usepackage[figurewithin=section]{caption}

\DeclareSymbolFontAlphabet{\mathbb}{AMSb}
\DeclareSymbolFontAlphabet{\mathbbl}{bbold}

\newenvironment{claimproof}[1][Proof of Claim]{\begin{proof}[#1]}{\end{proof}}

\newtheorem{lemma}{Lemma}[section]
\newtheorem*{lemma*}{Lemma}
\newtheorem{theorem}[lemma]{Theorem}
\newtheorem*{theorem*}{Theorem}
\newtheorem{corollary}[lemma]{Corollary}
\newtheorem{proposition}[lemma]{Proposition}

\newtheorem*{proposition*}{Proposition}

\newtheorem*{fact*}{Fact}
\newtheorem{notation}[lemma]{Notation}
\newtheorem*{notation*}{Notation}
\newtheorem*{conventions*}{Conventions}
\newtheorem{remark}[lemma]{Remark}
\newtheorem*{remark*}{Remark}
\newtheorem{lemmacorollary}[lemma]{Corollary}
\newtheorem*{corollary*}{Corollary}
\newtheorem{conjecture}{Conjecture}
\newtheorem*{conjecture*}{Conjecture}
\newtheorem{problem}{Problem}
\newtheorem*{problem*}{Problem}
\newtheorem{question}{Question}
\newtheorem*{question*}{Question}

\newtheorem{assumption*}{Assumption}

\newtheorem{theoremC}{Theorem}[chapter]
\newtheorem{definitionC}{Definition}[chapter]

\theoremstyle{definition}
\newtheorem{example}{Example}
\newtheorem*{example*}{Example}
\newtheorem{definition}[lemma]{Definition}
\newtheorem*{definition*}{Definition}

\theoremstyle{remark}
\newtheorem{claim}{Claim}
\newtheorem*{claim*}{Claim}

\newtheorem*{construction*}{Construction}

\newtheorem*{exercise*}{Exercise}

\numberwithin{equation}{section}

\newcommand{\N}{\mathbb{N}}

\newcommand{\Q}{\mathbb{Q}}

\newcommand{\set}[1]{\left\{ #1 \right\}}

\newcommand{\bs}{\backslash}
\newcommand{\ceil}[1]{\left \lceil #1 \right \rceil }
\newcommand{\floor}[1]{\left \lfloor #1 \right \rfloor}
\newcommand{\meet}{\wedge}
\newcommand{\join}{\vee}

\renewcommand\AA{{\mathcal A}}

\newcommand\CC{{\mathcal C}}

\newcommand\EE{{\mathcal E}}

\newcommand\GG{{\mathcal G}}
\newcommand\HH{{\mathcal H}}

\newcommand\KK{{\mathcal K}}
\newcommand\LL{{\mathcal L}}
\newcommand\MM{{\mathcal M}}

\newcommand\QQ{{\mathcal Q}}
\newcommand\RR{{\mathcal R}}

\newcommand\TT{{\mathcal T}}
\newcommand\UU{{\mathcal U}}

\newcommand\bbzero{\mathbbl 0} 
\newcommand\bbone{\mathbbl 1}  

\newcommand\into{\hookrightarrow}
\newcommand\includedin{\subseteq}

\newcommand\union{\cup}

\newcommand\fraisse{Fra\"\i ss\'e }

\newcommand{\rightcurve}[1]{%
\setlength{\unitlength}{0.03\DiagramCellWidth}
\begin{picture}(0,0)(0,0)
\qbezier(-58,3)(40,-10)(65,130)
\put(40,30){\makebox(0,0)[t]{$\scriptstyle {#1}$}}
\end{picture}
}

\newcommand{\leftcurve}[1]{%
\setlength{\unitlength}{0.03\DiagramCellWidth}
\begin{picture}(0,0)(0,0)
\qbezier(-65,130)(-40,-10)(58,3)
\put(-60,30){\makebox(0,0)[t]{$\scriptstyle {#1}$}}
\end{picture}
}
\newarrow {Line}          -----
\newarrow{Dashto}       {}{dash}{}{dash}>
\newarrow{Dash}           {}{dash}{}{dash}{}
\newarrow{Double}<--->

\newcommand\op{\text{op}}
\newcommand\tp{\text{tp}}

\DeclareMathOperator{\Lift}{Lift}

\newcommand{\arr}{\xrightarrow{E}}
\newcommand{\ra}{\xrightarrow}
\newcommand{\la}{\xleftarrow}

\newcommand{\To}{{\Rightarrow}}
\newcommand{\From}{{\Leftarrow}}

\usepackage{amsrefs}
\usepackage{setspace}

\begin{document} 
\phd 
\title{Infinite Limits of Finite-Dimensional Permutation Structures, and their Automorphism Groups: Between Model Theory and Combinatorics} 
\author{Samuel Walker Braunfeld} 
\program{Mathematics} 
\director{Gregory Cherlin} 
\approvals{4} 
\submissionyear{2018} 
\submissionmonth{May}
\abstract{In the course of classifying the homogeneous permutations, Cameron introduced the viewpoint of permutations as structures in a language of two linear orders \cite{Cameron}, and this structural viewpoint is taken up here. The majority of this thesis is concerned with Cameron's problem of classifying the homogeneous structures in a language of finitely many linear orders, which we call finite-dimensional permutation structures. Towards this problem, we present a construction that we conjecture produces all such structures. Some evidence for this conjecture is given, including the classification of the homogeneous 3-dimensional permutation structures.

We next consider the topological dynamics, in the style of Kechris, Pestov, and Todor\v{c}evi\'c, of the automorphism groups of the homogeneous finite-dimensional permutation structures we have constructed, which requires proving a structural Ramsey theorem for all the associated amalgamation classes. Because the $\emptyset$-definable equivalence relations in these homogeneous finite-dimensional permutation structures may form arbitrary finite distributive lattices, the model-theoretic algebraic closure operation may become quite complex, and so we require the framework recently introduced by Hubi\v{c}ka and Ne{\v{s}}etril \cite{HN}.

Finally, we turn to the interaction of model theory with more classical topics in the theory of permutation avoidance classes. We consider the decision problem for whether a finitely-constrained permutation avoidance class is atomic, or equivalently, has the joint embedding property. As a first approximation to this problem, we prove the undecidability of the corresponding decision problem in the category of graphs. Modifying this proof also gives the undecidability, in the category of graphs, of the corresponding decision problem for the joint homomorphism property, which is of interest in infinite-domain constraint satisfaction problems.

The results in the first 8 chapters of this thesis largely appeared in the previous articles \cite{Lattice}, \cite{RamExp}, and \cite{3dim}. In many places the arguments and context have been expanded upon, and in the case of some arguments from \cite{Lattice}, they have been simplified.
} 
\beforepreface 
\acknowledgements{
The first 8 chapters of this thesis largley appeared in the previous articles \cite{Lattice}, \cite{RamExp}, and \cite{3dim}. In particular, Chapters 4 and 5 mostly appeared in \cite{Lattice}, Chapter 6 in \cite{3dim}, Chapter 8 in \cite{RamExp}, and parts of Chapter 3 appeared in all 3 articles.

I thank my committee members, Gregory Cherlin, Grigor Sargsyan, Simon Thomas, and Michael Laskowski, for looking over my thesis and appearing at my defense, and I further thank Professors Sargsyan and Thomas for their teaching and availability through my years at Rutgers.

I thank Manuel Bodirsky for suggesting I consider the decision problem for the joint homorphism property, and Jan Hubi\v{c}ka for various discussions, particularly concerning the framework behind Chapter 8.

Most of all, I thank my advisor, Gregory Cherlin. For countless hours of discussion, for your sense of humor, and for guiding me toward becoming a mathematician and imparting your sense of the subject.

}
\dedication{
\begin{center}
\vspace{2 cm}
\textit{\Large To my parents}
\end{center}
}
\figurespage

\tablespage

\afterpreface 

\chapter{Introduction}



Homogeneous structures have provided a conduit from model theory to combinatorics at least since the work of Lachlan on the classification of finite homogeneous structures, which uncovered deep connections with stability theory. In the other direction, the generalization of this theory to smoothly approximable structures, suggested by Lachlan, led to the introduction of some fundamental notions in what is now known as neostability theory, in joint work by Cherlin and Hrushovski. The subject of homogeneity, and more generally $\omega$-categoricity, is  experiencing a golden age of activity at present, with applications to algorithmic problems, notably constraint satisfaction problems
in the work of Bodirsky and others, links to topological dynamics through the influential work of Kechris, Pestov,
and Todor\v{c}evi\'c, and the classification of stable homogeneous structures being continued in the neostable setting. 

Combinatorial applications of model theory bring the study of infinite structures to bear on the theory of finite structures, generally by taking an appropriate limit such as an ultraproduct or homogeneous \fraisse limit, or in some cases
suggest fruitful analogies between the two. One underdeveloped area from this perspective is the large and active subject
of permutation avoidance classes. Here a start has been made by Cameron \cite{Cameron}, who observed that the category of permutations
and embeddings is best viewed as the category of structures with two linear orders, in which case isomorphism types will be
permutation patterns. He then classified the homogeneous permutations,
or equivalently the permutation avoidance classes with the amalgamation property and joint embedding. He also posed the analogous
problems for generalized permutation avoidance classes, consisting of structures with a specified finite number of linear orders,
which, as Waton lays out in some detail in the final chapter of his thesis \cite{Waton}, are the natural higher dimensional analog of permutations.

\begin{problem}
Classify, for each $m$, the homogeneous structures in a language consisting of $m$ linear orders.
\end{problem}

This thesis is concerned mainly with Cameron's problem concerning homogeneous finite-dimensional permutation classes
and related topics, such as the associated structural Ramsey theory, but looks toward a richer interaction of model theory
with more classical topics in the theory of permutation avoidance classes, notably associated algorithmic problems.

In addition to Cameron's problem, another motivation for this classification arose from the classification of the homogeneous ordered graphs. The classification of various homogeneous structures equipped with a linear order is of particular importance for identifying new examples in structural Ramsey theory, and it was for this reason the classification of homogeneous ordered graphs was undertaken.  However, after a lengthy classification, no new examples were uncovered, as it turned out every homogeneous ordered graph arises in simple fashion from an unordered homogeneous structure (either a graph, tournament, or partial order) \cite{MHG}. To be more precise, these structures are interdefinable with generic expansions by a linear order of homogeneous graphs or tournaments, or generic linear extensions of homogeneous partial orders. It is natural to ask whether a similar statement might hold more generally, which would greatly simplify classifying homogeneous ordered structures.

\begin{question} \label{q:ordexp}
Is every homogeneous ordered structure interdefinable with an expansion of a homogeneous proper reduct by a linear order?

Is every \textbf{primitive} homogeneous ordered structure that is the \fraisse limit of a strong amalgamation class interdefinable with an expansion of a homogeneous proper reduct by a \textbf{generic} linear order?
\end{question}

The minimal case to test this would be to start with a structureless set, iteratively add linear orders, and observe the homogeneous structures that appear, which brings us back to investigating the homogeneous finite-dimensional permutation structures.

The first step toward a classification is the production of a catalog, or census, of examples occurring ``in nature.'' This is undertaken in Chapters \ref{chap:lambda}-\ref{chap:construction}, which include an amalgamation construction producing all known homogeneous finite-dimensional permutation structures, including many new imprimitive examples. 

That construction is based on \textit{$\Lambda$-ultrametric spaces}, which are of independent interest, and provide a presentation of structures consisting of equivalence relations using an analog of a metric taking values in a lattice $\Lambda$. The analogy with metric spaces provides an amalgamation procedure iff $\Lambda$ is distributive. Recently, Hubi\v{c}ka, Kone{\v{c}}n{\'y}, and Ne{\v{s}}et\v{r}il have introduced a common generalization of $\Lambda$-ultrametric spaces and Conant's generalized metric spaces \cite{ConMet}, which seems to provide a unified treatment of most homogeneous structures in a binary language with constraints of size at most 3 \cite{HKN2}.

Our construction of homogeneous finite-dimensional permutation structures proceeds roughly as follows. One starts with a fully generic $\Lambda$-ultrametric space. This structure is then expanded by linear orders so that every equivalence relation is convex with respect to at least one definable order, and the equivalence relations are then interdefinably exchanged for additional linear orders. However, we do not work directly with linear orders, but more generally with certain partial orders which we call \textit{subquotient orders}. This allows our expansion to be fully generic in a natural sense. We then introduce the notion of a \textit{well-equipped lift}, which captures when a set of subquotient orders is interdefinable with a set of linear orders. Our catalog is then as follows.

\begin{theoremC} \label{thm:introconstruction}
For every finite distributive lattice $\Lambda$, any well-equipped lift of the class of all finite $\Lambda$-ultrametric spaces is an amalgamation class, and its \fraisse limit is interdefinable with a finite-dimensional permutation structure.
\end{theoremC}  

Once a coherent catalog has been found, it is standard to present it as a conjectural classification, although the appearance of further sporadic examples would be unsurprising. Chapters \ref{chap:completeness}-\ref{c:3dim} are then concerned with the resulting conjecture.

\begin{conjecture} \label{conj:sqoconj}
The construction from Theorem \ref{thm:introconstruction} produces all homogeneous finite-dimensional permutation structures.
\end{conjecture}

It is also worth separately stating the primitive case of this conjecture. Two orders $<_1, <_2$ are said to be \textit{identified, up to reversal} if $<_1 = <_2$ or $<_1^{opp} = <_2$.

\begin{conjecture}[Primitivity Conjecture] \label{conj:primitivity}
Every primitive homogeneous finite dimensional permutation structure can be constructed by the following procedure.
\begin{enumerate}
\item Identify certain orders, up to reversal.
\item Take the \fraisse limit of the resulting amalgamation class, getting a fully generic structure, possibly in a simpler language.
\end{enumerate}
\end{conjecture} 

In personal communication just before the submission of this thesis, Pierre Simon confirmed a proof of the Primitivity Conjecture \cite{Simon}, as an elaboration on the work presented in the talk \cite{Simon2}.

Each of Chapters \ref{chap:completeness} and \ref{c:3dim} provides some evidence for the correctness of Conjecture \ref{conj:sqoconj}. In particular, in Chapter \ref{c:3dim} we classify the homogeneous 3-dimensional permutation structures, and obtain the following result.

\begin{theoremC}
   Every homogeneous 3-dimensional permutation structure is interdefinable with the \fraisse limit of some well-equipped lift of the class of all finite $\Lambda$-ultrametric spaces, for some distributive lattice $\Lambda$.
\end{theoremC}

There is little general theory that can be applied in such a classification. There is the fact that in a binary language, whether a finitely-constrained class is an amalgamation class is decidable from the forbidden substructures. There is also a Ramsey argument due to Lachlan \cite{LW}, of seemingly broad applicability, but amalgamation of linear orders is so constrained that we may proceed more directly.

In Chapter \ref{chap:Ramsey}, we consider the structural Ramsey property for finite-dimensional permutation structures and $\Lambda$-ultrametric spaces (for an introduction to structural Ramsey theory and the related topological dynamics, see Chapter \ref{section:RamseyBackground}). The step in our construction in which linear orders are added so that every $\emptyset$-definable equivalence relation becomes convex with respect to one such looks suspiciously like what one would try in order to produce a Ramsey expansion of the generic $\Lambda$-ultrametric space. This prompts the question whether the classes thus produced are in fact Ramsey classes, and the related question of describing the universal minimal flow of the automorphism group of the generic $\Lambda$-ultrametric space. 

\begin{theoremC}
Let $\Lambda$ be a finite distributive lattice and $\Gamma$ be the generic $\Lambda$-ultrametric space. For every meet-irreducible $E \in \Lambda$, expand $\Gamma$ by a generic subquotient order from $E$ to its successor, let $\vec \Gamma^{min} = (\Gamma, (<_{E_i})_{i=1}^n)$ be structure thus obtained, and $\vec \AA_\Lambda^{\min}$ its finite substructures. Then 
\begin{enumerate}
\item $\vec \AA_\Lambda^{\min}$ is a Ramsey class and has the expansion property relative to $\AA_\Lambda$.
\item The logic action of $\text{Aut}(\Gamma)$ on $\overline{\text{Aut}(\Gamma) \cdot (<_{E_i})_{i=1}^n}$ is the universal minimal flow of $\text{Aut}(\Gamma)$.
\end{enumerate}
\end{theoremC}
 
 The theorem above gives an explicit description of the universal minimal flow of $\text{Aut}(\Gamma)$ as isomorphic to the logic action on the full product of the following factors indexed by $i \in [n]$, where $n$ is the number of meet-irreducibles in $\Lambda$: the space of linear orders on $\Gamma$ if $<_{E_i}$ satisfies a certain condition (its top relation is $\bbone$), and otherwise the full infinite Cartesian power of that space. In particular, the universal minimal flow is metrizable.
 
 Most of the work in the proof of this theorem is directed toward establishing the Ramsey property, which we do for all the homogeneous finite-dimensional permutation structures produced by our construction.
 
  \begin{theoremC}
	All classes produced by the construction from Theorem \ref{thm:introconstruction} are Ramsey.
  \end{theoremC}
  
  This Ramsey theorem is proven using tools from Hubi\v{c}ka and Ne{\v{s}}etril \cite{HN}. In particular, we use a combination and generalization of various encoding techniques that were used there in the proofs of Ramsey theorems for the free product of Ramsey classes and for structures that have a chain of definable equivalence relations. 
  
  It is interesting to note that our original version of the main construction from Theorem \ref{thm:introconstruction}, given in \cite{Lattice}, used linear orders rather than subquotient orders, and we did not give a satisfactory catalog of the homogeneous finite-dimensional permutation structures there. Although, in hindsight, subquotient orders did implicitly appear in that paper, they only came to the foreground, allowing us to complete our goal of a satisfactory catalog, while considering the Ramsey theory of these structures.
  
  Chapter \ref{chap:JEP} is not concerned with the amalgamation or Ramsey properties, but instead deals with the much weaker joint-embedding property. Given a hereditary class, the joint-embedding property is equivalent to the existence of a countable structure that is finitely universal for that class, i.e. such that its finite substructures are the finite structures of the class.
  
  In the case of permutation avoidance classes, the joint-embedding property is equivalent to an indecomposability property called atomicity. We are inspired by the following question of Ru\v{s}kuc \cite{Rusk}.
  
  \begin{question} \label{question:permJEP}
  Is there an algorithm that, given finite set of forbidden permutations, decides whether the corresponding permutation avoidance class has the joint-embedding property?
  \end{question}
  
  This question is reminiscent of decidability questions that have been raised in model-theoretic contexts. The first is a question of Lachlan concerning the decidability of amalgamation given the forbidden substructures \cite{LachICM}, as well as the variants of this question presented in \cite{MHG}. The second is the question of the decidability of the existence of a countably universal countable graph given finitely many forbidden non-induced subgraphs, or equivalently whether the corresponding theory is small, as discussed in \cite{WQOU}.

  As a first step towards a possible solution to Ru\v{s}kuc's question, we work instead in the language of graphs. Via a reduction from the tiling problem, we obtain the following theorem.
  
\begin{theoremC} \label{introthm:jepUndecidable}
  There is no algorithm that, given a finite set of forbidden induced subgraphs, decides whether the corresponding hereditary graph class has the JEP.
\end{theoremC}

Proving this theorem requires representing a 2-dimensional grid in various hereditary graph classes. However, there seem to be obstacles to representing a grid in permutation avoidance classes, at least other than the class of all permutations, which prevent an analogous argument from being carried out there, although this issue seems to disappear in higher-dimensional permutation avoidance classes.

By modifying the proof of Theorem \ref{introthm:jepUndecidable}, we obtain an analogous result for the \textit{joint homomorphism property}, which is of significant interest in constraint satisfaction problems.

\begin{theoremC}
Work in a language with one binary relation, which will be interpreted as edge relation of a graph. Then there is no algorithm that, given a finite set of forbidden induced subgraphs, decides whether the corresponding hereditary graph class has the joint homomorphism property.
\end{theoremC} 


Throughout this introduction, we have been considering hereditary classes of structures and their limits. The hereditary condition is equivalent to the theory of the class having a particularly simple form, namely being a universal theory. For such theories, the existentially closed structures provide a natural class of limit objects. Although it seems that much of the machinery of model theory can be adapted to the class of existentially closed structures regardless, one may find it particularly attractive when this class is first-order axiomatizable, in which case the corresponding theory is called the model companion of the original theory. If we further have that the model companion is $\omega$-categorical, then we have a canonical countable limit structure for the class. Furthermore, in the case of graphs, where these issues have been extensively considered, the existence of an $\omega$-categorical model companion seems intimately tied to the existence of a countable structure that is \textit{countably universal} for the class \cite{WQOU}, a significantly stronger limit notion than the finite universality provided by the joint-embedding property. The model companion is automatically model-complete, i.e. every formula is equivalent to an existential formula, and if in addition to $\omega$-categoricity we demand full quantifier elimination, we recover the notion of \fraisse limit. 

Thus, investigating the existence and model-theoretic properties of the model companion for permutation avoidance classes provides a generalization of the \fraisse theory considered in this thesis, which we begin looking towards in the final chapter.  
\label{chap:intro}

\chapter{\fraisse Theory, Structural Ramsey Theory, and Topological Dynamics}\label{chap:background}
\section{Introduction}
In this chapter, we present the correspondence between hereditary classes of finite structures with the amalgamation property and homogeneous structures, as well as the correspondence between the structural Ramsey theory of such classes and the topological dynamics of the automorphism group of the corresponding homogeneous structure. A general framework for constructing homogeneous structures is given by Fra\"\i ss\'e's theory of amalgamation classes, reviewed in the next section. While this can be used to produce many examples of homogeneous structures, it also provides a means for classification, as the amalgamation property strongly constrains classes of structures and its implications can be analyzed quite directly.

Many of the notions around homogeneous structures and their automorphism groups can be considered in the broader class of $\omega$-categorical structures, which by Ryll-Nardzewski-Svenonius-Engeler theorem are precisely those structures whose automorphism groups are oligomorphic, i.e. induce finitely many orbits on $n$-tuples for each $n$. For more about homogeneity and $\omega$-categoricity, see Macpherson's survey \cite{Mac}.

In addition to results on classification, there has been a significant amount of work relating to the automorphism groups of homogeneous structures. This includes questions about the reconstruction of a homogeneous structure from its automorphism group, and properties such as \textit{automatic continuity} and the existence of \textit{generic automorphisms}.

Of particular note is the correspondence between  the topological dynamical property of \textit{extreme amenability} for the automorphism group and the structural Ramsey property for the corresponding amalgamation class. The first hint of this appeared in Pestov's proof of the extreme amenability of $(\Q, <)$ using the classical Ramsey theorem \cite{Pest}, and was fully developed in \cite{KPT}.

At the time, few examples of extremely amenable groups were known. Structural Ramsey theory had been considered at least since \cite{NR}, which introduced the powerful amalgamation-based \textit{partite-method} for proving the structural Ramsey property. This immediately yielded many new examples of extremely amenable groups and increased interest in structural Ramsey theory and the demand for interesting homogeneous structures to work on.

\section{Homogeneity Background}

\subsection{Fra\"\i ss\'e's Theorem}

For this chapter, let $L$ be a finite relational language and $M$ a countable $L$-structure.

\begin{definition}
$M$ is \textit{homogeneous} if any isomorphism between finite substructures of $M$ extends to an automorphism of $M$.
\end{definition}

Let $Age(M)$ be the class of finite $L$-structures isomorphic to a substructure of $M$. Note $Age(M)$ satisfies the following properties:
\begin{enumerate}[(i)]
\item  $Age(M)$ is closed under isomorphism and substructure.
\item $Age(M)$ has countably many isomorphism types.
\item Given, $B_1, B_2 \in Age(M)$, there is a $C \in Age(M)$ such that $B_1, B_2$ embed in $C$.
\end{enumerate}

The last of the above properties is called the \textit{joint-embedding property (JEP)}.

\begin{definition}
A class $\KK$ of finite $L$-structures has the \textit{amalgamation property} (and will be called an \textit{amalgamation class}) if, given $A, B_1, B_2 \in \KK$ with embeddings $f_i: A \to B_i$, there exist a $C \in \KK$ and embeddings $g_i: B_i \to C$ such that $g_1 \circ f_1 = g_2 \circ f_2$.
\end{definition}

\goodbreak

\begin{theorem*} [Fra\"\i ss\'e]
\hspace{.1cm}
\begin{enumerate}[(a)]
\item Let $M$ be a homogeneous structure. Then $Age(M)$ has the amalgamation property.
\item Let $\KK$ be a collection of finite $L$-structures satisfying $(i)-(iii)$ from above as well as the amalgamation property. Then, up to isomorphism, there is a unique countable, homogeneous $L$-structure $M$ with $Age(M) = \KK$.
\end{enumerate}
\end{theorem*}

The structure $M$ with $Age(M) = \KK$ from part $(b)$ is called the \textit{Fra\"\i ss\'e limit} of $\KK$.

Fra\"\i ss\'e's theorem provides a way to construct homogeneous structures by instead constructing amalgamation classes, which will be used significantly in Chapter \ref{chap:construction}. It also constrains homogeneous structures, as amalgamation is a strong condition whose implications can be concretely analyzed, and this aspect will be present in Chapters \ref{chap:completeness} and \ref{c:3dim}.

\subsection{Amalgamation Problems and Amalgamation Diagrams} \label{section:diagrams}
 In an amalgamation problem, one is asked to verify the amalgamation condition for specific structures $A, B_1, B_2$, and embeddings $f_i: A \hookrightarrow B_i$. This problem can be represented by an amalgamation diagram. 

In such a diagram, the points of the base $A$ are represented by points, while the points of $B_i \bs A$, which we call ``extension points'', are represented by circled points. We sometimes wish to depict an arbitrary finite set, in which case we use a large circle instead of individual points. The extension points of the first factor $B_1$ are placed on the left side of the diagram, while those of the second factor $B_2$ are placed on the right side.

When we are only considering binary languages, the relations are given by putting an arrow, or edge in the case the language consists of symmetric relations, between any pair of points in one of the $B_i$ and labeling it with the 2-type between these points (i.e. the isomorphism
type of that pair taken in order); a solution to the problem consists of determining the 2-types between the extension points in distinct $B_i$, which may then be placed on a dotted line between the points. 

Examples of amalgamation diagrams, both with and without solutions, may be found throughout, beginning with Proposition \ref{prop:genericLambdaUltrametric}.

It is worth noting that in order to verify that some class satisfies the amalgamation property, it suffices to verify a weaker form called 2-point amalgamation, in which each $B_i$ contains one extension point. Although this is true in general, the proof is particularly simple for the cases in this thesis since the following hold.
\begin{itemize}
\item The languages we consider are binary.
\item In the amalgamation strategies we consider, to determine the 2-type between two extension points, only those extension points and the base are used.
\end{itemize}
By the first point, a general amalgamation problem only requires determining the 2-types between extension points in separate factors. By the second point, each of these 2-types can be determined independently by solving the 2-point amalgamation problem containing the same base and the two relevant extension points. 

\begin{example}
Let $\GG$ be the class of all finite graphs. Then we may complete amalgamation diagrams by taking the disjoint union of the two factors over the base. The \fraisse limit of this class is the \textit{infinite random graph}.
\end{example}


\begin{example}
Let $\LL$ be the class of all finite linear orders. Then we may complete amalgamation diagrams by first adding any relation forced by transitivity, and then completing arbitrarily to a linear order. The \fraisse limit of this class is $(\Q, <)$.
\end{example}


\begin{example}
Let $\KK$ be the class of all finite $n$-dimensional permutation structures, for some $n$. Since linear orders can be amalgamated, and the amalgamation strategies can be carried out independently, $\KK$ is an amalgamation class. Thus it has a Fra\"\i ss\'e limit, called the \textit{fully generic} $n$-dimensional permutation structure.
\end{example}

\subsection{Strong Amalgamation and Related Topics}

\begin{definition}
An amalgamation class $\AA$ has \textit{strong amalgamation} if every amalgamation problem can be solved without identifying points from the two factors.
\end{definition}

All the examples from the Section \ref{section:diagrams} are strong amalgamation classes. A non-example would be the amalgamation class corresponding to a finite homogeneous structure.

We now introduce a model-theoretic concept that, by a standard result, gives an alternative characterization of strong amalgamation.

\begin{definition}
Let $M$ be a structure and $A \subset M$. Then $b \in M$ is \textit{algebraic} over $A$ if there is some first-order formula $\phi(x)$ with parameters from $A$ such that $\phi(M)$ is finite and $\phi(b)$ holds.

The \textit{algebraic closure of A in M}, denoted $acl(A)$, is the set of all $b \in M$ algebraic over $A$. 
\end{definition}

\begin{proposition}
Let $M$ be a homogeneous structure and $\AA$ its class of finite substructures. Then $\AA$ is a strong amalgamation class iff $acl(A) = A$ for all $A \subset M$.
\end{proposition}

Finally, we define the notion of a generic linear order.

\begin{definition}
Let $M$ be the \fraisse limit of a strong amalgamation class $\AA$. Then the class $\AA^<$ of all expansions of structures in $\AA$ by a linear order is also a strong amalgamation class. The \fraisse limit of $\AA^<$ is said to be the expansion of $M$ by a \textit{generic linear order}.
\end{definition}

\section{Structural Ramsey Theory} \label{section:RamseyBackground}
\begin{definition}
Let $\KK$ be a class of linearly ordered structures, closed under isomorphism. Given $A, B \in \KK$, let $\binom B A$ denote the set of substructures of $B$ that are isomorphic to $A$. We will say $\KK$ is a \textit{Ramsey class} if for any $n \in \N$ and $A, B \in \KK$, there is a $C \in \KK$ such that if $\binom C A$ is colored with $n$ colors, there is a $\widehat B \in \binom C B$ such that $\binom {\widehat B} A$ is monochromatic (we will often just say $\widehat B$ is monochromatic).

A homogeneous structure $\MM$ will be called a \textit{Ramsey structure} if its class of finite substructures is a Ramsey class.
\end{definition}

We give the above definition in terms of linearly ordered structures and coloring substructures. One could drop the assumption of a definable linear ordering and work with embeddings instead of substructures. However, it is fairly easy to see that the Ramsey property with embeddings implies the structures must be rigid. Furthermore, in the cases we will be considering, the Ramsey class will have a \fraisse limit $\Gamma$, and one can prove, via the connection with topological dynamics discussed in Section \ref{sect:dynamicsBackground}) that $\Gamma$ must have a $\emptyset$-definable linearly ordering. In particular, let $Aut(\Gamma)$ act on the compact space of linear orders of $\Gamma$ via the logic action; by extreme amenability, this action has a fixed point, which corresponds to a definable linear order.

The following theorem connects \fraisse theory and structural Ramsey theory.

\begin{theorem} [\cite{Nes}]
Let $\KK$ be a hereditary Ramsey class with the joint-embedding property. Then $\KK$ has the amalgamation property.
\end{theorem}

While this statement narrows down the search for Ramsey classes, it seems to go in only one direction. In order to find Ramsey structures, the most efficient method seems to be to start with a homogeneous structure and then add relations to obtain a Ramsey class. In the known cases, it is usually sufficient to add one or more linear orders, subject to the condition that every $\emptyset$-definable equivalence relation should be convex with respect to some order, which may be viewed as also considering elements from $M^{eq}$. 

Perhaps the main open question in the area is the following.
\begin{question}
 Does every homogeneous structure in a finite relational language admit an expansion to a Ramsey structure that is also homogeneous in a finite relational language? If so, describe that expansion.
 \end{question}
 
 An example of Evans, using a Hrushovski construction, shows the hypothesis of homogeneity cannot be relaxed to $\omega$-categoricity \cite{EHN2}.

\subsection{Topological Dynamics} \label{sect:dynamicsBackground}
A \textit{G-flow} is a continuous action of a topological group $G$ on a compact Hausdorff space. A \textit{G-flow} is minimal if it has no non-trivial subflows, or equivalently if every orbit is dense. By a theorem of Ellis, every topological group $G$ has a unique minimal flow that maps homomorphically onto any other minimal $G$-flow, and is called the \textit{universal minimal flow}. A topological group $G$ is called \textit{extremely amenable} if its universal minimal flow is a single point, or equivalently if every $G$-flow has a fixed point.

The following seminal theorem connects structural Ramsey theory and topological dynamics.

\begin{theorem}[\cite{KPT}]
Let $M$ be a homogeneous structure and $G = Aut(M)$, equipped with the pointwise convergence topology. Then $G$ is extremely amenable iff $M$ is a Ramsey structure.
\end{theorem}

This is an extreme case of the phenomenon that the universal minimal flow of the automorphism group of a homogeneous structure encodes its obstructions to Ramseyness. When the structure is not already Ramsey, finding a suitable Ramsey expansion is equivalent to describing the universal minimal flow of its automorphism group. This correspondence is laid out more fully in later in this section.

More broadly than the case of extreme amenability, it is of interest when the universal minimal flow is in some sense small. The condition that has received the most attention is metrizability, which is equivalent to second-countability as the universal minimal flow is compact Hausdorff. 

As discussed below, given a \textit{precompact} Ramsey expansion with the \textit{expansion property} of a homogeneous structure $\MM$, one can describe the universal minimal flow of $Aut(\MM)$, and it will be metrizable. In fact, for any homogeneous structure $\MM$, if the universal minimal flow of $Aut(\MM)$ is metrizable, it can always be identified in this fashion \cite{Zucker}.

We now give the details of identifying the universal minimal flow of the automorphism group of a homogeneous structre. The following definitions and theorem are from \cite{NVT}, extending the work of \cite{KPT}.

\begin{definition}
Given an $L$-structure $F$, we let $\text{Age}(F)$ denote the set of all $L$-structures isomorphic to a finite substructure of $F$.
\end{definition}

\begin{definition}
Let $L$ be a relational language, and $L^*$ a countable relational expansion of $L$. Let $F$ be a homogeneous $L$-structure. Then an $L^*$-expansion $F^*$ of $F$ is \textit{precompact} if any $A \in \text{Age}(F)$ has only finitely many $L^*$-expansions in $\text{Age}(F^*)$.
\end{definition}

\begin{definition}
Let $F$ be a homogeneous structure, and $F^*$ a precompact relational expansion of $F$. Then $\text{Age}(F^*)$ has the \textit{expansion property} relative to $\text{Age}(F)$ if for every $A \in \text{Age}(F)$ there exists a $B \in \text{Age}(F)$ such that, for any $L^*$-expansions $A^*$ of $A$ and $B^*$ of $B$ in $\text{Age}(F^*)$, $A^*$ embeds into $B^*$.
\end{definition}

\begin{definition}
Let $L$ be a language, and $L^* = L \cup \set{R_i}_{i \in I}$ be a relational expansion. Let $a(i)$ be the arity of $R_i$. Given a homogeneous $L$-structure $F$, we define $P^*$ as
$$P^* = \prod_{i \in I}\set{0, 1}^{F^{a(i)}}$$

We may define a group action of $\text{Aut}(F)$ on a given factor as follows: for $g \in \text{Aut}(F)$ and $S_i \in F^{a(i)}$, $g \cdot S_i(y_1, ..., y_{(a_i)}) \Leftrightarrow S_i(g^{-1}(y), ..., g^{-1}(y_{a(i)}))$.

Finally, we may define the \textit{logic action} of $\text{Aut}(F)$ on $P^*$ as given by applying the above action componentwise.
\end{definition}

\begin{theorem} [\cite{NVT}*{Theorems 4, 5}] \label{theorem:NVT}
Let $L$ be a language, $L^* = L \cup \set{R_i}_{i \in I}$ be a countable relational expansion, and $F$ a homogeneous $L$-structure. Let $F^* = (F, \vec R^*)$ be a precompact $L^*$-expansion of $F$. Then we have the following equivalence, with the closure taken in $P^*$.
\begin{enumerate}
\item The logic action of $Aut(F)$ on $\overline{Aut(F) \cdot \vec R^*}$ is minimal.
\item $\text{Age}(F^*)$ has the expansion property relative to $\text{Age}(F)$.
\end{enumerate}

Furthermore, the following are also equivalent.

\begin{enumerate}
\item The logic action of $Aut(F)$ on $\overline{Aut(F) \cdot \vec R^*}$ is its universal minimal flow.
\item $\text{Age}(F^*)$ has the Ramsey property and the expansion property relative to $\text{Age}(F)$.
\end{enumerate}
\end{theorem}


\chapter{Homogeneous Permutations, $\Lambda$-Ultrametric Spaces, and Subquotient Orders}\label{chap:lambda}
\section{Introduction}
This chapter contains preparatory material for the main construction of homogeneous finite-dimensional permutation structures in the next chapter. First, we present Cameron's classification of the homogeneous permutations, where the ideas of the construction and the corresponding catalog are present in a simple form.

We then present the notions of \textit{$\Lambda$-ultrametric spaces} and \textit{subquotient orders}. $\Lambda$-ultrametric spaces give a convenient a presentation of structures in a language of equivalence relations such that transitivity becomes analogous to the metric triangle inequality. Although they are more general, we will primarily use subquotient orders to provide an alternative means of presenting linear orders convex with respect to various equivalence relations. Such linear orders may interact with each other in complex ways in our construction, but by interdefinably moving to subquotient orders these interactions vanish and everything becomes generic.   

\section{Homogeneous Permutations}
Up to interdefinability, there are 3 homogeneous permutations. We divide them using the following notion.

\begin{definition}
A structure is \textit{primitive} if it has no non-trivial $\emptyset$-definable equivalence relation.
\end{definition}

\begin{theorem}[\cite{Cameron}]
Every homogeneous permutation is interdefinable with one of the following.
\begin{itemize}
\item Imprimitive: Lexicographic ordering (more naturally viewed as equipped with a single order and an equivalence relation with convex classes, as explained below)
\item Primitive:
\begin{itemize}
\item Degenerate: The orders agree, up to reversal
\item Fully generic, in the sense of the Fra\"\i ss\'e theory of amalgamation classes
\end{itemize}
\end{itemize}
\end{theorem}

The imprimitive structure is naturally presented as $(\Q^2, E, <_{lex})$, where $E$ is the equivalence relation for agreement in the first coordinate and $<_{lex}$ is the standard lexicographic order. Note that $E$-classes are $<_{lex}$-convex. This structure can be presented as a permutation by interdefinably replacing $E$ by another order $<^*$, which agrees with $<_{lex}$ within $E$-classes and disagrees with $<_{lex}$ between $E$-classes. 

This construction will be generalized in order to generate our catalog of homogeneous finite-dimensional permutation structures. We will start with a homogeneous structure consisting solely of equivalence relations. This will then be expanded by enough orders so that the classes of each equivalence relation are convex with respect to some $\emptyset$-definable order. Finally, we may interdefinably replace the equivalence relations with additional orders.

\section{$\Lambda$-Ultrametric Spaces}
\setcounter{figure}{0}

In this section, we set up a language that will be more convenient for our amalgamation arguments than the language of equivalence relations. One point that is somewhat ungainly from the point of view of equivalence relations is that when passing to a substructure, for example a single point, various equivalence relations can collapse, but we would like to keep the lattice $\Lambda$ of equivalence relations fixed.

\begin{definition} Let $\Lambda$ be a lattice. A \textit{$\Lambda$-ultrametric space} is a metric space where the metric takes values in $\Lambda$ and the triangle inequality uses the join rather than addition. Analogous to a pseudometric space, we also define a \textit{$\Lambda$-ultrapseudometric space} as a $\Lambda$-ultrametric space without the requirement that the metric assign non-zero distance to distinct points.
\end{definition}

\begin{example}
If $\Lambda$ is a chain, then $\Lambda$-ultrametric spaces are ultrametric spaces in the standard sense.
\end{example}

These may be viewed in a relational language by using a binary relation for each possible distance.

As in metric spaces, quotienting out by the relation $d(x,y) = 0$ in a $\Lambda$-ultrapseudometric space yields a $\Lambda$-ultrametric space, and we also have a path variant of the the triangle inequality: $d(x,y)$ is no greater than the join of the distances between points on any path from $x$ to $y$.

\begin{theorem}\label{theorem:isomorphism}
Fix a finite lattice $\Lambda$. Let $\MM_\Lambda$ be the category of $\Lambda$-ultrametric spaces, with isometries as morphisms. Let $\EE\QQ_\Lambda$ be the category of structures consisting of a set equipped with a family of not-necessarily-distinct equivalence relations $\set{E_\lambda| \lambda \in \Lambda}$ satisfying the following conditions, with embeddings as morphisms.
\begin{enumerate}
\item $\set{E_\lambda}$ forms a lattice.
\item The map $L: \lambda \mapsto E_\lambda$ is meet-preserving. In particular, if $\lambda_1 \leq \lambda_2$, then $E_{\lambda_1} \leq E_{\lambda_2}$.
\item $E_\bbzero$ is equality and $E_\bbone$ is the trivial relation. 
\end{enumerate}

Then $\EE\QQ_\Lambda$ is isomorphic to $\MM_\Lambda$. Furthermore, the functors of this isomorphism preserve homogeneity.

\end{theorem}
\begin{proof}
We first define the functors $m: \EE\QQ_\Lambda \to \MM_\Lambda$ and $e: \MM_\Lambda \to \EE\QQ_\Lambda$ giving this isomorphism.

Given a $A \in \EE\QQ_\Lambda$, we define $m(A)$ by taking the same universe and defining $d(x, y) = \bigwedge \set{\lambda \in \Lambda | x E_\lambda y}$. In the reverse direction, given $M \in \MM_\Lambda$, we get define $e(M)$ by taking the same universe and defining $E_\lambda = \set{(x,y) | d(x,y) \leq \lambda}$.

We must check that the images of these maps lie in the specified codomains, that they are inverses of each other, and that morphisms are preserved.

\begin{claim}
Let $A \in \EE\QQ_\Lambda$. Then $m(A) \in \MM_\Lambda$.
\end{claim} 
\begin{claimproof}
Symmetry is clear. As the finest equivalence relation between $x$ and itself is equality, and we have assumed $E_\bbzero$ is equality, we have $d(x,x) = \bbzero$. Conversely, if $d(x,y) = \bbzero$, then the finest equivalence relation holding between $x$ and $y$ is equality.

We now check the triangle inequality. Fix $x, y, z \in A$. Note $E_{d(x, y)}$ is the finest equivalence relation holding between $x, y$, and similarly $E_{d(x,z)}$ and $E_{d(z,y)}$ for $x,z$ and $z,y$, respectively. Transitivity of equivalence relations and then the fact that $L$ is order preserving imply  $E_{d(x, y)} \leq E_{d(x, z)} \join E_{d(z, y)} \leq E_{d(x,z) \join d(z,y)}$. Thus $d(x, y) = \bigwedge \set{\lambda | xE_\lambda y} \leq d(x,z) \join d(z, y)$.
\end{claimproof}

\begin{claim}
Let $M \in \MM_\Lambda$. Then $e(M) \in \EE\QQ_\Lambda$.
\end{claim}
\begin{claimproof}
We first check the $E_\lambda$ are equivalence relations. For any $x \in M$, $d(x,x) = \bbzero$, and so in $e(M)$, $xE_\lambda x$ for every $\lambda$. Symmetry is clear.

For transitivity, suppose $xE_\lambda y E_\lambda z$. Then $d(x, y), d(y, z) \leq \lambda$. Thus $d(x,z) \leq \lambda \join \lambda = \lambda$, and so $x E_\lambda z$.

For preservation of meets, note $E_{\lambda \meet \lambda'}(x, y) \iff d(x, y) \leq \lambda \meet \lambda' \iff d(x, y) \leq \lambda, \lambda' \iff E_\lambda(x, y), E_{\lambda'}(x, y)$.

We have $E_\bbzero(x,y) \iff d(x, y) \leq \bbzero \iff d(x,y) = \bbzero \iff x=y$. Also, $E_\bbone(x, y) \iff d(x, y) \leq \bbone$, which is true for all $x,y$.

As $\set{E_\lambda}$ forms a finite bounded meet-semilattice, it must form a lattice.
\end{claimproof}

\begin{claim}
Let $A \in \EE\QQ_\Lambda$ and $M \in \MM_\Lambda$. Then $e(m(A)) = A$ and $m(e(M)) = M$.
\end{claim}
\begin{claimproof}
In $m(A)$, $d(x, y)$ is the least $L$-preimage of the finest equivalence relation to hold between $x$ and $y$. In $e(m(A))$, we have $E_\lambda(x,y) \iff \lambda \geq d(x,y)$, which thus gives $A$.

In $e(M)$, we have $E_{\lambda(x,y)} \iff \lambda \geq d(x,y)$. In $m(e(M))$, we have $d(x,y)$ is the least $L$-preimage of the finest equivalence relation to hold between $x$ and $y$, which will be $\lambda$.
\end{claimproof}

\begin{claim}
Let $M_1, M_2 \in M_\Lambda$, and $A_1 = e(M_1), A_2 = e(M_2)$. Then $f: M_1 \to M_2$ is an isometry iff $f$ induces an embedding from $A_1$ to $A_2$.
\end{claim}
\begin{claimproof}
Suppose $d(x,y) = d(f(x), f(y))$, and suppose $E_\lambda(x,y)$ in $A_1$. Then $\lambda \geq d(x,y) = d(f(x), f(y))$. Thus $E_\lambda(f(x),f(y))$.

Now suppose $E_\lambda(x,y) \iff E_\lambda(f(x),f(y))$. As $d(a,b)$ is determined by the $\lambda$ such that $E_\lambda(a,b)$, we get $d(x,y) = d(f(x),f(y))$.
\end{claimproof}

Finally, to see the functors preserve homogeneity, note they preserve the notion of (partial) isomorphism.
\end{proof}

There is a well known amalgamation strategy for metric spaces called \textit{shortest path completion}, which determines the distance between extension points in distinct factors by taking the length of the shortest path between them. We here give the analog for $\Lambda$-ultrametric spaces.

\begin{definition}
Consider an amalgamation diagram of $\Lambda$-ultrametric spaces with base $B$. Let $x$ and $y$ be extension points in different factors, and for each $b_i \in B$ let $d(x, b_i) = e_i$ and $d(y, b_i) = e'_i$. \textit{Pre-canonical amalgamation} is the amalgamation strategy assigning $d(x, y) = \bigwedge_i(e_i \vee e'_i)$. 
\end{definition}

\begin{definition}
Consider an amalgamation diagram of $\Lambda$-ultrametric spaces. Suppose that pre-canonical amalgamation assigns $d(x, y) = \bbzero$ only if $x$ and $y$ realize the same 1-type over the base (in particular, if pre-canonical amalgamation yields a $\Lambda$-ultrapseudometric space). Then we may define \textit{canonical amalgamation} as the strategy of pre-canonical amalgamation, followed by identifying $x$ and $y$ if $d(x, y) = \bbzero$. 
\end{definition}

Two-point pre-canonical amalgamation is shown in the figure below (for guidance on the interpretation of amalgamation diagrams in this thesis, see Section \ref{section:diagrams}). Note that by the triangle inequality, we must have $d(x, y) \leq e_i \vee e'_i$ for each $i$. Thus pre-canonical amalgamation attempts to make $d(x, y)$ maximal while respecting these instances of the triangle inequality. The next proposition provides a condition on when this is sufficient to ensure the resulting amalgam satisfies the triangle inequality, with the proof a straightforward adaptation of the argument for metric spaces. 

\begin{figure}[h]
\begin{diagram}
x \odot & & \rDash^{\bigwedge_i(e_i \vee e'_i)} & & \odot y \\
& \luLine_{\set{e_i}_{i \leq |B|}} &  & \ruLine_{\set{e'_i}_{i \leq |B|}} & \\
& & \underset{B}{\bigcirc} & & \\
\end{diagram}
\caption{}
\label{fig:1}
\end{figure}
\goodbreak

\begin{proposition} \label{prop:genericLambdaUltrametric}
Let $\Lambda$ be a distributive lattice, and let $\KK$ be the class of all finite $\Lambda$-ultrametric spaces. Then $\KK$ is an amalgamation class, and any amalgamation diagram can be completed by canonical amalgamation.
 \end{proposition}

\begin{proof} 
It suffices to check that pre-canonical amalgamation in a 2-point amalgamation diagram produces a $\Lambda$-ultrapseudometric space. In other words, we check that Figure \ref{fig:1} satisfies the triangle inequality, given the distributivity of $\Lambda$. 

Fix some $b_i \in B$ and consider the corresponding triangle.

We have $d(x, y) =  \bigwedge_j (e_j \vee e'_j) \leq (e_i \vee e_i') = d(x, b_i) \vee d(b_i, y)$, so this side satisfies the triangle inequality by definition. 

The remaining two sides are handled symmetrically, so we only check $d(b_i, y) \leq d(b_i, x) \vee d(x, y)$. We have $d(b_i, x) \vee d(x, y) = e_i \vee (\bigwedge_j (e_j \vee e'_j)) = \bigwedge_j (e_i \vee e_j \vee e'_j)$, by distributivity. We now use the path variant of the triangle inequality, which the diagram satisfies before the distance between $x$ and $y$ is determined. Going from $b_i$ to $y$, for each $j$ we get $e_i \vee e_j \vee e'_j \geq e_i'$ (see Figure \ref{fig:2}), giving $d(b_i, x) \vee d(x, y) \geq \bigwedge_j e'_i =  d(b_i, y)$.

\begin{figure}[h]
\begin{diagram}
& & \overset{b_j}{\bullet} & & \\
& \ldLine^{e_j} &  & \rdLine^{e'_j} & \\
x \odot & & & & \odot y \\
& \luLine_{e_i} &  & \ruLine_{e'_i} & \\
& & \underset{b_i}{\bullet} & & \\
\end{diagram}
\caption{}
\label{fig:2}
\end{figure}
\end{proof}

\section{Subquotient Orders}
Rather than expanding $\Lambda$-ultrametric spaces directly by linear orders, there are technical advantages to working with certain partial orders on substructures of quotients, which we call \textit{subquotient orders}. The idea behind the change of language is that when a linear order on $\Gamma$ is convex with respect to an equivalence relation $E$, it is better viewed as two partial orders: one that orders points within any given $E$-class, and one that encodes an order on $\Gamma/E$. Working with subquotient orders increases the variety of structures produced by our construction (see Examples \ref{ex:complex3} and \ref{ex:fullQ2}, and the following remark), and significantly simplifies certain proofs.

\begin{definition}
Let $X$ be a structure, and $E \leq F$ equivalence relations on $X$. A \textit{subquotient-order from $E$ to $F$} is a partial order on $X/E$ in which two $E$-classes are comparable iff they lie in the same $F$-class (note, this pulls back to a partial order on $X$). Thus, this partial order provides a linear order of $C/E$ for each $C \in X/F$. We call $E$ the \textit{bottom relation} and $F$ the \textit{top relation} of the subquotient-order. 
\end{definition}

Depending on the context, we will switch between considering a given subquotient order as a partial order on equivalence classes, or its pullback to a partial order on points. A special case of this is when the subquotient order has equality as its bottom relation, which amounts to equating $X$ with $X/=$.

The following theorem is the main result of the next section, but as we use it in this section, we state it now.

\begin{definition}
Let $\Lambda$ be a lattice, and $x \in \Lambda$. Then $x \in \Lambda$ is \textit{meet-irreducible} if there do not exist $y,z \neq x$ such that $x = y \meet z$.
\end{definition}

\begin{theorem} 
Let $\Lambda$ be a finite distributive lattice, and $\Gamma$ the generic $\Lambda$-ultrametric space. For each meet-irreducible $E \in \Lambda$, fix a function $f_E: \set{F \in \Lambda| E<F} \to \N$. Then there is a homogeneous expansion of $\Gamma$, generic in a natural sense, adding, for each meet-irreducible $E \in \Lambda$ and $F>E$, $f_E(F)$ subquotient orders $\set{<_{E,F,i}}_{i=1}^{f_E(F)}$ from $E$ to $F$, generic in the following sense. 
\end{theorem}

\begin{remark}
The meaning of ``generic in a natural sense'' is made more precise in Theorem \ref{theorem:amalg}. It is the natural analogue of an expansion by a generic linear order.
\end{remark}

We now define two useful constructions with subquotient orders, and then give two examples of homogeneous finite-dimensional permutation structures that cannot be produced by the construction of \cite{Lattice}, but which can be produced using the theorem above.

For the remainder of this section, if $x$ is an $E$-class, and $F$ an equivalence relation above $E$, then $x/F$ will represent the $F$-class containing $x$.

\begin{definition}
Let $<_{E, F}$ be a subquotient order with bottom relation $E$ and top relation $F$, and let  $<_{F, G}$ be a subquotient order with bottom relation $F$ and top relation $G$. Then the \textit{composition} of $<_{F, G}$ with $<_{E, F}$, denoted $<_{F, G}[<_{E, F}]$, is the subquotient order with bottom relation $E$ and top relation $F$ given by $x <_{F, G}[<_{E, F}] y$ iff either of the following holds.
\begin{enumerate}
\item $x$ and $y$ are in the same $F$-class, and $x <_{E, F} y$
\item $x$ and $y$ are in distinct $F$-classes, and $x/F <_{F, G} y/F$.
\end{enumerate}
\end{definition}

\begin{definition}
Let $<_{E, F}$ be a subquotient order with bottom relation $E$ and top relation $F$, and let $G$ be an equivalence relation lying between $E$ and $F$. Then the \textit{restriction of $<_{E, F}$ to $G$}, denoted $<_{E, F} \upharpoonright_G$, is the subquotient order with bottom relation $E$ and top relation $G$ given by $x <_{E, F} \upharpoonright_G y$ iff $x$ and $y$ are in the same $G$-class and $x <_{E, F} y$.
\end{definition}

\begin{example} \label{ex:complex3}
Let $\AA$ be the amalgamation class consisting of all finite structures in the language $\set{E, <_1, <_2}$, where $E$ is an equivalence relation, $<_1$ is a linear order, and $<_2$ is an $E$-convex linear order that agrees with $<_1$ on $E$-classes. 

 Let $\AA'$ be the class of all finite structures in the language $\set{E', <'_1, <'_2}$, where $E'$ is an equivalence relation, $<'_1$ is a subquotient order from $=$ to $\bbone$, and $<'_2$ a subquotient order from $E'$ to $\bbone$. This is also an amalgamation class, and its Fra\"\i ss\'e limit $\Gamma'$ is interdefinable with the Fra\"\i ss\'e limit $\Gamma$ of $\AA$. 

To define $\Gamma$ from $\Gamma'$, let $<_1 = <'_1$, and let $<_2 = <'_2[<'_1 \upharpoonright_E]$. To define $\Gamma'$ from $\Gamma$, let $<'_1 = <_1$, and let $x <'_2 y$ iff $\neg xEy$ and $x <_2 y$.
\end{example}

For the next example, we will use the following lemma, which also enters into the proof of Lemma \ref{lemma:sqoproduct}.

\begin{lemma} \label{lemma:movesqo}
Let $\Gamma$ be the generic $\Lambda$-ultrametric space. Let $E \in \Lambda$, with $E = F_1 \meet F_2$. Then a subquotient order $<_{F_1, F_1 \join F_2}$ on $\Gamma$ with bottom relation $F_1$ and top relation $F_1 \join F_2$ induces a definable subquotient order with bottom relation $E$ and top relation $F_2$.
\end{lemma}
\begin{proof}
We wish to define an order on $E$-classes within $F_2$-classes. Since $E = F_1 \meet F_2$, within a given $F_2$-class each $E$-class is in a distinct $F_1$-class, and they are all in the same $(F_1 \join F_2)$-class. Thus, we can define a subquotient order $<_{E, F_2}$ with bottom relation $E$ and top relation $F_2$ by $x <_{E, F_2} y \Leftrightarrow x/F_1 <_{F_1, F_1 \join F_2} y/F_1$.
\end{proof}

\begin{example} \label{ex:fullQ2}
For a more complex example of the use of subquotient orders, consider the full product $\Q^2$. This is a homogeneous structure with universe $\Q^2$ in the language $\set{E_1, E_2, <_1, <_2}$, where $E_1$ and $E_2$ are the relations defined by agreement in the first and second coordinates, respectively, $<_1$ is a generic subquotient order from $E_1$ to $\bbone$, and $<_2$ is a generic subquotient order from $E_2$ to $\bbone$.

Since $E_1 \meet E_2 = \bbzero$, we see that $<_1$ defines a linear ordering on each $E_2$-class, and $<_2$ defines a linear ordering on each $E_1$-class. Thus, the composition (abusing notation slightly) $<_1[<_2]$ defines an $E_1$-convex linear order, and $<_2[<_1]$ defines an $E_2$-convex linear order. 
\end{example}

\begin{remark}
These examples cannot be produced by our construction when using linear orders rather than subquotient orders, as we only forbid substructures of order 3. However, in Example \ref{ex:complex3}, we must forbid a substructure of order 2 to force $<_1$ and $<_2$ to agree between $E$-related points. In Example \ref{ex:fullQ2}, we must forbid the following substructure of order 4 (as well as another symmetric substructure):
\begin{enumerate}
\item $x_1 E_1 x_2$, $y_1 E_1 y_2$, $\neg x_1 E_1 y_1$
\item $x_1 E_2 y_1$, $x_2 E_2 y_2$, $\neg x_1 E_2 x_2$ 
\item $x_1 <_1 x_2$, $y_2 <_1 y_1$
\end{enumerate}
\end{remark}

If we try to generalize our construction using linear orders by allowing side conditions of the type appearing in the examples, it becomes unclear which such combinations should be allowed. We instead translate our construction to work with subquotient orders and then later move back to the language of linear orders. For this last translation step, the following definition is essential.

\begin{definition}
Let $\Lambda$ be a finite distributive lattice, and let $L$ be a language consisting of relations for the distances in $\Lambda$ and finitely many subquotient orders, labeled with their top and bottom relations. We say that the language $L$ is \textit{$\Lambda$-well-equipped} if for each $E \in \Lambda$, $E$ appears as the bottom relation of some subquotient order in $L$ with distinct bottom and top relations iff $E$ is meet-irreducible.

If $\AA_\Lambda$ is the class of all finite $\Lambda$-ultrametric spaces, and $L$ a $\Lambda$-well-equipped language, we will call $\vec \AA_\Lambda$ a \textit{well-equipped lift} of $\AA_\Lambda$ if it consists of all finite $\Lambda$-ultrametric spaces equipped with subquotient orders from $L$.
\end{definition}

\begin{lemma} \label{lemma:sqoproduct}
Let $\Lambda$ be a finite distributive lattice, $\AA_\Lambda$ the class of all finite $\Lambda$-ultrametric spaces, and $\vec \AA_\Lambda$ a well-equipped lift of $\AA_\Lambda$, with Fra\"\i ss\'e limit $\vec \Gamma$. Then for every $E<F \in \Lambda$, $\vec \Gamma$ has a definable subquotient order with bottom-relation $E$ and top-relation $F$.
\end{lemma}

\begin{proof}
We prove this by downward induction in $\Lambda$. Take an arbitrary $E \in \Lambda$, and assume the claim is true for every element above $E$. 

We first note that it is sufficient, for every $F' \in \Lambda$ covering $E$, to construct a definable subquotient order $<_{E, F'}$ with bottom relation $E$ and top relation $F'$. Indeed, by the induction hypothesis, there is some definable subquotient order $<_{F', F}$ with bottom relation $F'$ and top relation $F$, and then the composition $<_{F', F}[<_{E, F'}]$ gives the desired definable subquotient order.  

First assume $E$ is meet-irreducible. Then there is a unique $F' \in \Lambda$ covering $E$. By assumption, there is some subquotient order $<_E$ with bottom relation $E$ and top relation some $F'' \geq F'$. Then the restriction $<_{E}\upharpoonright_{F'}$ is as desired.

Now assume that $E$ is meet-reducible, and let $F'$ cover $E$. Since $E$ is meet-reducible, there is some $F''> E$ such that $E = F' \meet F''$. By the induction hypothesis, there is a definable subquotient order $<_{F'', F' \join F''}$ with bottom relation $F''$, and top relation $F' \join F''$. Then Lemma \ref{lemma:movesqo} provides a definable subquotient order with bottom relation $E$ and top relation $F'$.
\end{proof}

\begin{corollary} \label{corollary:extendsqo}
Let $\Lambda$ be a finite distributive lattice, $\AA_\Lambda$ the class of all finite $\Lambda$-ultrametric spaces, and $\vec \AA_\Lambda$ a well-equipped lift of $\AA_\Lambda$, with Fra\"\i ss\'e limit $\vec \Gamma$. Given any subquotient order $<_{E}$ on $\vec \Gamma$ with bottom relation $E$ and top relation $F$, we can define on $\vec \Gamma$ a subquotient order $<'_{E}$ with bottom relation $E$ and top relation $\bbone$, in such a way that $x <_{E} y$ iff $x <'_{E} y$ and $x, y$ are in the same $F$-class. 
\end{corollary}
\begin{proof}
By Lemma \ref{lemma:sqoproduct}, there is a definable subquotient order $<_{F}$ with bottom relation $F$ and top relation $\bbone$. Then the composition ${<_{F}[<_{E}]}$ is as desired.
\end{proof}

\begin{remark} \label{rem:choice}
We will later find it useful to have made concrete choices when applying Lemma \ref{lemma:sqoproduct} and Corollary \ref{corollary:extendsqo}. In particular, given an enumeration ${(<_{G, i})}$ of the subquotient orders with bottom relation $G$ for every $G \in \Lambda$, we may always use subquotient orders that have 1 in the second index, with the possible exception of the specified subquotient order $<_E$ in Corollary \ref{corollary:extendsqo}. 
\end{remark}

\begin{proposition} \label{prop:sqolininter}
Let $\Lambda$ be a finite distributive lattice, $\AA_\Lambda$ be the class of all finite $\Lambda$-ultrametric spaces, and $\vec \AA_\Lambda$ a well-equipped lift of $\AA_\Lambda$, with Fra\"\i ss\'e limit $\vec \Gamma$. Then the relations of $\vec \Gamma$ are interdefinable with a set of linear orders.
\end{proposition}

\begin{proof}
For each $E \in \Lambda$, and each subquotient order $<_{E, i}$ in the language with bottom relation $E$, let $<'_{E, i}$ be a subquotient order as in Corollary \ref{corollary:extendsqo}. By Lemma \ref{lemma:sqoproduct}, let $<_{\bbzero, E}$ be a definable subquotient order with bottom relation equality and top relation $E$, and let $<''_{E, i}$ be the linear order given by the composition $<'_{E, i}[<_{\bbzero, E}]$.

Then, in the language consisting of the equivalence relations $E \in \Lambda$, the set of subquotient orders is interdefinable with the set of corresponding linear orders produced above. Note that $<''_{E, i}$ is $E$-convex. Then, each $E \in \Lambda$ can be interdefinably replaced with a linear order $<^*_E$ as defined below.

For each $E \in \Lambda$:
\begin{enumerate}
\item Let $<_E$ be the definable linear order such that $E$ is $<_E$-convex
\item Let $<^*_E$ agree with $<_E$ within $E$-classes, and agree with the reverse of $<_E$ between $E$-classes.
\end{enumerate} 
\end{proof}


\chapter{Generating a Catalog of Homogeneous Finite-Dimensional Permutation Structures}\label{chap:construction}
\section{Introduction}
\setcounter{figure}{0}
In this chapter, we present the amalgamation argument needed for our catalog, which consists of all \fraisse limits of well-equipped lifts of generic $\Lambda$-ultrametric spaces. By Proposition \ref{prop:sqolininter}, these structures are quantifier-free interdefinable with homogeneous finite-dimensional permutation structures. 

It is an interesting problem to determine the minimum number of orders required to represent a given structure in our catalog, and we address this point in section \ref{subsection:OrderNumber}.

We also note that we have a positive answer to Question \ref{q:ordexp} in the case of the structures we produce, that we may realize any finite distributive lattice $\Lambda$ as the lattice of all $\emptyset$-definable equivalence relations in some homogeneous finite-dimensional permutation structure, and that Conjecture \ref{conj:primitivity}, the Primitivity Conjecture, is truly a special case of Conjecture \ref{conj:sqoconj}.

\section{The Main Construction}

We first repeat the theorem we wish to prove. Note that this will cover some structures not in our catalog, since we don't require every meet-irreducible be the bottom relation for some subquotient order.

\begin{theorem} \label{theorem:expandLambda}
Let $\Lambda$ be a finite distributive lattice, and $\Gamma$ the generic $\Lambda$-ultrametric space. For each meet-irreducible $E \in \Lambda$, fix a function $f_E: \set{F \in \Lambda| E<F} \to \N$. Then there is a homogeneous expansion of $\Gamma$, which is generic in a natural sense, adding, for each meet-irreducible $E \in \Lambda$ and $F>E$, $f_E(F)$ subquotient orders from $E$ to $F$.
\end{theorem}

In particular, if $f_E$ is identically $0$ for every $E$, we recover $\Gamma$ itself, and we produce a well-equipped lift if $f_E$ is non-zero for every $E$.

After the following preparatory definition, we restate this theorem in terms of the amalgamation property, and prove it.

\begin{definition}
Let $X$ be a structure, equipped with a binary relation $R$ and an equivalence relation $E$. We say that $E$ is a \textit{$R$-congruence} if $E(x, x')$ and $E(y, y')$ implies that $R(x, y)$ iff $R(x', y')$.
\end{definition}

\begin{theorem} \label{theorem:amalg}
Let $\Lambda$ be a finite distributive lattice. Let $\AA^*$ be the class of finite structures $(A,d,\set{<_{E_i}}_{i=1}^n)$ satisfying the following conditions.
\begin{itemize}
\item $(A,d)$ is  a $\Lambda$-ultrametric space.
\item $<_{E_i}$ is a subquotient order with bottom relation $E_i$, for some meet-irreducible $E_i \in \Lambda$, and top relation $F_i \in \Lambda$.
\end{itemize}
Then $\AA^*$ is an amalgamation class.
\end{theorem}

\begin{proof}[Proof of Theorem \ref{theorem:amalg}]
\setcounter{claim}{0}
Like linear orders, subquotient orders may be amalgamated independently, so we may assume $n=1$, and we will call the only subquotient order $<_E$ with bottom relation $E$ and top relation $F$.

We first introduce some notation. We define the relations $\preceq_E$ and $\arr$ on $\AA^*$-structures by 
\begin{enumerate}
\item $a \preceq_E b \Leftrightarrow (d(a, b) \leq E) \vee (a <_E b)$
\item $a \arr b \Leftrightarrow \exists x(a \preceq_E x \preceq_E b) \wedge (d(a, b) \not\leq E)$.
\end{enumerate}

We will make use of the following properties of $\preceq_E$ on $A^*$-structures.
\begin{enumerate}
\item If $a \preceq_E b <_E c$ or $a <_E b \preceq_E c$, then $a <_E c$.
\item $\preceq_E$ is transitive.
\item If $a \preceq_E b \preceq_E c$ and $d(a, c) \leq E$, then $d(a, b), d(a, c) \leq E$.
\item If $a \preceq_E b \preceq_E a$, then $d(a, b) \leq E$.   
\end{enumerate}

Property $(1)$ follows from the fact that $E$ is a $\preceq_E$-congruence. Properties $(2)$ and $(3)$ follow from $(1)$, and $(4)$ is a special case of $(3)$.

It suffices to show that $\AA^*$ contains solutions to all two-point amalgamation problems
$A_0^*\includedin A_1^*,A_2^*$, $A_i^*=A_0^*\union \{a_i\}$ for $i=1,2$.

Let $A$ be the extension of the free amalgam given by determining $d(a_1, a_2)$ by pre-canonical amalgamation. Either $<_E$ is already a subquotient order with bottom relation $E$ and top relation $F$, or we need to extend it to one by determining either $a_1 <_E a_2$ or $a_2 <_E a_1$. We break this into three cases.

\begin{claim}
Suppose $d(a_1, a_2) \leq E$. Then for $x \in A^*_0$, we have
$$a_1 <_E x \Longleftrightarrow a_2 <_E x$$
In particular, $<_E$ is a subquotient order on $A$ from $E$ to $F$.
\end{claim}
\begin{claimproof}
Since $E$ is meet-irreducible, if pre-canonical amalgamation yields $d(a_1, a_2) \leq E$, then there is a $y \in A_0^*$ such that $d(a_1, y), d(a_2, y) \leq E$.

By the fact that $E$ is a $<_E$-congruence, we get $a_1 <_E x \Longleftrightarrow y <_E x \Longleftrightarrow a_2 <_E x$. This proves the first part of the claim, and the second part follows immediately.
\end{claimproof}

We also note that if $d(a_1, a_2) = \bbzero$, then by the claim above $A_1 \cong A_2$, so we may amalgamate by identifying $a_1$ with $a_2$.

\begin{claim}
Suppose $d(a_1, a_2) \not\leq F$. Then $<_E$ is a subquotient order on $A$ from $E$ to $F$.
\end{claim}
\begin{claimproof}
This is clear, as $a_1$ and $a_2$ lie in distinct $F$-classes in $A$.
\end{claimproof}

\begin{claim}
Suppose $d(a_1, a_2) \in (E, F]$. On $A$, define $<^*_E = <_E \cup \arr$. Then
\begin{enumerate}
\item $a_1 \arr a_2$ and $a_2 \arr a_1$ cannot both hold.
\item $E$ is a $<^*_E$-congruence.
\end{enumerate}
\end{claim}
\begin{claimproof}
\hspace{1 cm} \newline
$(1)$ Suppose $a_1 \arr a_2 \arr a_1$. Then there exist $x_1, x_2$ such that $a_1 \preceq_E x_1 \preceq_E a_2$, and $a_2 \preceq_E x_2 \preceq_E a_1$.

In particular, $x_1 \preceq_E x_2 \preceq_E x_1$, so $d(x_1, x_2) \leq E$. As $d(a_1, a_2) \not \leq E$, we may suppose $d(a_1, x_2) \not \leq E$. 

But $x_2 \preceq_E a_1$, so $x_2 <_E a_1 \preceq_E x_1$. Thus $x_2 <_E x_1$, which contradicts $x_2 \preceq_E x_1$.

$(2)$ We check that $E$ is a $<^*_E$-congruence. Since $d(a_1, a_2) \not \leq E$, it suffices without loss of generality to consider some $x \in A^*_0$ such that $d(a_1, x) \leq E$, $d(a_2, x) \in (E, F]$.

In this case, we claim
$$a_1 \arr a_2 \Longleftrightarrow x <_E a_2 \hspace{1 cm} a_2 \arr a_1 \Longleftrightarrow a_2 <_E x$$

The implications from right to left hold by the definition of $\arr$.

For the implication from left to right, we consider only the case $a_1 \arr a_2$, since the other is similar. By definition, there exists some $y$ such that $a_1 \preceq_E y \preceq_E a_2$. Then $x \preceq_E a_1 \preceq_E y$, so $x \preceq_E y$. Since $y \preceq_E a_2$, then $x \preceq_E a_2$. Since $d(x, a_2) \not\leq E$, we have $x <_E a_2$.
\end{claimproof}

Claims 1 and 2 dispose of the cases in which $d(a_1, a_2) \not\in (E, F]$. By Claim 3, if $d(a_1, a_2) \in (E, F]$ and $a_1 \arr a_2$, we may complete amalgam by determining $a_1 <_E a_2$, and vice versa if $a_2 \arr a_1$. If $d(a_1, a_2) \in (E, F]$ and neither $a_1 \arr a_2$ nor $a_2 \arr a_1$, we may complete the amalgam by arbitrarily determining either $a_1 <_E a_2$ or $a_2 <_E a_1$.
\end{proof}

\section{Observations on the Catalog}

In this section, we make some observations about the structures produced in Theorem \ref{theorem:expandLambda}. The first is that they fall into the regime suggested by Question \ref{q:ordexp}. Next, we show the lattice of $\emptyset$-definable equivalence relations is not affected by the expansion by subquotient orders, and so each finite distributive lattice $\Lambda$ appears as such a lattice for some homogeneous finite-dimensional permutation structure. Finally, we observe the structures produced by this method satisfy the Primitivity Conjecture.

Before addressing Question \ref{q:ordexp} in this context, we note that the reduct should also be a strong amalgamation class in order for there to be an expansion by a generic linear order. However, this follows from the hypotheses.

\begin{lemma}
Let $\Gamma$ be the \fraisse limit of a strong amalgamation class and $\Gamma^{red}$ a homogeneous reduct. Then $\Gamma^{red}$ is also the \fraisse limit of a strong amalgamation class.
\end{lemma}
\begin{proof}

We will use the equivalence that a homogeneous structure is the \fraisse limit of a strong amalgamation class iff it has trivial algebraic closure, i.e. $acl(A) = A$ for any subset $A$. However, any formula witnessing some non-trivial algebraic closure in $\Gamma^{red}$ will still witness it in $\Gamma$.
\end{proof}

We now address Question \ref{q:ordexp}.

\begin{proposition} 
Let $\Lambda$ be a finite distributive lattice, and $\vec \Gamma$ the \fraisse limit of a well-equipped lift of the class of all finite $\Lambda$-ultrametric spaces. Then $\vec \Gamma$ is the expansion by a linear order of a homogeneous proper reduct. Furthermore, if $\vec \Gamma$ is primitive, the linear order may be taken to be generic.
\end{proposition}
\begin{proof}
In $\vec \Gamma$, there are definable subquotient orders $<_{\bbzero, E}$ from $\bbzero$ to $E$ and $<_{F, \bbone}$ from $F$ to $\bbone$. We now take the reduct $\vec \Gamma^{red}$ forgetting a single generic subquotient order $<_{E, F}$ from $E$ to $F$. Then the expansion of $\vec \Gamma^{red}$ by the linear order $<_{F, \bbone}[<_{E, F}[<_{\bbzero, E}]]$ yields $\vec \Gamma$.

If $\vec \Gamma$ were primitive then we would have $E = \bbzero$ and $F = \bbone$. Thus the linear order we add back at the end would just be $<_{E, F}$, which would be a generic linear order.
\end{proof}

By iterating the above proof, we get the corollary that any homogeneous finite-dimensional permutation structure produced by our construction, with a lattice $\Lambda$ of $\emptyset$-definable equivalence relations, can be produced by starting with the generic $\Lambda$-ultrametric space and iteratively adding linear orders, with each step yielding a homogeneous structure producible by our construction. However, the corresponding result replacing the generic $\Lambda$-ultrametric space with a structureless set is not true. To see this, consider the full product $\Q^2$ (see Example \ref{ex:fullQ2}), which may be presented in a language of four linear orders. However, any reduct to three linear orders is no longer homogeneous (one can either directly find a violation of homogeneity, or check the classification of homogeneous 3-dimensional permutation structures in Chapter \ref{c:3dim}).

\begin{lemma} \label{lemma:expansionlattice}
Let $\Lambda$ be a finite distributive lattice, and $\Gamma$ the generic $\Lambda$-ultrametric space. Let $\AA^*$ be as in Theorem \ref{theorem:amalg}, with \fraisse limit $\Gamma^*$. Then the lattice of $\emptyset$-definable equivalence relations in $\Gamma^*$ is isomorphic to $\Lambda$; in particular, it is distributive.
\end{lemma}

\begin{corollary}[Representation Theorem]\label{theorem:Representation}
Let $\Lambda$ be a finite distributive lattice. 
Then there is a homogeneous finite dimensional permutation
structure whose lattice of $\emptyset$-definable equivalence 
relations is isomorphic to $\Lambda$.
\end{corollary}

\begin{proof}[Proof of Lemma \ref{lemma:expansionlattice}]
Any $2$-type $p$ realized in $\MM$ by the ordered pair $(a,b)$ is encoded by the data $(E_p,p_{ord})$ where
$E_p=E_{d(a,b)}$ and $p_{ord}$ is the type of $(a,b)$ in the language restricted to the subquotient orders
$<_{E_i}$, which records whether $a<_{E_i} b$, $b<_{E_i} a$, or $a, b$ are $<_{E_i}$-incomparable, for each $i \in [n]$.

We may consider such a $2$-type as a minimal nontrivial $\emptyset$-definable
binary relation on $\MM$. Let $E_p'$ denote the smallest equivalence relation containing
the relation $p$: i.e., the transitive and reflexive closure of the symmetrized type $p\union p^{\op}$,
where $p^{\op}$ is the type of $(b,a)$.  

\begin{claim*}
Let $p$ be a $2$-type realized in $\MM$. Then 
$$E_p'=E_p$$
\end{claim*}

Given the claim, consider an arbitrary $\emptyset$-definable equivalence relation $E$ on $\MM$.
This is the union of the $2$-types $p$ contained in $E$, and hence is the join of the
equivalence relations $E_p'$ generated by those types. Since $E_p'=E_p$, this join
lies in $\Lambda$.

\begin{claimproof}

If $a=b$ then $E_p$ and $E_p'$ are both equality, so suppose $a\ne b$. Since the pair satisfies $E_p(a,b)$, it follows that $E_p'\includedin E_p$.

Conversely, suppose that we have a pair $c,d$ satisfying $E_p(c,d)$.
Let $q$ be the type of $(c,d)$. We extend $(c,d)$ to a triangle $(a,c,d)$ by setting
$\tp(a,c)=\tp(a,d)=p$. 
If this triangle belongs to $\AA^*$, then by homogeneity it embeds into 
$\Gamma^*$ over $(c,d)$, so $E_p'(c,d)$ holds and we are done. 
(And the proof shows $E_p=p\circ p^{\op}$.)

So let us check the conditions on the triangle $(a,c,d)$, to show it is in $\AA^*$.

\goodbreak
\begin{figure}[h]
\begin{diagram}
c \bullet & & \rTo{(E_q, q_{ord})} & & \bullet d \\
& \luTo_{(E_p, p_{ord})} &  & \ruTo_{(E_p, p_{ord})} & \\
& & \underset{a}{\bullet} & & \\
\end{diagram}
\caption{}
\label{fig:3}
\end{figure}

1.~The $\Lambda$-metric triangle inequality holds. The labels are $(E_p,E_p,E_q)$ where $E_q=d(c,d)\le E_p$.
So this is clear.

2.~The relations $<_{E_i}$ are subquotient orders. We only need to check transitivity, and since for each $i$ we have $a<_{E_i} c,d$, or $c,d <_{E_i} a$, or $a$ is $<_{E_i}$-incomparable with $c,d$, this is clear.
\end{claimproof}
\end{proof}

\begin{proposition} \label{prop:sqoprim}
Let $\Gamma$ be the generic $n$-dimensional permutation structure, in the language $\set{<_1, ..., <_n}$, and let $<$ be a definable linear order on $\Gamma$. Then there is an $i \in [n]$ such that $< = <_i$ or $< = <_i^{opp}$.  
\end{proposition}
\begin{proof}
Note that $<$ must be a union of 2-types $\cup q_i$, and for each 2-type, exactly one of it and its opposite must be appear as some $q_i$. We may assume $q_0 \vdash x <_i y$ for all $i$. If the conclusion is false, then for each $i \in [n]$, there must be a type $p_i = q_j$ for some $j$, such that $p_i \vdash y <_i x$. 

Now consider the \textit{partial structure} on $X = \set{x_1, ..., x_{n+1}}$ given by setting $p_i(x_i, x_{i+1})$ for each $i \in [n]$. For each $i \in [n]$, looking at $<_i$ on $X$ gives a directed acyclic graph, whose transitive closure is a partial order in which $x_1$ and $x_n$ are $<_i$-incomparable. This can then be completed to a linear order in which $x_n <_i x_1$.
 
 Each $<_i$ is a linear order in the resulting structure, which is thus a substructure of $\Gamma$. However, we have $x_1 < ... < x_n$ but $x_n < x_1$. Thus $<$ is not transitive on this structure, and so does not define a linear order on $\Gamma$.
\end{proof}

\begin{corollary}
Suppose Conjecture \ref{conj:sqoconj} is true. Then Conjecture \ref{conj:primitivity}, the Primitivity Conjecture, is true.
\end{corollary}
\begin{proof}
Let $\Gamma$ be a primitive homogeneous $n$-dimensional permutation structure. Assuming Conjecture \ref{conj:sqoconj}, $\Gamma$ is interdefinable with $\Gamma'$, the fully generic $m$-dimensional permutation structure for some $m \leq n$. If $m=n$, we are finished. If not, then by Proposition \ref{prop:sqoprim} the additional linear orders in $\Gamma$ must be equal to one of the $m$ liner orders in $\Gamma'$, up to reversal, proving the Primitivity Conjecture.
\end{proof}

\section{The Number of Orders Needed for the Representation Theorem}\label{subsection:OrderNumber}
Although we have finished the proof of the Representation Theorem, the translation from equivalence relations and subquotient orders to linear orders via Proposition \ref{prop:sqolininter} made no attempt to minimize the number of linear orders used. Before continuing, we take this section to provide a better bound on the number of orders required.

\begin{question} \label{qu:represent1}
Given a homogeneous finite-dimensional permutation structure $\Gamma$ presented in a language of equivalence relations and subquotient orders, what is the minimal $n$ such that $\Gamma$ is quantifier-free interdefinable with an $n$-dimensional permutation structure? 
\end{question}

A first remark is that it is not true that one needs at most $n$ orders to represent a structure with at most $2^n$ 2-types. This is illustrated by the full product $\Q^2$ (see Example \ref{ex:fullQ2}), which only has 8 non-trivial 2-types, but requires 4 linear orders to represent. This fact can be seen by noting this structure does not appear in a our classification of homogeneous 3-dimensional permutation structures in Chapter \ref{c:3dim}. However, we will now give a direct argument, which shows how increasing the number of linear orders adds more freedom beyond simply increasing the number of 2-types, due to transitivity constraints.

\begin{lemma} \label{lemma:linearLattice}
Let $\Gamma$ be a homogeneous 3-dimensional permutation structure. Then $\Gamma$ does not have incomparable $\emptyset$-definable equivalence relations.
\end{lemma}
\begin{proof}
Suppose $E, F$ are incomparable $\emptyset$-definable equivalence relations in $\Gamma$. Then $E$ contains some 2-type $p$ and its opposite, and $F$ contains some distinct 2-type $q$ and its opposite. By possibly permuting and reversing orders, we may assume $p \vdash x <_i y$ for all $i$, and $q \vdash \set{x <_1 y, x <_2 y, y <_3 x}$.

Consider the following amalgamation diagram, with $p=tp(x,b)$, $q = tp(b, y)$.

\begin{figure}[h]
\begin{diagram}
 x \odot & &  & & \odot y  \\
& \rdTo_{p} &  & \ruTo_{q} & \\
& & {\bullet}_b & & \\
\end{diagram}
\caption{}
\label{fig:4}
\end{figure}

Transitivity forces $x <_1 y$ and $x <_2 y$. If $x <_3 y$, then $tp(x, y) = p$, and if $y <_3 x$ then $tp(x, y) = q$. Either choice leads to a violation of transitivity of these equivalence relations.
\end{proof}

 As a final note, evidence from Chapter \ref{c:3dim} hints that in the case when the lattice of $\emptyset$-definable equivalence relations is linear, the trivial bound of $n$ orders to represent structures with at most $2^n$ 2-types might be attainable.

In this section, we concerned with the following subquestion, asking how many linear orders are needed for Corollary \ref{theorem:Representation}, the Representation Theorem.

\begin{question} \label{qu:represent2}
Given a lattice $\Lambda$, what is the minimal $n$ such that $\Lambda$ is isomorphic to the lattice of $\emptyset$-definable equivalence relations of some homogeneous $n$-dimensional permutation structure?
\end{question}

The basic idea is that when the lattice is a chain, the encoding can be done with nearly maximal efficiency. We thus partition the meet-irreducibles into chains and encode each chain separately.

\begin{proposition}\label{Proposition:Countchain}
Let $\set{E_i}$ be a chain of equivalence relations of height at most $2^n-1$, and let $<$ be a linear order such that each $E_i$ is $<$-convex. Then there exist $n$ linear orders $\set{<_j}$, such that each $<_j$ is quantifier-free definable in $((E_i), <)$, and each $E_i$ is quantifier-free definable in $(< , (<_j))$.
\end{proposition}

\begin{proof}
We may suppose the chain has height exactly $2^n-1$, and does not contain equality or the universal relation, since those are already definable. Extend the chain to length $2^n+1$ by letting $E_0$ be equality, and $E_{2^n}$ be the universal relation. In the language $(<, (<_j)_{j=1}^n)$, where the $<_j$ are binary relations, enumerate the non-trivial quantifier-free 2-types containing the formula $x<y$ (and so these 2-types merely specify whether $x <_j y$ for each $j$) as $(p_i : 1 \leq i \leq 2^{n})$. We will use each pairing $p_i \cup p_i^{op}$ to produce an equivalence relation.

Define the relation $R_j(x, y) \Leftrightarrow \bigvee_{\set{i | (x <_j y) \in p_i}} (E_i(x, y) \wedge \neg E_{i-1}(x, y))$. We now define $<_j$ to be the canonical irreflexive, asymmetric extension of $(x < y) \wedge R_j(x, y)$, i.e. $x <_j y \Leftrightarrow ((x< y) \wedge R_j(x, y)) \vee ((y < x) \wedge \neg R_j(y, x))$. Thus each $<_j$ is quantifier-free definable from $((E_i), <)$.

Conversely, we see $x(E_i \bs E_{i-1})y \Leftrightarrow (x < y \wedge type(x,y)=p_{i}) \vee (y < x \wedge type(x,y)=p_{i}^{op})$, so each $E_i$ is quantifier-free definable from $(< , (<_j))$.

It remains to check that the relations $<_j$ are actually linear orders. Clearly $<_j$ is irreflexive and asymmetric, so we check that it is total and has no cycle. To see it is total, first assume $x < y$. Then if $x \not <_j y$, we must have $\neg R_j(x, y)$. But then $y <_j x$.

 We now show there is no cycle. Suppose $x <_j y <_j z <_j x$. Up to a change of notation (cyclically permuting the variables, and reversing $<$ if needed), we may assume $x < y < z$. Let $d(x,y) = \min (E_i : E_i(x,y))$, and define $d(x,z)$ and $d(y,z)$ similarly. By the $<$-convexity of the $E_i$, we have $d(x,y), d(y,z) \leq d(x,z)$, and so the triangle inequality gives either $d(x,z) = d(x,y)$ or $d(x,z) = d(y,z)$. In the first case, our definition of $<_j$ gives $x <_j z$ iff $x<_j y$, since $type(x,z) = type(x, y)$ in the language $(E_i), <)$, and similarly in the second case it gives $x <_j z$ iff $y <_j z$.
\end{proof}

$\Q^n$ with the lexicographic order can naturally be expressed in a language of one order $<$ and a chain of $n-1$ $<$-convex equivalence relations $E_i$, $1 \leq i \leq n-1$, given by $x E_i y$ iff $x$ and $y$ agree in the first $i$ coordinates. The lexicographic $\Q^2$ requires two orders to define, and the lexicographic $\Q^3$ requires three. One might expect each new convex equivalence relation to require an additional order, but we already see the exponential growth implied by the above proposition illustrated by the lexicographic $\Q^4$, which also only requires three orders.

\begin{lemmacorollary}\label{Corollary:CountLattice}
Let $\Lambda$ be a finite distributive lattice, $\Lambda_0$ the poset of meet-irreducibles of $\Lambda \bs \set{\bbzero, \bbone}$, and $\LL$ a set of chains covering $\Lambda_0$. Then the dimension of the permutation structure needed for the Representation Theorem is at most $|\LL| + \sum_{L \in \LL} \ceil{\log_2(|L|+1)}$.
\end{lemmacorollary}
\begin{proof}
We construct a well-equipped lift of the generic $\Lambda$-ultrametric space. For each $L \in \LL$, enumerate its elements from least to greatest as $\set{\lambda_{L, i}}_{i=1}^{n_L}$ and perform the following.
\begin{enumerate}
\item For each $i$, expand by a generic subquotient order from $E_{L,i}$ to $E_{L.i+1}$.
\item Expand by a generic subquotient order from $E_{L, n_L}$ to $\bbone$.
\item If $\bbzero$ is meet-irreducible, expand by a generic subquotient order from $\bbzero$ to $E_{L,1}$.
\end{enumerate}

In the case $\bbzero$ is meet-reducible, then for each $L \in \LL$ we may use Lemma \ref{lemma:sqoproduct} to define a subquotient order from $\bbzero$ to $E_{L, 1}$. Now, for each $L \in \LL$, by taking compositions of the respective subquotient orders, we may define a linear order $<_L$ convex with respect to each $E_{L, i} \in L$. We have so far added $|\LL|$ linear orders.

Now each $L \in \LL$ considered with the order $<_L$ satisfies the hypotheses of Proposition \ref{Proposition:Countchain}, and so the equivalence relations labeled by elements of $L$ are definable after the addition of $\ceil{\log_2(|L|+1)}$ linear orders. Thus all the elements of $\Lambda_0$ are definable a further expansion by $\sum_{L \in \LL} \ceil{\log_2(|L|+1)}$ linear orders, and all of $\Lambda$ is definable from the elements of $\Lambda_0$. The quantifier-free-definability conditions from Proposition \ref{Proposition:Countchain} ensure that the structure obtained by adding these linear orders and removing the equivalence relations is still homogeneous with the same lattice of $\emptyset$-definable equivalence relations, and so we obtain the bound of the statement.
\end{proof}

When considering the possible optimality of this bound, note that in a homogeneous $\Lambda$-ultrametric space, if $E, F \in \Lambda$ are incomparable then every $E$-class must meet every $F$-class in the same $E \join F$ class. Thus no linear order can be convex with respect to incomparable equivalence relations.

Consider $\Q^2$ with two equivalence relations, each given by equality in one of the coordinates. There are no further non-trivial $\emptyset$-definable equivalence relations, and so the above corollary gives a bound of four orders to define this structure, which is in fact the number needed. We do not know if the bound of Corollary \ref{Corollary:CountLattice} is tight in general.
\chapter{Towards the Completeness of the Catalog}\label{chap:completeness}
\section{Introduction}

In this chapter, we provide some evidence for the conjecture that our construction produces all homogeneous finite-dimensional permutation structures. In Chapter \ref{c:3dim}, we provide further evidence for this conjecture by classifying the homogeneous 3-dimensional permutation structures.

 The construction of Chapter \ref{chap:construction} suggests the following conjecture. In particular, the backward direction is \ref{theorem:Representation}, and Lemma \ref{lemma:expansionlattice} shows Conjecture \ref{conj:sqoconj} would imply the forward direction.

\begin{conjecture} [The Distributivity Conjecture] \label{Conjecture:Distributivity}
A finite lattice is isomorphic to the lattice of $\emptyset$-definable equivalence relations in some homogeneous finite dimensional permutation structure iff it is distributive.
\end{conjecture}

The main result of the first section is a partial result in the forward direction, which actually follows from a more general result involving the \textit{infinite index property} (see Definition \ref{Definition:IIP} and Lemma \ref{lemma:IIPDIS}). 

\newtheorem*{theorem:Distributivity}{\bf{Theorem \ref{theorem:Distributivity}}}
\begin{theorem:Distributivity} 
\emph{Let $\Lambda$ be the lattice of $\emptyset$-definable equivalence relations in a homogeneous finite dimensional permutation structure $\MM$. If the reduct of $\MM$ to the language of equivalence relations from $\Lambda$ is homogeneous, then $\Lambda$ is distributive.}
\end{theorem:Distributivity} 

This implies that if there is a homogeneous finite-dimensional permutation structure with a non-distributive lattice of $\emptyset$-definable equivalence relations, it must have a completely different method of construction. We cannot simply start with a homogeneous $\Lambda$-ultrametric space and take an expansion preserving the equivalence relation structure. The more general result alluded to above, involving the infinite index property, sharpens this considerably.

In the next section of this chapter, we consider homogeneous permutation structures in which all minimal forbidden substructures are of order 2. Such a structure is necessarily primitive, since defining equivalence relations requires forbidding 3-types, and thus the Primitivity Conjecture predicts its form. This is confirmed, in this special case, by the following proposition.

\newtheorem*{Proposition:2-types}{\bf{Proposition \ref{Proposition:2-types}}}
\begin{Proposition:2-types}
\emph{Let $\KK$ be an amalgamation class of $n$-dimensional permutation structures. If no 3-type compatible with the allowed 2-types is forbidden, then the forbidden 2-types collectively specify that certain orders agree up to reversal.}
\end{Proposition:2-types}

\section{Towards the Necessity of Distributivity}

We first define a property that will be satisfied by the lattice of $\emptyset$-definable equivalence relations in any homogeneous finite-dimensional permutation structure. 

\begin{definition}\label{Definition:IIP}
Let $\MM$ be a structure with a transitive automorphism group. 
\begin{enumerate}
\item For $\emptyset$-definable equivalence relations $F \subset E$ on $\MM$, set $$[E:F] = |C/F|$$ for $C$ any $E$-class, and call this the \textit{index of $F$ in $E$}.
\item Let $\Lambda$ be the lattice of $\emptyset$-definable equivalence relations on $\MM$. Then $\Lambda$ has the \textit{infinite index property (IIP)} if whenever $F \subset E$ for $E, F \in \Lambda$, $[E:F]$ is infinite.
\end{enumerate}
\end{definition}

\begin{lemma}\label{lemma:IIP}
 Given a homogeneous finite-dimensional permutation structure $\MM$, let $\Lambda$ be the lattice of $\emptyset$-definable equivalence relations. Then $\Lambda$ satisfies the IIP.
 \end{lemma}

\begin{proof}
Let $E,F \in \Lambda$ with $F < E$, and let $p$ be a 2-type of orders that is realized in $E \bs F$. Then $p$ is an intersection of linear orders, and so gives a partial order on a given $E$-class. Since the structure's automorphism group is transitive, this partial order has no maximal elements, and so contains an infinite linear order $L$. The 2-type between any pair of elements of $L$ is $p$, and thus every pair is $E$-related but not $F$-related. Thus $[E: F] \geq |L|$ is infinite.
\end{proof}

The following lemma is reminiscent of Neumann's lemma that a group cannot be covered by finitely many cosets of subgroups of infinite index. We generalize from the group-theoretic setting, replacing the equivalence relations induced by subgroups with equivalence relations in some lattice, but impose the stronger condition that this lattice must satisfy the IIP.

\begin{lemma}\label{lemma:Neumann}
Let $\Lambda$ be a finite lattice of equivalence relations satisfying the IIP, and let $E \in \Lambda$ with $C$ an $E$-class. Let $\set{B_i}_{i \in I}$ be a finite set of equivalence classes of certain equivalence relations in $\Lambda$ such that for each of the corresponding equivalence relations $E_{i} \in \Lambda$, $E_{i} \not \geq E$. Then there exists some $c \in C \bs \bigcup_{i\in I} B_i$.
\end{lemma}

\begin{proof}
We proceed by induction on the height of $E$ in the lattice. In the base case, $E$ is equality, and the claim is vacuous. 

Now assume $E$ is higher up. We wish to work entirely below $E$, so we replace each $B_i$ with $B_i \cap C$, and replace each $E_i$ with $E_i \cap E$. Let $E'$ be a maximal equivalence relation strictly below $E$. Then, for any $E_{i}$, we cannot have $E_{i} \geq E'$ unless $E_{i} = E'$.  Since we are trying to avoid finitely many equivalence classes, and by the IIP there are infinitely many $ E'$-classes in $C$, we may pick an $E'$-class $C' \subset C$ that is not equal to any of the $B_i$. Then, letting $I' = \set{i \in I | E_i \neq E'}$, by induction we can find a $c \in C'\bs \bigcup_{i \in I'}(B_i \cap C) \subseteq C \bs \bigcup_{i\in I} B_i$.
\end{proof}

We now use the above lemma to prove a one-point extension property for homogeneous $\Lambda$-ultrametric spaces where $\Lambda$ satisfies the IIP. However, we restate the property using amalgamation classes and diagrams.

\begin{lemma}\label{lemma:Extension}
Let $\KK$ be an amalgamation class of $\Lambda$-ultrametric spaces with Fra\"\i ss\'e limit $\MM$, and suppose the lattice of $\emptyset$-definable equivalence relations in $\MM$ satisfies the IIP. Let $K \in \KK$ with $K = X \cup \set{b}$, $e \in \Lambda$, and $B = \set{b, y}$ with $d(b,y) = e$. Then the canonical amalgam of $K$ and $B$ is in $\KK$. (Alternatively, the amalgamation diagram below, with arbitrary first factor, a single point in the base, and a single extension point in the second factor, can be completed by canonical amalgamation.)
\end{lemma}

\begin{figure}
\begin{diagram}
X \bigcirc & &  & & \odot y \\
& \luLine_{\set{e_i}_{i \leq |X|}} &  & \ruLine_{e} & \\
& & \underset{b}{\bullet} & & \\
\end{diagram}
\caption{}
\label{fig:5}
\end{figure}

\begin{proof}
Identify the elements of $\Lambda$ with the corresponding $\emptyset$-definable equivalence relations in $\MM$. We choose $y \in \MM$ using Lemma \ref{lemma:Neumann} with $E=e$, $C$ the $e$-class containing $b$, and $\set{B_i}$ the set of equivalence classes containing $b$ for every equivalence relation below $e$, as well as, for every $x_i \in X$, the equivalence classes containing $x_i$ for equivalence relations not above $e$. Note that the first group of $B_i$ ensures that $d(b,y) = e$.

Now fix an $x_i \in X$, let $d(x_i, b) = e_i$, and let $d(x_i, y) = e'$. From the second group of $B_i$, we have $e \leq e'$. Thus, using the triangle inequality for the upper bound, we have $e \leq e' \leq e \vee e_i$. Then, $e_i \leq e' \vee e = e'$, so $e_i \leq e'$ as well. Thus $e' = e \vee e_i$.
\end{proof}

\begin{lemma}\label{lemma:IIPDIS}
Let $\MM$ be a homogeneous structure equipped with a set of equivalence relations forming a finite lattice $\Lambda$ satisfying the IIP. Then $\Lambda$ is distributive.
\end{lemma}
\begin{proof}
\setcounter{claim}{0}
Let $\KK$ be the amalgamation class corresponding to $\MM$, viewed as a class of $\Lambda$-ultrametric spaces.
\begin{claim}
Suppose both factors of the amalgamation diagram shown in Figure \ref{fig:6} are contained in an amalgamation class of $\Lambda$-ultrametric spaces, for every $e,f,g \in \Lambda$. Then $\Lambda$ is distributive.
\end{claim}

\begin{figure}[h]
\begin{diagram}
  x \odot & &  & & & & & & \odot y \\
  & \rdLine_e \rdLine(6,2)^f& & & & & & \ldLine_f \ldLine(6,2)^e &  \\
  & &v \bullet & & \hLine_{e \vee f}& & \bullet w \\
  & & &\rdLine_{e \vee g} & &\ldLine_{f \vee g} & \\
  & & \leftcurve{(e \vee g) \wedge (f \vee g)} & & \underset{u}{\bullet} & & \rightcurve{g}  & &  \\
\end{diagram}
\caption{}
\label{fig:6}
\end{figure}
\begin{claimproof}
Let $d(x, y) = h$ in the completed diagram. Then $h \leq e$ and ${h \leq f}$ by the triangle inequality. Going from $x$ to $u$ via $y$, the triangle inequality gives $(e \vee g) \wedge (f \vee g) \leq h \vee g \leq (e \wedge f) \vee g$, and so the claim follows.
\end{claimproof}

\begin{claim}
$\KK$ contains both factors of the amalgamation diagram shown in Figure \ref{fig:6} for every $e,f,g \in \Lambda$.
\end{claim}
\begin{claimproof}
Because $\Lambda$ satisfies the IIP, we may use the one-point extension property in Lemma $\ref{lemma:Extension}$ to build the factors of Figure \ref{fig:6} one point at a time. For the second factor (omitting $x$), we start with $y$ as a base point, and proceed as in Figure \ref{fig:7}, adding $v$, $w$, and $u$, in order.

\begin{figure}[h]
\begin{diagram}
v\odot &&   v\odot & & \rDash^{e \vee f} & & \odot w&& w\odot \\
\vLine_e && & \rdLine_{e} & & \ruLine_{f} & &&  \vLine^{e \vee f} & \rdLine(2,4)_f \rdDash(4,2)^{f \vee g}\\
y \bullet && & & y \bullet & & && v\odot & & \rDash^{e \vee g} & &\odot u\\
&&&&&&&& & \rdLine_e & & \ruLine_g\\
&&&&&&&& & & y \bullet \\
\end{diagram}
\caption{}
\label{fig:7}
\end{figure}

The construction of the first factor proceeds similarly, starting with $x$ as a base point, but in the step corresponding to the last diagram of Figure \ref{fig:7}, we put $d(x,u) = (e \vee g) \wedge (f \vee g)$ instead of $d(x, u) = g$. Since the diagram is then completed by canonical amalgamation, in we must check that $f \vee ((e \vee g) \wedge (f \vee g)) = f \vee g$ and $e \vee ((e \vee g) \wedge (f \vee g)) = e \vee g$. Since the arguments are identical, we will only consider the first identity.

Clearly, $f \vee ((e \vee g) \wedge (f \vee g)) \leq f \vee g$, and $f \vee ((e \vee g) \wedge (f \vee g)) \geq f \vee g$, since  $(e \vee g) \wedge (f \vee g) \geq g$.
\end{claimproof}
\end{proof}

The following example of a homogeneous $\Lambda$-ultrametric space, where $\Lambda$ is the non-distributive diamond lattice $M_3$, shows that some IIP-like property is necessary in Lemma \ref{lemma:IIPDIS}.

\begin{example} \label{ex:non-dis}
Let $\Lambda$ be a copy of $M_3$ with non-trivial equivalence relations $E_1, E_2,$ and $E_3$. Let $\MM$ be an $\Lambda$-ultrametric space on 4 points $\set{a,b,c,d}$, defined as follows. 

There are 2 $E_i$-classes, for each $i$. The $E_1$-classes are $\set{a, b}$ and $\set{c, d}$. The $E_2$-classes are $\set{a, c}$ and $\set{b, d}$. The $E_3$-classes are $\set{a, d}$ and $\set{b, c}$.

It is easy to check that $\MM$ is homogeneous, yet its lattice of $\emptyset$-definable equivalence relations is $\Lambda$.
\end{example}

The above example can be viewed as taking the equivalence relations to be parallel classes of lines in the affine plane over ${F}_2$.

\begin{corollary} \label{cor:metric-dis}
Let $\Lambda$ be a finite lattice, and suppose the class of all finite $\Lambda$-ultrametric spaces is an amalgamation class. Then $\Lambda$ is distributive.
\end{corollary}
\begin{proof}
For each pair $E, F \in \Lambda$ with $E \leq F$, it is easy to construct a finite $\Lambda$-ultrametric space consisting of a single $F$-class containing arbitrarily many $E$-classes.
 As all these structures embed into the \fraisse limit, its lattice of $\emptyset$-definable equivalence relations must satisfy the IIP. The results follows from Lemma \ref{lemma:IIPDIS}.
\end{proof}

\begin{theorem}\label{theorem:Distributivity}
 Let $\Lambda$ be the lattice of $\emptyset$-definable equivalence relations in a homogeneous finite dimensional permutation structure $\MM$. If the reduct of $\MM$ to the language of equivalence relations from $\Lambda$ is homogeneous, then $\Lambda$ is distributive.
\end{theorem} 

\begin{proof}
By Lemma \ref{lemma:IIP}, $\Lambda$ satisfies the IIP. Thus, we may apply Lemma \ref{lemma:IIPDIS} to conclude.
\end{proof}

\section{Forbidding Only 2-Types}

In Cameron's homogeneous permutations, whenever a 2-type is forbidden it forces one order to be equal to another, up to reversal. This need not always be the case in higher dimensional homogeneous permutation structures. Consider $\Q^4$ as a lexicographic order. As discussed following Proposition \ref{Proposition:Countchain}, this structure can be defined using 3 orders. Since each $\emptyset$-definable equivalence relation is separated by a single 2-type and its opposite, forbidding a 2-type and its opposite causes one equivalence relation to collapse to one beneath it, and the resulting structure is the lexicographic $\Q^3$. Since we only forbid one pair of 2-types, we cannot have made one order equal to another, up to reversal.

We know two ways to forbid 2-types in homogeneous finite-dimensional permutation structures: collapse one order to another, up to reversal, or collapse one equivalence relation to another. In the second case, we must forbid some 3-types compatible with the allowed 2-types, i.e. 3-types respecting the transitivity constraints such that the restriction to any 2 variables is an allowed 2-type. The following result hints that these two constructions may be typical to some degree.

\begin{proposition}\label{Proposition:2-types}
\setcounter{claim}{0}
Let $\KK$ be an amalgamation class of $n$-dimensional permutation structures. If no 3-type compatible with the allowed 2-types is forbidden, then the forbidden 2-types collectively specify that certain orders agree up to reversal.
\end{proposition}

The following diagram will reappear prominently in Chapter \ref{c:3dim}, so we extract its definition and an associated lemma from the proof of Proposition \ref{Proposition:2-types}.

\begin{definition} \label{definition:majoritydiagram}
Give 2-types $p,q,r$, the \textit{$(p,q,r)$-majority diagram} is the following amalgamation diagram, where $x_1 \ra{q} x_3$ holds (and follows from $x_1 \ra{q} x_2 \ra{q} x_3$), but is not drawn.
\begin{figure}[h]
\begin{diagram}
& & \overset{x_1}{\bullet} & & \\
& \ruTo^{p} & \dTo_q  & \rdTo^{q} & \\
a_1 \odot &\rTo^p & \bullet_{x_2}  & \rTo^r & \odot a_2 \\
& \rdTo_{q} & \dTo_q & \ruTo_{r} & \\
& & \underset{x_3}{\bullet} & & \\
\end{diagram}
\caption{}
\label{fig:8}
\end{figure}
\end{definition}

\begin{lemma} \label{lemma:majoritydiagram}
There is a unique solution to the $(p,q,r)$-majority diagram, given by $a_1 <_i a_2$ iff $<_i$ is true in at least two of $p,q,$ and $r$. 
\begin{proof}
Note that every pair of $p,q,$ and $r$ appears on a path of length two from $a_1$ to $a_2$. Thus if $<_i$ is true in some pair, the path containing that pair will force $a_1 <_i a_2$, and vice versa if $<_i$ is false in some pair.
\end{proof}
\end{lemma}

\begin{proof}[Proof of Proposition \ref{Proposition:2-types}]
 First, assume that no two orders are equal up to reversal, since otherwise we could pass to a reduct in which this is the case. Then we must show that any 2-type is realized.
 
  We will use the notation $[n] = \set{1, ..., n}$ and let $\binom{[n]}{k}$ denote the set of $k$-subsets of $[n]$. Given 2-types $p$ and $q$ and $X \subset [n]$, we say $p$ is an $X$-approximation to $q$ if $p$ and $q$ agree on the orders indexed by elements of $X$. 
  
  Given a 2-type $t$, we will prove $t$ is realized using induction on the size of approximations to $t$. By assumption, given any $X \in \binom{[n]}{2}$, there is a non-forbidden $X$-approximation to $t$; otherwise the orders in $X$ would have to be equal up to reversal. 

\begin{claim*}
 Let $p,q,r$ be 2-types realized in $\KK$. Then both factors of the $(p,q,r)$-majority diagram are in $\KK$. 
 \end{claim*}
\begin{claimproof}
We consider only the first factor, since there is a symmetric argument for the second factor. First, note that by transitivity of 2-types, the first factor is the unique amalgam of two triangles, as shown in the following diagram: 

\begin{figure}[h]
\begin{diagram}
& & \overset{x_3}{\odot} & & \\
& \ruTo^{p} & \dTo_q   \\
a_1 \bullet &\rTo^p & \bullet_{x_2}  \\
& \rdTo_{q} & \dTo_q\\
& & \underset{x_1}{\odot} & & \\
\end{diagram}
\caption{}
\label{fig:9}
\end{figure}


Both factors of this last diagram are in $\KK$, because no 3-types compatible with the allowed 2-types are forbidden, so the only constraint is transitivity. However, because each triangle has two equal sides pointing from or to the same point, all transitivity constraints are satisfied.
\end{claimproof}

We are now ready to treat the inductive step of our argument. Suppose the 2-type $t$ has a non-forbidden $X$-approximation for every $X \in \binom{[n]}{k}$. Fix $Y \in \binom{[n]}{k+1}$. Without loss of generality, we may assume $Y = [k+1]$. Let $p$ be a $\set{1, ..., k}$-approximation to $t$, $q$ a $\set{1, ..., k-1, k+1}$-approximation to $t$, and $r$ a $\set{k,k+1}$-approximation to $t$. Then, by Lemma \ref{lemma:majoritydiagram}, the solution to the corresponding $(p,q,r)$-majority diagram will be a $Y$-approximation to $t$.
\end{proof}

As mentioned in the introduction, Proposition \ref{Proposition:2-types} confirms a special case of the Primitivity Conjecture.
\chapter{The Classification of Homogeneous 3-Dimensional Permutation Structures}\label{c:3dim}
\section{Introduction}
   \newtheorem*{lemma:trianglereduce}{\bf{Lemma \ref{lemma:trianglereduce}}}

In this chapter, we continue to explore the conjectural completeness of our catalog of homogeneous finite-dimensional permutation structures. The main result of this chapter is the classification of the homogeneous 3-dimensional permutation structures, confirming Conjecture \ref{conj:sqoconj} in this case. In particular the Primitivity Conjecture holds, as does the Distributivity conjecture; in fact, the lattice of $\emptyset$-definable equivalence relations is linear. Most of the structures in the catalog are \textit{compositions} of simpler structures, in the following sense.

\begin{definition}
Given structures $\Gamma_1, \Gamma_2$ in disjoint languages, the \textit{composition of $\Gamma_1$ with $\Gamma_2$}, denoted  $\Gamma_1[\Gamma_2]$, is the structure obtained by expanding $\Gamma_1$ with an equivalence relation $E$, and replacing the points of $\Gamma_1$ by $E$-classes that are copies of $\Gamma_2$.
\end{definition}

If $\Gamma_1$ and $\Gamma_2$ both have distributive lattices of $\emptyset$-definable equivalence relations, then as the lattice of $\emptyset$-definable equivalence relations in $\Gamma_1[\Gamma_2]$ is the lattice sum, it is also distributive.

\begin{definition}
We use $\Gamma^{(g)}_i$ to denote the generic $i$-dimensional permutation structure; in particular $\Gamma^{(g)}_0$ is a set equipped only with equality.
\end{definition}

\begin{theorem}[The Classification] \label{thm:catalog}
Let $(\Gamma, <_1, <_2, <_3)$ be a homogeneous 3-dimensional permutation structure. Then $\Gamma$ is quantifier-free interdefinable with one of the following 16 structures.

\begin{enumerate}
\item $\Gamma$ has no non-trivial $\emptyset$-definable congruence
	\begin{enumerate}
	\item $\Gamma$ is primitive: $\Gamma = \Gamma^{(g)}_1, \Gamma^{(g)}_2,$ or $\Gamma^{(g)}_3$.
	\item $\Gamma$ is imprimitive: $\Gamma$ is the expansion  by a generic linear order of $\Gamma^{(g)}_1[\Gamma^{(g)}_j]$, for $j \in \set{0, 1}$.
	\end{enumerate}
\item $\Gamma$ has a non-trivial $\emptyset$-definable congruence
	\begin{enumerate}
	\item $\Gamma$ is a repeated composition of primitive structures: For any multisubset $I \subset \set{1,2}$ such that $|I|>1$ and $\sum_{i \in I} 2^i \leq 8$, $\Gamma$ is the composition in any order of $\Gamma^{(g)}_i$ for $i \in I$.
	\item $\Gamma$ is a composition of primitive and imprimitive structures: Let $\Gamma^*$ denote the structure from $(1b)$ with $j=0$. Then $\Gamma = \Gamma^*[\Gamma^{(g)}_1]$ or $\Gamma^{(g)}_1[\Gamma^*]$.
	\end{enumerate}
\end{enumerate}

\end{theorem}

The classification proceeds in two stages. First, we confirm the Primitivity Conjecture for 3 orders using explicit amalgamation arguments. Then for the imprimitive case, we pick a minimal non-trivial equivalence relation $E$. The Primitivity Conjecture makes it fairly clear what happens on $E$-classes, and some analysis of the type structure between $E$-classes eventually allows us to carry out an inductive classification.

\begin{corollary} \label{cor:confirmsqo}
   Every homogeneous 3-dimensional permutation structure is interdefinable with the \fraisse limit of some well-equipped lift of the class of all finite $\Lambda$-ultrametric spaces, for some distributive lattice $\Lambda$.
\end{corollary}

Despite the fact that assuming the correctness of Conjecture \ref{conj:sqoconj} gives a simple description of \textit{all} finite-dimensional permutation structures, it is difficult to determine the corresponding catalog for a \textit{fixed number} of linear orders. (This problem is discussed Chapter \ref{subsection:OrderNumber}).

Thus, Corollary \ref{cor:confirmsqo} is not proven by first producing a conjectural classification and then confirming it. Rather, it is proven by observing that all the structures appearing in the classification may be presented appropriately.

Finally, although we have a confirmation of Conjecture \ref{conj:sqoconj} in the case of 3 orders, a plausible exceptional imprimitive structure arises in the analysis (see Lemma \ref{lemma:dense}) that is ultimately shown not to exist. However, the proof of non-existence makes use of the limited type structure with 3 orders, and it seems possible similar structures will appear in the richer languages afforded by more orders. 

\section{The Primitive Case}
In this section, we classify the primitive homogeneous 3-dimensional permutation structures, obtaining the following.

\begin{theorem}\label{theorem:primclass}
The primitive homogeneous 3-dimensional permutation structures are as predicted by the Primitivity Conjecture.
\end{theorem}

The main results needed for the proof are Proposition \ref{Proposition:2-types} and the lemmas below.

\begin{lemma}\label{lemma:trianglereduce}
Suppose $(\Gamma; <_1, <_2, <_3)$ is homogeneous and realizes all 3-types. Then $\Gamma$ is generic.
\end{lemma}

\begin{lemma} \label{lemma:primtarget}
Let $\Gamma$ be a primitive homogeneous 3-dimensional permutation structure. Then all 3-types involving realized 2-types are realized.
\end{lemma}

\begin{proof}[Proof of Theorem \ref{theorem:primclass}] 
By Lemma \ref{lemma:primtarget} all 3-types involving realized 2-types are realized. Thus, if no 2-types are forbidden, all 3-types are realized, and so by Lemma \ref{lemma:trianglereduce}, the resulting structure is generic. If some 2-types are forbidden, then by Proposition \ref{Proposition:2-types}, some orders agree up to reversal. Thus the resulting structure is quantifier-free interdefinable with a primitive homogeneous 2-dimensional permutation structure. By the classification in \cite{Cameron}, these satisfy the Primitivity Conjecture.
\end{proof}

A variant of Lemma \ref{lemma:trianglereduce}, generalizing beyond 3 linear orders but weakening the bound, appears as Lemma \ref{lemma:4gen} below. Generalizing Lemma \ref{lemma:trianglereduce} while maintaining the bound at realizing all 3-types would be a major step toward proving the Primitivity Conjecture.

The proof of Lemma \ref{lemma:primtarget} is lengthy, and relies on simple and explicit amalgamation arguments given in Section \ref{section:amalgLemmas}. Repeated use of these amalgamation lemmas and a suitable case division eventually yield the desired result.

\subsection{Reduction to 3-types} 
The following lemma strengthens an argument appearing in the proof of \cite{Cameron}*{Theorem 1}, and is the first step toward proving Lemma \ref{lemma:trianglereduce}.

We remark that proving Lemma \ref{lemma:trianglereduce} seems to be the most promising place for an application of some version of Lachlan's Ramsey argument.

\begin{lemma}\label{lemma:4gen}
Let $\Gamma$ be a homogeneous $k$-dimensional permutation structure that realizes all configurations on $n-1$ points, where $n$ satisfies
$$\frac{n!}{(n-\ell)!} > 2^\ell k \text{ for } \ell = \floor{n/2}$$
Then $\Gamma$ is generic.

More precisely, any configuration on $N \geq n$ points is contained in the unique amalgam of $(N-1)$-point configurations.
\end{lemma}

Note this does not apply in the case $n=4, k=3$, so proving Lemma \ref{lemma:trianglereduce} will require additional argument.

\begin{proof}
Let $A$ be a structure on $N$ points. Let a \textit{pairing} be an $\ell$-set of unordered pairs of points from $A$, with each point appearing in at most one pair. A pairing is \textit{separated} if, for every $i \leq k$, there is a pair $(a_i, a'_i)$ such that $a_i$ and $a'_i$ are not $<_i$-adjacent; otherwise the pairing is \textit{unseparated}.

\begin{claim*}
 There is at least one separated pairing on $A$.
\end{claim*} 

\begin{claimproof}The number of pairings is given by $\binom{n}{2\ell} \binom{2\ell}{2_1, 2_2, ..., 2_\ell}/\ell! = \frac{n!}{2^\ell \ell! (n-2\ell)!}$. We will now show the number of unseparated pairings is at most $k \binom{n-\ell}{\ell}$. First suppose $k=1$. If $N$ is even, there is only 1 unseparated pairing. If $N$ is odd, the pairing is determined after choosing any one of the odd-indexed points to be omitted, so there are $\ceil{n/2}$ such pairings. In both cases, there are $\binom{n-\ell}{\ell}$. For larger $k$, note that an unseparated pairing must be unseparated with respect to at least one order, so there are at most $k\binom{n-\ell}{\ell}$ such pairings. By the inequality in the hypothesis, we are done.
\end{claimproof}

Now let $P$ be a separated pairing. By extending $A$ by a single point, we may, in every order, make one pair non-adjacent. Thus, after extending $A$ by at most $\ell-1$ points, an extension we will denote by $A^*$, we may assume that every pair in $P$ is non-adjacent in every order.

Let $(a_1, a'_1)$ be a pair from $P$, and let $F_1 = A^* \bs{\set{a_1}}, F'_1 = A^* \bs{\set{a'_1}}$, and $B_1 = A^* \bs{\set{a_1, a'_1}}$. By assumption, for every $i \leq k$, there is a point $b_i \in B_1$ that is $<_i$-between $a_1$ and $a'_1$. Thus $A^*$ is the unique amalgam of $F_1$ and $F'_1$ over $B$.

We may recursively continue this process on each factor until we have gone through all the pairs in $P$. At the end, each factor will look like a copy of $A^*$ with $\ell$ points removed, and so have size $N-1$. Now $A^*$ is the unique amalgam of all these factors, and $A \subset A^*$.
\end{proof}

\begin{lemma:trianglereduce}
Suppose $(\Gamma; <_1, <_2, <_3)$ is homogeneous and realizes all 3-types. Then $\Gamma$ is generic.
\end{lemma:trianglereduce}
\begin{proof}
By Lemma \ref{lemma:4gen}, if $\Gamma$ realizes all 4-point configurations, it is generic.

Let $A = (\set{a,b,c,d}; <_1, <_2, <_3)$ be a substructure of $\Gamma$. There are three possible pairings: $P_1 = \set{\set{a,b}, \set{c,d}}, P_2 = \set{\set{a,c}, \set{b,d}}, P_3 = \set{\set{a,d}, \set{b,c}}$. Each order can be unseparated in at most one pairing, so if all the pairings are unseparated, each must be so with respect to a different order. By possibly relabeling the points, we may assume that $a <_1 b <_1 c <_1 d$, and by relabeling orders we may assume that $P_i$ is unseparated with respect to $<_i$.

Thus, we have that $a,c$ and $b,d$ must be $<_2$-adjacent, and $a,d$ and $b,c$ must be $<_3$-adjacent.

We may extend $A$ by a single point, $e$, that lies between $a$ and $b$ with respect to $<_1$, lies between $a$ and $c$ with respect to $<_2$, and lies between $b$ and $c$ with respect to $<_3$, and label the resulting structure $A^*$. Then, viewing $B = \set{e, c, d}$ as the base of an amalgamation digram with $F = B \cup \set{a}$ the first factor and $F' = B \cup \set{b}$ the second, we have that $A^*$ is the unique amalgam.

We now show that $F$ and $F'$ have separated pairings, and so are contained in the unique amalgam of certain 3-types. For $F$, $P = \set{\set{a,c}, \set{e,d}}$ is separated, since $e$ and $d$ are never $<_2$-adjacent and only $<_3$-adjacent if $b$ and $d$ were $<_3$-adjacent, in which case $a$ and $c$ not $<_3$-adjacent. For $F'$, $P' = \set{\set{b,c}, \set{e,d}}$ is separated, since $e$ and $d$ are never $<_3$-adjacent and only $<_2$-adjacent if $a$ and $d$ were $<_2$-adjacent, in which case $b$ and $c$ are not $<_2$-adjacent.
\end{proof}

\subsection{Notation}
We now begin preparing for the proof of Lemma \ref{lemma:primtarget}.

Given three linear orders, there are 8 2-types, which we may associate with the vertices of the unit cube $\set{\pm 1}^3$ based on whether or not $<_i$ holds in the 2-type. The unit cube is bipartite, with one part consisting of the following four types at Hamming distance 2, while the other part consists of their opposites.

$$0: \la{123} \hspace{1 cm} 1: \ra{23} \hspace{1 cm} 2: \ra{13} \hspace{1 cm} 3: \ra{12}$$

We now introduce notation for three families of 3-types that will recur in our analysis. From left to right, the 3-types below will be denoted $p \To q$, $p \From q$, and $C_3(p,q,r)$.

\begin{figure}[h]
\begin{diagram}
 \bullet & & \rTo^{q} & & \bullet  &  \bullet & & \rTo^{q} & & \bullet & \bullet & & \rTo{q} & & \bullet  \\
& \luTo_{p} &  & \ruTo_{p}  & & &  \rdTo_p &  & \ldTo_p & & &  \luTo_p &  & \ldTo_r & \\
& & \underset{}{\bullet} & & & & &  \underset{}{\bullet} & & & & & \underset{}{\bullet} & & \\
\end{diagram}
\caption{}
\label{fig:10}
\end{figure}

\subsection{Amalgamation Lemmas} \label{section:amalgLemmas}

In the following three lemmas, $(p,q,r,s)$ is taken to be some permutation of the four 2-types $(0,1,2,3)$, defined in the previous section. These lemmas will be the basis of the proof of Lemma \ref{lemma:primtarget}.

\begin{lemma} \label{lemma:1}
Suppose $p \To q$ is forbidden. Then one of each of the following combinations of 3-types is forbidden.
\begin{enumerate}[(A)]
\item $(p \To r$ \text{and} $C_3(p,r,s))$ or $(r \From q$ \text{and} $C_3(p,q,r))$
\item $p \From q$ or $q \From p$
\end{enumerate}
\end{lemma}
\begin{proof}
For $(A)$, we amalgamate one 3-type from each pair over an edge of type $r$. In the diagrams below, we assume $p \To r$ is realized; the arguments assuming $C_3(p,r,s)$ is realized are similar.

\begin{figure}[h]
\begin{diagram}
& & {\bullet} & & & & & & {\bullet} & & \\
& \ruTo^{p} &  & \luTo^{r} & & & & \ruTo^{p} &  & \rdTo^{p} & \\
x\odot & & \uTo{r} & & \odot y  & & x\odot & & \uTo{r} & & \odot y \\
& \rdTo_{p} &  & \ldTo_{q} & & & & \rdTo_{p} &  & \ldTo_{q} & \\
& & {\bullet} & & & & & & {\bullet} & &\\
\end{diagram}
\caption{}
\label{fig:11}
\end{figure}

We wish to argue that the only way to complete either diagram is to take $tp(x, y) = p$. This is clear by transitivity for the diagram on the right. For the diagram on the left, note that since $p$ and $r$ are at Hamming distance 2, $p$ and $r^{opp}$ agree on exactly 2 orders, as do $p$ and $q^{opp}$. Thus, by transitivity, $tp(x, y)$ must agree with $p$ in all 3 orders.

For $(B)$, we use the $(p^{opp}, q, p)$-majority diagram (see Definition \ref{definition:majoritydiagram}), and then take Lemma \ref{lemma:5} into account.
\end{proof}

\begin{lemma} \label{lemma:4}
Suppose $p \To q$, $C_3(p,q,r)$, and $C_3(p,q,s)$ are forbidden. If $p$ and $q$ are realized, then $q \From p$ is realized.
\end{lemma}
\begin{proof}
We complete the following amalgamation diagram.

\begin{figure}[h]
\begin{diagram}
 x \odot & &  & & \odot y  \\
& \rdTo_{p} &  & \ruTo_{q} & \\
& & {\bullet} & & \\
\end{diagram}
\caption{}
\label{fig:12}
\end{figure}

By assumption $tp(x, y) \neq p, r^{opp}, s^{opp}$. The remaining types, except $q$, are ruled out by transitivity.

\end{proof}

\begin{lemma} \label{lemma:5}
Suppose $p \To q$ is forbidden. If $p$ and $q$ are realized, then $q \To p$ is realized.
\end{lemma}
\begin{proof}
We complete the following amalgamation diagram.

\begin{figure}[h]
\begin{diagram}
 x \odot & &  & & \odot y  \\
& \rdTo_{p} &  & \ldTo_{q} & \\
& & {\bullet} & & \\
\end{diagram}
\caption{}
\label{fig:13}
\end{figure}

By assumption $tp(x, y) \neq p$.  The remaining types, except $q^{opp}$, are ruled out by transitivity.
\end{proof}

\subsection{Case Division}
The proof of Lemma \ref{lemma:primtarget} proceeds by consideration of several cases. However, the following lemma provides a uniform point of departure.

\begin{lemma}
If $\Gamma$ is primitive and omits a 3-type then without loss of generality it omits the 3-type of type $(0 \Rightarrow 1)$ while realizing the 2-types $0,1$.
\end{lemma}
\begin{proof}
We may assume that at least 3 of the 2-types $0,1,2,3$, say $0,1,2$ after relabeling, are realized, since otherwise we reduce to the case of fewer linear orders. 

If any 3-type $p \Rightarrow q$ or $p \Leftarrow q$ is forbidden while $p$ and $q$ are realized, then by reversing the orders and changing the language we may assume that $0 \Rightarrow 1$ is forbidden. So assume this is not the case. 

By the above paragraphs, we may construct the standard $(0,1,2)$-majority diagram, which shows $3$ is realized as well. Up to a change of language, the forbidden 3-type must be of the form $C_3(0,1,2)$. But this is a substructure of the unique solution to the $(1^{opp}, 0^{opp}, 2^{opp})$-majority diagram.
\end{proof}

\begin{remark}
Although we can assume the 2-type $1$ is realized, we may not want to, since it breaks the symmetry between $1,2,$ and $3$. Thus, this will not be assumed unless otherwise noted.
\end{remark}

We now divide into cases the ways Lemma \ref{lemma:primtarget} might fail.
\begin{list}{}{}
\item \textbf{Case 1:} All 3-types of type $0 \To p$ ($p=1,2,3)$ are forbidden, and $0$ is realized.
\item \textbf{Case 2:} For a given pair of 2-types $p,q$ at Hamming distance 2, at most 2 3-types of type $p \To q$ are forbidden, and $0,1$ are realized.
\begin{list}{}{}
\item \textbf{Case 2.1:} There exist $p,q,r$ at Hamming distance 2 such that $p \To q$ and $p \To r$ are forbidden.
\item \textbf{Case 2.2:} For any $p,q,r$ at Hamming distance 2, at most one of $p \To q$ and $p \To r$ is forbidden.
\end{list}  
\end{list}

We also wish to divide Case 1 into subcases. If $0 \To p$ is forbidden, for $p = 1,2,3$, consider the directed graph with vertex set $\set{1,2,3}$, and an edge $(p,q)$ when the type $C_3(0, p, q)$ is forbidden. By Lemma \ref{lemma:1}, for any arrangement $(p,q,r)$ of the vertices, either $(p,q)$ or $(q,r)$ is an edge. Thus, $D$ contains a symmetric edge $p \leftrightarrow q$, which we may assume is $1 \leftrightarrow 2$, and $D$ has at least 4 edges.

We now subdivide Case 1 as follows.

\begin{list}{}{}
\item \textbf{Case 1.1:} $D$ has 6 edges.
\item \textbf{Case 1.2:} $D$ has 5 edges.
\item \textbf{Case 1.3:} $D$ has 4 edges.
\end{list}

\subsection{Proof of Lemma \ref{lemma:primtarget}}
The proof proceeds by starting with the assumptions of one of the subcases and then repeatedly applying the amaglamation lemmas \ref{lemma:1} (A,B), \ref{lemma:4}, and \ref{lemma:5} until reaching a contradiction. This contradiction could either be that a structure is both forbidden and realized, or could be a violation of the primitivity constraint by the appearance of a definable equivalence relation. 

More explicitly, the 2-types $p_1, ..., p_k$ generate a definable equivalence relation if every 3-type on points $x,y,z$ satisfying the following two conditions is forbidden.
\begin{enumerate}
\item $tp(x, y), tp(y, z) \in \set{p_1, ..., p_k}$
\item $tp(x, z) \not\in \set{p_1, ..., p_k}$
\end{enumerate} 

The proofs are presented in tables. In each line, some 3-type is shown to be realized or forbidden. The reason is given; if the reason is one of the amalgamation lemmas then the assignment of $(p,q,r,s)$ is given; finally the previous lines used are given. When one of the amalgamation lemmas is used with opposite types, so for example $p \From q$ is assumed forbidden rather than $p \To q$, an ``R'' (for ``reversed'') is appended to the name of the lemma.

In the tables, we assume all 2-types are realized; after each table is a remark noting the alterations required if some 2-type is forbidden.

\begin{longtable}{|c|c|c|c|c|c|}
\caption{Case 1.1} \\
\hline 
\multicolumn{1}{|c|}{Line} &
 \multicolumn{1}{c|}{Realized} & 
 \multicolumn{1}{c|}{Forbidden} &
\multicolumn{1}{c|}{Reason} &
 \multicolumn{1}{c|}{(p,q,r,s)} & 
 \multicolumn{1}{c|}{Used} \\ \hline 
\endfirsthead
   1. & & $0 \To p$ & Case 1 & & \\
   2. & & $C_3(0,p,q)$ & Case 1.1 & & \\
   3. & $p \From 0$ & & \ref{lemma:4} & & 1,2 \\
   4. & & $0 \From p$ & \ref{lemma:1}B & &1,3 \\
   \hline
      \multicolumn{6}{|l|}{{Now the 2-type $0$ generates an equivalence relation, contradicting primitivity.}} \\ \hline
\end{longtable}

\begin{remark}
This proof works with some 2-type forbidden. If $1,2$, or $3$ is forbidden, then the corresponding case of line 4 follows without needing line 3.
\end{remark}

The treatment of the remaining cases follows the same scheme at somewhat greater length, and shows that the amalgamation lemmas given previously suffice to complete the analysis in the primitive case. Other methods will be required in the next section.

\begin{longtable}{|c|c|c|c|c|c|} 
\caption{Case 1.2} \\
\hline 
\multicolumn{1}{|c|}{Line} &
 \multicolumn{1}{c|}{Realized} & 
 \multicolumn{1}{c|}{Forbidden} &
\multicolumn{1}{c|}{Reason} &
 \multicolumn{1}{c|}{(p,q,r,s)} & 
 \multicolumn{1}{c|}{Used} \\ \hline 
\endfirsthead
   1. & & $0 \To p$ & Case 1 & & \\
   2. &  &$C_3(0, p, q)$    & Case 1.2 & & \\
      &  &except $C_3(0,3,2)$   &  & & \\
   3. &$C_3(0, 3, 2)$  &  & Case 1.2 & & \\
   4. &$1 \From 0$  &   & \ref{lemma:4} &(0,1,2,3) & 1,2 \\
   5. &$2 \From 0$  &   & \ref{lemma:4} &(0,2,1,3) & 1,2 \\
   6. & &$0 \From 1$   & \ref{lemma:1}B &(0,1,2,3) & 1,4 \\
   7. & &$0 \From 2$   & \ref{lemma:1}B &(0,2,1,3) & 1,5 \\
   8. &$0 \From 3$ &   & Primitivity & & 1,6,7 \\
   9. & &$3 \From 0$   & \ref{lemma:1}B &(0,3,1,2) & 1,8 \\
   10. & &$3 \From 2$   & \ref{lemma:1}AR &(3,0,2,1) & 5,9 \\
   11. & &$3 \To 2$  or $2 \To 3$  & \ref{lemma:1}BR &(3,2,0,1) & 10 \\
   12. & $2 \To 0$ & & \ref{lemma:5} & (0,2,1,3) & 1 \\
   13. & $3 \To 0$ & & \ref{lemma:5} & (0,3,1,2) & 1 \\
   14. & &$3 \To 0$ or $2 \To 0$   & \ref{lemma:1}A &(3,2,0,1) & 3,8,11 \\
   & & & &or (2,3,0,1) & \\
   \hline 
   \multicolumn{6}{|l|}{{Now line 14 contradicts lines 12 and 13.}} \\ \hline
\end{longtable}

\begin{remark}
This proof works with some 2-types forbidden. By assumption, the types $2$ and $3$ are realized. The assumption that the type $1$ is realized only appears in line $4$, which becomes unnecessary if $1$ is forbidden since line 4 is only used for line $6$.
\end{remark} 

Case 1.3 requires further subdivision according to our assumptions on the directed graph $D$. We draw the directed graph $D$ corresponding to each of the further subcases.

\pagebreak

\begin{figure}[h]
\begin{diagram}
1 \bullet & & \rDouble^{1.3.1} & & \bullet 2 & 1 \bullet & & \rDouble^{1.3.2} & & \bullet 2 & 1 \bullet & & \rDouble^{1.3.3} & & \bullet 2 \\
& \luDouble &  &  & & &  \rdTo &  & \ldTo & & &  \luTo &  & \ruTo & \\
& & \underset{3}{\bullet} & & & & &  \underset{3}{\bullet} & & & & & \underset{3}{\bullet} & & \\
\end{diagram}
\caption{}
\label{fig:14}
\end{figure}

\begin{longtable}{|c|c|c|c|c|c|} 
\caption{Case 1.3.1} \\
\hline 
\multicolumn{1}{|c|}{Line} &
 \multicolumn{1}{c|}{Realized} & 
 \multicolumn{1}{c|}{Forbidden} &
\multicolumn{1}{c|}{Reason} &
 \multicolumn{1}{c|}{(p,q,r,s)} & 
 \multicolumn{1}{c|}{Used} \\ \hline 
\endfirsthead
   1. & & $0 \To p$ & Case 1 & & \\
   2. &  &$C_3(0, 1, p)$,   & Case 1.3.1 & & \\
   & &$C_3(0,p,1)$ & & & \\
   3. &$C_3(0, 3, 2)$,&  & Case 1.3.1 & & \\
   & $C_3(0,2,3)$   & & & & \\
   4. &$1 \From 0$  &   & \ref{lemma:4} &(0,1,2,3) & 1,2 \\
   5. & &$0 \From 1$   & \ref{lemma:1}B &(0,1,2,3) & 1,4 \\
   6. &$0 \From 2$ or $0 \From 3$ & & Primitivity & & 1,5\\
   7. &$0 \From 3$ & & W.l.o.g & & 6\\
   8. & &$3 \From 0$   & \ref{lemma:1}B &(0,3,1,2) & 1,8 \\
   9. & &$3 \From 2$   & \ref{lemma:1}AR &(3,0,2,1) & 3,8 \\
   10. & &$3 \To 2$ or $2 \To 3$   & \ref{lemma:1}BR &(3,2,0,1) & 9 \\
   11. & $2 \To 0$ & & \ref{lemma:5} & (0,2,1,3) & 1 \\
   12. & $3 \To 0$ & & \ref{lemma:5} & (0,3,1,2) & 1 \\
   13. & &$3 \To 0$ or $2 \To 0$   & \ref{lemma:1}A &(3,2,0,1)& 3,10 \\
   & & & &  or (2,3,0,1) & \\
   \hline
      \multicolumn{6}{|l|}{{Now line 13 contradicts lines 11 and 12.}} \\ \hline
\end{longtable}

\begin{remark}
This proof works with some 2-types forbidden. By assumption, the types $2$ and $3$ are realized. The assumption that the type $1$ is realized only appears in line $4$, which becomes unnecessary if $1$ is forbidden since line 4 is only used for line $5$.
\end{remark}

\begin{longtable}{|c|c|c|c|c|c|} 
\caption{Case 1.3.2} \\
\hline 
\multicolumn{1}{|c|}{Line} &
 \multicolumn{1}{c|}{Realized} & 
 \multicolumn{1}{c|}{Forbidden} &
\multicolumn{1}{c|}{Reason} &
 \multicolumn{1}{c|}{(p,q,r,s)} & 
 \multicolumn{1}{c|}{Used} \\ \hline 
\endfirsthead
   1. & & $0 \To p$ & Case 1 & & \\
   2. &  &$C_3(0, 1, p)$, $C_3(0,2,p)$   & Case 1.3.2 & & \\
   3. &$C_3(0, 3, p)$  &  & Case 1.3.2 & & \\
   4. &$1 \From 0$  &   & \ref{lemma:4} &(0,1,2,3) & 1,2 \\
   5. &$2 \From 0$  &   & \ref{lemma:4} &(0,2,1,3) & 1,2 \\
   6. & &$0 \From 1$   & \ref{lemma:1}B &(0,1,2,3) & 1,4 \\
   7. & &$0 \From 2$   & \ref{lemma:1}B &(0,2,1,3) & 1,5 \\
   8. & &$0 \From 3$   & \ref{lemma:1}AR &(0,1,3,2) & 3,6 \\
   \hline
      \multicolumn{6}{|l|}{{Now $0$ generates an equivalence relation.}} \\ \hline
\end{longtable}

\begin{remark}
By assumption, all 2-types are realized.
\end{remark}

\begin{longtable}{|c|c|c|c|c|c|} 
\caption{Case 1.3.3} \\
\hline 
\multicolumn{1}{|c|}{Line} &
 \multicolumn{1}{c|}{Realized} & 
 \multicolumn{1}{c|}{Forbidden} &
\multicolumn{1}{c|}{Reason} &
 \multicolumn{1}{c|}{(p,q,r,s)} & 
 \multicolumn{1}{c|}{Used} \\ \hline 
\endfirsthead
   1. & & $0 \To p$ & Case 1 & & \\
   2. &$C_3(0, p, 3)$  &  & Case 1.3.3 & & \\
   3. & &$2 \From 1$   & \ref{lemma:1}A &(0,1,2,3) & 1,2 \\
   4. & &$1 \From 2$   & \ref{lemma:1}A &(0,2,1,3) & 1,2 \\
   5. &$1 \From 2$ &   & \ref{lemma:5}R &(2,1,0,3) & 3 \\
   \hline
      \multicolumn{6}{|l|}{{However, line 5 contradicts line 4.}} \\ \hline
\end{longtable}

\begin{remark}
By assumption, all 2-types are realized.
\end{remark}

For Case 2.1, we may assume that $0 \To 1$ and $0 \To 2$ are forbidden, and thus $0 \To 3$ is realized.

\begin{longtable}{|c|c|c|c|c|c|} 
\caption{Case 2.1} \\
\hline 
\multicolumn{1}{|c|}{Line} &
 \multicolumn{1}{c|}{Realized} & 
 \multicolumn{1}{c|}{Forbidden} &
\multicolumn{1}{c|}{Reason} &
 \multicolumn{1}{c|}{(p,q,r,s)} & 
 \multicolumn{1}{c|}{Used} \\ \hline 
\endfirsthead
   1. & &$0 \To 1$, $0 \To 2$ & Case 2.1 & & \\
   2. &$0 \To 3$  & & Case 2 & & \\
   3. &$1 \To 0$  & & \ref{lemma:5} & (0,1,2,3) &1 \\
   4. &$2 \To 0$  & & \ref{lemma:5} & (0,2,1,3) &1 \\
   5. &  &$3 \From 1$, $C_3(0, 1, 3)$  & \ref{lemma:1}A &(0,1,3,2) &1,2 \\
   6. &  &$3 \From 2$, $C_3(0, 2, 3)$  & \ref{lemma:1}A &(0,2,3,1) &1,2 \\
   7. &$3 \From 0$&   & Case 2 & & 5,6 \\
   8. &$1 \From 3$ &   & \ref{lemma:5}R& (3,1,0,2)& 5 \\
   9. &$2 \From 3$ &   & \ref{lemma:5}R& (3,2,0,1)& 6 \\
   10. & &$1 \From 0$ or $0 \From 1$   & \ref{lemma:1}B & (0,1,2,3) & 1 \\
   \hline   
\end{longtable}

We now split into cases based on line 10.
\begin{list}{}{}
\item \textbf{Case 2.1.1:} $1 \From 0$ is forbidden.
\item \textbf{Case 2.1.2:} $0 \From 1$ is forbidden.
\end{list}

\begin{longtable}{|c|c|c|c|c|c|} 
\caption{Case 2.1.1} \\
\hline 
\multicolumn{1}{|c|}{Line} &
 \multicolumn{1}{c|}{Realized} & 
 \multicolumn{1}{c|}{Forbidden} &
\multicolumn{1}{c|}{Reason} &
 \multicolumn{1}{c|}{(p,q,r,s)} & 
 \multicolumn{1}{c|}{Used} \\ \hline 
\endfirsthead
   11. & &$1 \From 0$   & Case 2.1.1 & & 10\\
   12. & &$1 \From 2$, $C_3(1,3,2)$   & \ref{lemma:1}AR &(1,0,2,3) & 4,11 \\
   13. &$1 \To 3$ &   & \ref{lemma:4}R& (3,1,0,2)& 5,12 \\
   14. & &$3 \To 1$   & \ref{lemma:1}BR &(3,1,0,2) & 5,13 \\
   15. & &$3 \To 0$   & \ref{lemma:1}AR &(1,0,3,2) & 8,11 \\
   16. &$3 \To 2$ &  & Case 2 & & 14,15 \\
   17. & &$1 \From 3$ or $3 \To 2$ & \ref{lemma:1}AR & (1,2,3,0) & 12\\
   \hline
      \multicolumn{6}{|l|}{{However, line 17 contradicts line 8 and line 16.}} \\ \hline
\end{longtable}

\begin{longtable}{|c|c|c|c|c|c|} 
\caption{Case 2.1.2} \\
\hline 
\multicolumn{1}{|c|}{Line} &
 \multicolumn{1}{c|}{Realized} & 
 \multicolumn{1}{c|}{Forbidden} &
\multicolumn{1}{c|}{Reason} &
 \multicolumn{1}{c|}{(p,q,r,s)} & 
 \multicolumn{1}{c|}{Used} \\ \hline 
\endfirsthead
   11. & &$0 \From 1$   & Case 2.1.2 & & 10\\
   12. & &$0 \From 2$ or $2 \To 1$, $C_3(0,2,1)$   & \ref{lemma:1}AR &(0,1,2,3) & 11 \\
\hline
\end{longtable}

We now split into cases based on line 12.
\begin{list}{}{}
\item \textbf{Case 2.1.2.1:} $0 \From 2$ is forbidden.
\item \textbf{Case 2.1.2.2:} $2 \To 1, C_3(0,2,1)$ is forbidden.
\end{list}

\begin{longtable}{|c|c|c|c|c|c|} 
\caption{Case 2.1.2.1} \\
\hline 
\multicolumn{1}{|c|}{Line} &
 \multicolumn{1}{c|}{Realized} & 
 \multicolumn{1}{c|}{Forbidden} &
\multicolumn{1}{c|}{Reason} &
 \multicolumn{1}{c|}{(p,q,r,s)} & 
 \multicolumn{1}{c|}{Used} \\ \hline 
\endfirsthead
   13. & &$0 \From 2$  & Case 2.1.2.1 & & 12 \\
   14. &$0 \From 3$ &   & Case 2& & 11,13 \\
   15. & &$3 \To 1$, $C_3(0,3,1)$   & \ref{lemma:1}AR &(0,1,3,2) & 11,14 \\
   16. & &$3 \To 2$, $C_3(0,3,2)$   & \ref{lemma:1}AR &(0,2,3,1) & 13,14 \\
      \hline
         \multicolumn{6}{|l|}{{Now, $0 \cup 3$ generates an equivalence relation.}} \\ \hline
\end{longtable}

\begin{longtable}{|c|c|c|c|c|c|} 
\caption{Case 2.1.2.2} \\
\hline 
\multicolumn{1}{|c|}{Line} &
 \multicolumn{1}{c|}{Realized} & 
 \multicolumn{1}{c|}{Forbidden} &
\multicolumn{1}{c|}{Reason} &
 \multicolumn{1}{c|}{(p,q,r,s)} & 
 \multicolumn{1}{c|}{Used} \\ \hline 
\endfirsthead
   13. & &$2 \To 1$, $C_3(0,2,1)$  & Case 2.1.2.2 & & 12 \\
   14. &$2 \From 0$ &   & \ref{lemma:4}& (0,2,1,3) & 1,6,13 \\
   15. & &$0 \From 2$   & \ref{lemma:1}B &(0,2,1,3) & 1,14 \\
         \hline
            \multicolumn{6}{|l|}{{Now $0 \From 2$ is forbidden, and we may finish as in Case 2.1.2.1.}} \\ \hline
\end{longtable}

\begin{remark}
This proof works with some 2-types forbidden. By assumption, $1$ and $3$ are realized. Assume $2$ is forbidden. Case 2.1.1 ends at line 15 with a contradiction of the Case 2 assumption, since $3 \To 1$, $3 \To 0$, and $3 \To 2$ will all be forbidden. Only lines 4 and 9 depend on 2 being realized, and those are only used in line 12, which would hold anyway if 2 were forbidden. Case 2.1.2.1 works as before, since lines 4 and 9 aren't used anywhere. Also, there is no need for Case 2.1.2.2, since we know $2 \From 0$ is forbidden.
\end{remark}

For Case 2.2, we may assume $0 \To 1$ is forbidden, and thus $0 \To 2$ and $0 \To 3$ are realized.

\begin{longtable}{|c|c|c|c|c|c|} 
\caption{Case 2.2} \\
\hline 
\multicolumn{1}{|c|}{Line} &
 \multicolumn{1}{c|}{Realized} & 
 \multicolumn{1}{c|}{Forbidden} &
\multicolumn{1}{c|}{Reason} &
 \multicolumn{1}{c|}{(p,q,r,s)} & 
 \multicolumn{1}{c|}{Used} \\ \hline 
\endfirsthead
   1. & &$0 \To 1$& Case 2.2 & & \\
   2. &$0 \To 2$, $0 \To 3$  & & Case 2.2 & & \\
   3. &  & $2 \From 1$, $C_3(0,1,2)$ & \ref{lemma:1}A & (0,1,2,3) &1,2 \\
   4. &  & $3 \From 1$, $C_3(0,1,3)$ & \ref{lemma:1}A & (0,1,3,2) &1,2 \\
   5. & $2 \From 0$, $2 \From 3$  & & Case 2.2 & &3\\
   6. & $3 \From 0$, $3 \From 2$  & & Case 2.2 & &4 \\
   7. & & $3 \To 1$, $C_3(2,3,1)$   & \ref{lemma:1}AR & (2,1,3,0) & 3,5 \\
   8. & & $2 \To 1$, $C_3(3,2,1)$   & \ref{lemma:1}AR & (3,1,2,0) & 4,6 \\
   9. & $3 \To 0$, $3 \To 2$  & & Case 2.2 & &7\\
   10. & $2 \To 0$, $2 \To 3$  & & Case 2.2 & &8\\
   11. & &$0 \From 1$, $C_3(3,1,0)$ & \ref{lemma:1}A & (3,1,0,2) &7,9\\
   12. & & $C_3(2,1,0)$ & \ref{lemma:1}A & (2,1,0,3) &8,10\\
   \hline
   \multicolumn{6}{|l|}{{Now $0 \cup 2 \cup 3$ generates an equivalence relation.}} \\ \hline
\end{longtable}

\begin{remark}
By assumption, all 2-types are realized.
\end{remark}

This completes the proof of Lemma \ref{lemma:primtarget}.

\section{The Imprimitive Case}
 Now suppose $\Gamma$ an imprimitive homogeneous 3-dimensional permutation structure, and let $E$ be a minimal non-trivial $\emptyset$-definable equivalence relation. We make the following initial case division of the imprimitive case.

\begin{list}{}{}
\item \textbf{Case 1:} $E$ is convex with respect to $<_1, <_2, <_3$, and thus is a congruence.
\item \textbf{Case 2:} $E$ is not convex with respect to at least one of $<_1, <_2, <_3$. Without loss of generality, we assume $E$ is not $<_1$-convex.
\end{list}

In Case 1, we may proceed inductively by factoring out $E$, noting that the resulting structure now omits a 2-type, and so $\Gamma$ will be a composition of a homogeneous 3-dimensional permutation structure with one fewer 2-type available and a primitive homogeneous 3-dimensional permutation structure. 

Our goal for Case 2 will be to show that $\Gamma$ is still determined by its restriction to $E$-classes and by the $E$-quotient of the reduct of $\Gamma$ forgetting all orders for which $E$ is non-convex.

While not necessary for our arguments, it is perhaps psychologically helpful to note that the lattice of $\emptyset$-definable equivalence relations in $\Gamma$ must be linear, by Lemma \ref{lemma:linearLattice}.

The following statement, which is immediate from Theorem \ref{theorem:primclass}, will be important for both cases.

\begin{lemma} \label{lemma:Cgeneric}
Let $E$ be a minimal non-trivial $\emptyset$-definable equivalence relation in a homogeneous 3-dimensional permutation structure, and $C$ be an $E$-class. Then the induced structure on $C$ is generic, modulo the agreement of certain orders up to reversal.
\end{lemma}

We will frequently use the following characterization of genericity.

\begin{proposition}
Let $\Gamma$ be a homogeneous $n$-dimensional permutation structure. Then $\Gamma$ is generic iff for any non-empty open intervals $I_i$ in each order, $<_i$, $\cap_{i=1}^n I_i \neq \emptyset$.
\end{proposition}
\begin{proof}
Genericity of $\Gamma$ is equivalent to the following one-point extension property: given a type $p$ over a finite set $A$ not realized in $A$, $p$ is realized in $\Gamma$ iff its restriction to each individual order is realized by an element not in $A$. The restriction of $p$ to an order $<_i$ specifies a $<_i$-interval with endpoints in $A \cup \set{\pm \infty}$, which is open since $p$ is not realized in $A$. This interval is non-empty exactly when the restriction has a realization not in $A$.
\end{proof}

\subsection{Convex Closure}
In this section, we show $E$-classes are $<_1$-dense in their $<_1$-convex closures, and the $<_1$-convex closure of $E$ is an equivalence relation. Our arguments for this section depend heavily on the limited type structure in the case $k=3$, although in a few cases a step where our argument as given depends on $k=3$ could have been carried out in greater generality.

\begin{lemma} \label{lemma:norev}
Let $E$ be a minimal non-trivial $\emptyset$-definable equivalence relation in a homogeneous finite-dimensional permutation structure, and $C, C'$ be distinct $E$-classes. Then no 2-type $p$ is realized in both $C \times C'$ and $C' \times C$.
\end{lemma}
\begin{proof}
Let $a, b \in C$, $a', b' \in C'$, such that $a \ra{p} b'$ and $a' \ra{p} b$. Let $b \ra{q} b'$, and note that $p \neq q$, since otherwise transitivity would force $a' \ra{p} b'$ and so $p \subset E$. By homogeneity, there is an automorphism sending $(a, b')$ to $(a', b)$. Thus there must be some $c \in C$ such that $b' \ra q c$. But then by transitivity $b \ra q c$, which is a contradiction.
\end{proof}

\begin{definition}
We define $\widetilde E$ to be the $<_1$-convex closure of $E$, i.e. $a \widetilde E b$ if there exists a $c$ such that $aEc$ and $b$ is $<_1$-between $a$ and $c$. For an $E$-class $C$, we define $\widetilde C$ to be the $<_1$-convex closure of $C$.
\end{definition}

\begin{notation*}
For the rest of this section, we fix notation, by reversing and switching orders as needed, so that the 2-type $\ra{123}$ is contained in $E$, and if $E$ contains another 2-type besides $\ra{123}$ and its opposite then it contains the 2-type $\ra{23}$.
\end{notation*}

\begin{lemma}\label{lemma:3cross}
Let $E$ be a minimal non-trivial $\emptyset$-definable equivalence relation in a homogeneous 3-dimensional permutation structure, and $C$ be an $E$-class. Let $a_1, a_2 \in C$, $b \not\in C$, and $a_1 <_1 b <_1 a_2$. Then $tp(a_1, b) = \ra{12}$, $ tp(b, a_2) = \ra{13}$, or $tp(a_1, b) = \ra{13}, tp(b, a_2) = \ra{12}$.
\end{lemma}
\begin{proof}
If $E$ contains $\ra{123}$ and $\ra{23}$, the conclusion follows by Lemma \ref{lemma:norev} and the fact that only 4 2-types remain.

Otherwise, we have that $a_1 \ra{123} a_2$. Since we cannot have $a_1 \ra{123} b$, there is some $i$ such that $b <_i a_1 <_i a_2$, and so $b \ra{1i} a_2$. Thus, there is a unique $j$ such that $a_2 <_j b$, so $a_1 <_j b$. Thus $a_1 \ra{1j} b$ and $b \ra{1i} a_2$.
\end{proof}

\begin{corollary}\label{cor:contain}
Let $E$ be a minimal non-trivial $\emptyset$-definable equivalence relation in a homogeneous 3-dimensional permutation structure, and $C$ be an $E$-class. Suppose $b \in \widetilde C \bs C$. Then $b/E \subset \widetilde C$.
\end{corollary}
\begin{proof}
Take $a_1, a_2 \in C$ such that $a_1 <_1 b <_1 a_2$. By Lemma \ref{lemma:3cross}, we may suppose without loss of generality that $tp(a_1, b) = \ra{12}$, $tp(b, a_2) = \ra{13}$.

Take $b' \in b/E$. If $b' >_1 a_2$, then $b <_1 a_2 <_1 b'$ and $tp(b, a_2) = \ra{13}$, so by Lemma \ref{lemma:3cross}, $tp(a_2, b') = \ra{12}$. By homogeneity, there is an automorphism $\phi$ sending $(a_1, b)$ to $(a_2, b')$, so $b'$ is $<_1$-between $a_2$ and $\phi(a_2)$. The case where $b' <_1 a_1$ is nearly identical.
\end{proof}

\begin{corollary}\label{cor:samecross}
Let $E$ be a minimal non-trivial $\emptyset$-definable equivalence relation in a homogeneous 3-dimensional permutation structure, and $C$ be an $E$-class. Let $a, a' \in C$, with $a' \ra{123} a$. For any $b \not\in C$, if $b <_1 a'$ or $a <_1 b$, then $tp(a, b) = tp(a', b)$
\end{corollary}
\begin{proof}
We only treat the case $a <_1 b$, since the other case is similar.

Suppose $tp(a, b) = \ra{1x}$. By transitivity, $a' <_1 b$, $a' <_x b$. Since we cannot have $tp(a', b) = \ra{123}$, we are done. 

Now suppose $tp(b, a) = \ra{23}$. By transitivity, $a' <_1 b$. However, we cannot have $tp(a', b) = \ra{12}$ or $\ra{13}$, since by Lemma \ref{lemma:3cross} there would be some $a'' \in C$ such that $a'' >_1 b$, and then applying Lemma \ref{lemma:3cross} again, we would have that $tp(a, b)$ would also be $\ra{12}$ or $\ra{13}$.
\end{proof}

\begin{corollary} \label{cor:min}
Let $E$ be a minimal non-trivial $\emptyset$-definable equivalence relation in a homogeneous 3-dimensional permutation structure. Then any non-trivial $\emptyset$-definable equivalence relation contains $E$.
\end{corollary}
\begin{proof}
Consider the equivalence relation generated by a 2-type $p$, and without loss of generality assume $<_1$ holds in $p$. Given $a, b$ such that $a \ra{p} b$, find $b'$ such that $b \ra{123} b'$. By Corollary \ref{cor:samecross}, $a \ra{p} b'$, so $p$ generates $\ra{123}$. 

If $\ra{23} \subset E$, so $p = \ra{1x}$, then run the above argument with $b \ra{23} b'$. By transitivity, $a <_x b'$, so $a \ra{1x} b'$, and $\ra{23}$ is generated by $p$ as well.
\end{proof}

We note that much of the proof of the following lemma is concerned with ruling out a plausible configuration in which given $E$-classes $C, C_1$ such that $C_1 \subset \widetilde C$, then $C_1$ defines a non-trivial $<_1$-Dedekind cut in $C$. Although the type structure is too constrained to allow this with 3 orders, it seems possible that it may occur with more orders.

\begin{lemma} \label{lemma:dense}
Let $E$ be a minimal non-trivial $\emptyset$-definable equivalence relation in a homogeneous 3-dimensional permutation structure, and $C$ be an $E$-class. Then $C$ is $<_1$-dense in $\widetilde C$.
\end{lemma}
\begin{proof}
Given an $E$-class $C$ and an element $a$, we denote by $\widehat a(C)$ the $<_1$-Dedekind cut defined by $a$ in $C$.

If $C$ is not $<_1$-dense in $\widetilde C$, then there are $a, b \in \widetilde C$ such that $\widehat a(C) = \widehat b(C)$.

Using the next two claims, we show that we may suppose $a/E = b/E$ and $a \ra{123} b$.

\begin{claim}
 $a,b \in \widetilde C \bs C$.
\end{claim}
\begin{claimproof}
 Trivially, we cannot have $a, b \in C$. Now assume only one of $a, b \in C$, say $a$. Then $a$ is a maximal element of the cut $\widehat b(C)$. But given any $d \in \widetilde C\bs C$, $\widehat d(C)$ has no maximal or minimal elements; otherwise, the elements of $C$ would realize at least 3 types over $d$, but there are only 2 realized types by Lemma \ref{lemma:3cross}. Thus $a, b \in \widetilde C \bs C$.
\end{claimproof}

Now suppose $a/E = b/E$, but $a \la{23} b$. By the genericity of $C$, there is a $b'$ in the $<_1$-interval $(a, b)_{<_1}$ such that $a \ra{123} b'$, so we may replace $b$ by $b'$.

\begin{claim}
Suppose $C_1 = a/E \neq b/E = C_2$. Then there exists $a' \in C_1$ such that $a' \ra{123} a$ and $\widehat a(C) = \widehat {a'}(C)$
\end{claim}
\begin{claimproof}[Proof of Claim 1.]
Let $a' \ra{123} a$. Since $a <_1 b$, by Corollary \ref{cor:samecross} $tp(a, b) = tp(a', b)$. Since by Lemma \ref{cor:contain}, $a/E \in \widetilde C$, there is a $c \in C$ such that $c <_1 a'$, so by Corollary \ref{cor:samecross} $tp(a, c) = tp(a', c)$. Thus $(a, b, c) \cong (a', b, c)$, so by homogeneity there is an automorphism fixing $c$ and taking $(a, b)$ to $(a', b)$. Thus $\widehat {a'}(C) = \widehat b(C) = \widehat a(C)$.
\end{claimproof}

In this case, we may then replace $a, b$ by $a', a$. 

Thus, we may now suppose that $a/E = b/E = C_1 \neq C$ and $a \ra{123} b$.

\begin{claim}
$\widehat c(C)$ is independent of the choice of $c \in C_1$.
\end{claim}
\begin{claimproof}[Proof of Claim 2.]
Consider $x, y \in C_1$, and, using the genericity of $C_1$, find $c_1, c_2 \in C_1$ such that $c_1 <_1 x,y <_1 c_2$ and $c_1 \ra{123} c_2$.

Take $z \in C$, with $z <_1 a,c_1$. By Corollary \ref{cor:samecross}, $(z, a, b) \cong (z, c_1, c_2)$. Thus, since $\widehat a(C) = \widehat b(C)$, we have $\widehat c_1(C) = \widehat c_2(C)$, so $\widehat x(C) = \widehat y(C)$.
\end{claimproof}

Without loss of generality, we now assume $C <_2 C_1$, so by Lemma \ref{lemma:3cross}, the types realized in $C \times C_1$ are $\ra{12}$ and $\la{13}$. Thus by homogeneity, given any $E$-classes $C, C'$, if $\ra{12}$ or $\la{13}$ is realized in $C \times C'$, then $C'$ defines a $<_1$-Dedekind cut in $C$; if neither these types nor their opposites are realized, then the only remaining types are $\ra{23}$ and $\la{23}$, and by Lemma \ref{lemma:norev} exactly one of them is realized, so neither class is in the $<_1$-convex closure of the other. In particular, $E$-classes are $<_2, <_3$-convex.

Note that if every $E$-class $C' \subset \widetilde C$ such that $C <_2 C'$ defined the same $<_1$-Dedekind cut in $C$, $C$ would have an $\emptyset$-definable partition, contradicting the minimality of $E$.

\begin{claim}
Both factors of the $(12,23,13)$-majority diagram, displayed below (with the edge $x \ra{23} z$ not drawn), are realized in $\Gamma$.

\begin{figure}[h]
\begin{diagram}
& & \overset{x}{\bullet} & & \\
& \ruTo^{12} & \dTo_{23}  & \rdTo^{23} & \\
a \odot &\rTo^{12} & \bullet_{y}  & \rTo^{13} & \odot b \\
& \rdTo_{23} & \dTo_{23} & \ruTo_{13} & \\
& & \underset{z}{\bullet} & & \\
\end{diagram}
\caption{}
\label{fig:15}
\end{figure}
\end{claim}
\begin{claimproof}[Proof of Claim 3.]
We only prove the first factor is realized, since the argument for the second is nearly identical. First, as shown below, the first factor is the unique amalgam of the following 3-types, so it suffices to show these are realized.

\begin{figure}[h]
\begin{diagram}
& & \overset{x}{\odot} & & \\
& \ruTo^{12} & \dTo_{23}   \\
a \bullet &\rTo^{12} & \bullet_{y}  \\
& \rdTo_{23} & \dTo_{23}\\
& & \underset{z}{\odot} & & \\
\end{diagram}
\caption{}
\label{fig:16}
\end{figure}

For the triangle $(a,x,y)$ from the diagram, let $a/E = C$.  Take distinct $E$-classes $C', C'' \subset \widetilde C$ such that $C <_2 C'<_2 C''$ and $C'$ and $C''$ define distinct $<_1$-cuts in $C$. Then there are $x \in C', y \in C''$ realizing the triangle $(a,x,y)$.

For the triangle $(a,y,z)$, we will show it is the unique amalgam of the following diagram.

\begin{figure}[h]
\begin{diagram}
 y \odot & &  & & \odot z  \\
& \luTo_{12} &  & \ruTo_{23} & \\
& & {\bullet}_a & & \\
\end{diagram}
\caption{}
\label{fig:17}
\end{figure}

By transitivity, $y <_3 z$ and $z <_1 y$, so the possible completions are $y \ra{23} z$ and $z \ra{12} y$. However, if $z \ra{12} y$, then $y/E$ defines a $<_1$-Dedekind cut in both $a/E$ and $z/E$, but $z/E <_1 a/E$, which is a contradiction. Thus the only allowed completion is $y \ra{23} z$.
\end{claimproof}

We are forced to complete the $(12,23,13)$-majority diagram by $a \ra{123} b$, so that $aEb$. However, $a \ra{12} x \ra{23} b$ violates the requirement that $E$-classes are $<_2$-convex. Thus $\Gamma$ is not homogeneous.
\end{proof}

\begin{proposition}
Let $E$ be a minimal non-trivial $\emptyset$-definable equivalence relation in a homogeneous 3-dimensional permutation structure, and $C$ be an $E$-class. Then $\widetilde E$ is an equivalence relation.
\end{proposition}
\begin{proof}

Let $C$ be an $E$-class. By Corollary \ref{cor:contain}, $\widetilde C$ is a union of $E$-classes. Now suppose $C' \subset \widetilde C$ is an $E$-class. By Lemma \ref{lemma:dense} there are $c_1, c_2 \in C'$ such that $\widehat c_1(C) \neq \widehat c_2(C)$, and applying Corollary \ref{cor:contain} again we see $C \subset \widetilde {C'}$. Thus $\widetilde E$ defines a partition.
\end{proof}

\begin{corollary} \label{cor:congconv}
$\widetilde E$ is a congruence, $E$-classes are $(<_2,<_3)$-convex, $<_2 = <_3$ on $E$-classes, and $<_2 = <_3^{opp}$ between $E$-classes in the same $\widetilde E$-class.
\end{corollary}

\subsection{Reduction via Quotients}
Since $\widetilde E$ is a congruence by Corollary \ref{cor:congconv}, it suffices to consider the case $\widetilde E = \bbone$, since we may otherwise consider the restriction ${\Gamma \upharpoonright \widetilde E}$. For this subsection, we work with $k$-dimensional permutation structures.

We now aim for the following two lemmas. The first implies that $\Gamma$ is determined by its restriction to $E$-classes and the reduct of $\Gamma/E$ forgetting all orders that are not $E$-convex. The second allows us to carry out our induction by showing that the specified reduct of $\Gamma/E$ must be homogeneous.

The first of the following lemmas is more naturally stated in the language of subquotient orders, but as it is the concluding step in the classification of certain permutation structures, we give it in a form appropriate for its intended application.

\begin{lemma} \label{lemma:uniqueextension}
Let $(\Gamma^*, <_1^*, ..., <_\ell^*)$ be homogeneous.Let $k \geq \ell$, and partition $[k]$ as $\cup_{i \leq m} I_i$, such that each $I_i$ contains at most one $j \geq \ell+1$.
 Then there exists a homogeneous structure $(\Gamma, E, <_1, ..., <_k)$, unique up to isomorphism, with the following properties.
\begin{enumerate}
\item $E$-classes are $<_1, ..., <_\ell$-convex and $<_{\ell+1}, ..., <_k$-dense.
\item $(\Gamma/E, <_1, ..., <_\ell) \cong (\Gamma^*, <_1^*, ..., <_\ell^*)$
\item $<_j \upharpoonright_E = <_{j'}\upharpoonright_E$ for $j, j'$ in a given $I_i$, and the induced structure on any $E$-class $C$ is fully generic, modulo the identification of orders in the same $I_i$.
\end{enumerate}
\end{lemma}


\begin{lemma} \label{lemma:homquotient}
Let $\Gamma$ be a homogeneous $k$-dimensional permutation structure.  Let $E$ be a minimal non-trivial $\emptyset$-definable equivalence relation in $\Gamma$, and suppose $E$-classes are $<_i$-convex for $1 \leq i \leq \ell$ and $<_i$-dense for $\ell+1 \leq i \leq k$. Suppose each $E$-class is generic, modulo the agreement of certain orders up to reversal. Then $(\Gamma/E, <_{1}, ..., <_\ell)$ is homogeneous.
\end{lemma}

The following lemmas prepare for the proof of Lemmas \ref{lemma:uniqueextension} and \ref{lemma:homquotient}. The first of these is not necessary for the case $k=3$, since there $E$ is only dense with respect to one order.

\begin{lemma}\label{lemma:denserestrict}
Suppose $(\Gamma, <_1, ..., <_k)$ is homogeneous. Let $E$ be a minimal non-trivial $\emptyset$-definable equivalence relation in $\Gamma$, and $C$ be an $E$-class. Suppose $C$ is generic, modulo the agreement of certain orders up to reversal. Further suppose that $C$ is $<_i$-convex for $1 \leq i \leq \ell$ and $<_i$-dense for $\ell+1 \leq i \leq k$. Then
\begin{enumerate}
\item If $C_1, C_2$ are $E$-classes, then $C_1$ remains homogeneous after naming $C_2$.
\item If $i, j \geq \ell+1$, and $<_i \upharpoonright_E = <_j \upharpoonright_E$, then $<_i = <_j$.
\end{enumerate}
\end{lemma}
\begin{proof} \hspace{.1 cm} \newline
$(1)$ Given a finite $A \subset C_1$ and $i \geq \ell+1$, let $B_i = \set{x \in C_2 | A <_i x}$. Each such $B_i$ is a $<_i$-terminal segment of $C_2$, so by genericity their intersection is non-empty.

Now, consider $A_1 \cong A_2$ finite substructures of $C_1$. Let $A = A_1 \cup A_2$, and choose a $b$ in the intersection of the corresponding $B_i$. Then $A_1b \cong A_2b$, and by homogeneity there is an automorphism taking $A_1b$ to $A_2b$ and fixing $b$, hence $C_2$.

$(2)$ Suppose the condition is false, as witnessed by $<_i, <_j$. We consider $E$-classes as ordered sets with respect to the common restriction of these orders.

Take $a, b$ with $a <_i b$ and $b <_j a$, and  let $C_1 = a/E$ and $C_2 = b/E$. Let
$$I_a = \set{x \in C_2| a<_i x, x<_j a}, J_a = \set{x \in C_1| I_a \cap I_x \neq \emptyset}$$
Note that these sets are intervals in $C_2$ and $C_1$, respectively.

\begin{claim*}
$J_a = \set{a}$
\end{claim*}
\begin{claimproof}
By density and genericity, there are $b_1, b_2 \in C_2$ such that $b_1 <_{i,j} a <_{i, j} b_2$, so $I_a \subset (b_1, b_2)$. Then find $a_1, a_2 \in C_1$ such that $a_1 <_{i, j} b_1, b_2 <_{i, j} a_2$, so $J_a \subset (a_1, a_2)$.

Thus $J_a$ is $(a, C_2)$-definable and $<_{i, j}$-bounded. By $(1)$, $J_a$ is $a$-definable in $C_1$ and $<_{i, j}$-bounded, so $J_a = \set{a}$ by genericity of $C_1$.
\end{claimproof}

If there were some $b' \in I_a$ with $b' \neq b$, then by density, we could find some $a' \in C_1$ $<_i$-between $b$ and $b'$, and so would have $a' \in J_a$. Thus $I_a = \set{b}$. But by density there is a $b' \in C_2$ $<_i$-between $a$ and $b$, so $b' \in I_a$, which is a contradiction.
\end{proof}

Given $(\Gamma, <_1, ..., <_k)$ homogeneous such that no orders agree up to reversal, with $E$-classes $<_i$-convex for $1 \leq i \leq \ell$ and $<_i$-dense for $\ell+1 \leq i \leq k$, we will prefer to work in the quantifier-free interdefinable reduct $\Gamma^{red} = (\Gamma, {<_{i_1}',} ..., {<_{i_m}',} {<_1'',} ...,{<_\ell'',} {<_{\ell+1},} ..., <_k)$ obtained as follows.
\begin{enumerate}
\item For $1 \leq i \leq \ell$, decompose $<_i$ into two subquotient orders: $<_i'$ from $\bbzero$ to $E$ and $<_i''$ from $E$ to $\bbone$. 
\item For each $i \geq \ell+1$, add the restriction $<_i \upharpoonright_E$ to the language as $<_i'$.
\item Consider the set of all $<_i'$. Many of these subquotient orders may be equal up to reversal, so pick one representative from each class and forget the rest. By Lemma \ref{lemma:denserestrict}, each class can contain at most one $<_i'$ with $i \geq \ell+1$, in which case this is taken as the representative.
\item Forget the $<_i'$ for $i \geq \ell+1$. 
\end{enumerate}

We now prove a 1-point extension property, which shows that to realize a type $p$ in an $E$-class $C$ , it is sufficient that the restriction of the type to each subquotient order is individually realized. 

\begin{lemma}\label{lemma:1point}
Let $(\Gamma, <_1, ..., <_k)$ be homogeneous such that no orders agree up to reversal. Let $E$ be a minimal non-trivial $\emptyset$-definable equivalence relation in $\Gamma$, and $C$ be an $E$-class. Suppose the induced structre on $C$ is generic, modulo the agreement of certain orders up to reversal. Suppose that $C$ is $<_i$-convex for $1 \leq i \leq \ell$ and $<_i$-dense for $\ell+1 \leq i \leq k$. We now work in $\Gamma^{red}$.

Let $A \subset \Gamma^{red}$ be finite, and $p$ a 1-type over $A$ not realized in $A$. Then $p$ is realized in a given $E$-class $C$ by a point not in $A$ iff the following hold.
\begin{enumerate}
\item $p \upharpoonright (<_1'', ..., <_\ell'')$ is realized by $C$ in $\Gamma/E$.
\item For each $<_i'$,  $(p \upharpoonright <_i') \upharpoonright A$ is realized in $C \bs A$.
\item For $j \geq \ell+1$, $p \upharpoonright <_j$ is realized by some element not in $A$.
\item $p$ does not contain the formula ``$x = a$'' for any $a \in A$.
\end{enumerate}
\end{lemma}
\begin{proof}
These conditions are clearly necessary. We will prove they suffice. By condition $(1)$, all of $C$ satisfies $p\upharpoonright (<_1'', ..., <_\ell'')$. List all the subquotient orders from $\bbzero$ to $E$ together with $<_i$ for $i \geq \ell+1$ as $<_1^*, ..., <_n^*$, and let $p_i = p \upharpoonright <_i^*$. It now suffices to show $p_i$ contains a non-empty open $<_i^*$-interval of $C$, since then by the genericity of $C$ there will be some point in their intersection, which thus realizes $p$. 

In the case $<_i^*$ is a subquotient order from $\bbzero$ to $E$, by condition $(2)$ some point in $C$ realizes $p_i$ restricted to parameters outside of $C$, and so all of $C$ does; again by condition $(2)$, $p_i$ restricted to parameters inside of $C$ then contains an open interval of $C$. In the case $<_i^* = <_j$ for $j \geq \ell+1$, condition $(3)$ implies $p_i$ contains a non-empty open interval in $\Gamma$; since $E$-classes are $<_i^*$-dense, this interval meets $C$ in a non-empty open interval.
\end{proof}

\begin{proof}[Proof of Lemma \ref{lemma:homquotient}]
Let $\bar A \cong \bar B$ be finite subsets of $(\Gamma/E, <_1, ..., <_\ell)$. We lift $\bar A$ to $A \subset (\Gamma, <_1, ..., <_k)$, and look for an automorphism moving $A$ to a set covering $B$.

We proceed by induction on $|\bar A|$, and so consider $A = A_0 \cup \set{a}$ with $\bar a \not\in \bar A_0$, $\bar B = \bar A_0 \cup \set{C}$ for some $E$-class $C \not\in \bar A_0$.

Let $p = tp(a/A_0)$. We will now work in $\Gamma^{red}$ and use Lemma \ref{lemma:1point} to find a realization of $p$ in $C$. Condition $(1)$ is equivalent to $\bar A \cong \bar B$. Since $\bar a \not\in \bar A_0$, $A_0 \cap C = \emptyset$, so $(p \upharpoonright <_i') \upharpoonright A$ simply says $x$ is not $<_i'$-related to any $a \in A$, which will be true for every $x \in C$. Finally, since $\bar a \not\in \bar A_0$, $a \not\in A_0$, so $a$ witnesses condition $(3)$.
\end{proof}

\begin{proof}[Proof of Lemma \ref{lemma:uniqueextension}]
For existence, let $\Gamma$ be the composition $\Gamma^*[C]$, where $C$ only carries the equality relation, and let $E$ be the corresponding equivalence relation. Note that each $<_i^*$ is now a subquotient order from $E$ to $\bbone$. For $1 \leq i \leq m-(k-\ell)$, add a generic subquotient order $<'_i$ from $\bbzero$ to $E$. For $\ell+1 \leq i \leq k$, add a generic linear order $<_i$. We may then define the specified convex orders $<_i$ for $1 \leq i \leq \ell$ as compositions of the $<_i^*$ with the $<_j'$ or the restrictions to $E$ of the $<_n$ for $\ell+1 \leq n \leq k$. 

For uniqueness, suppose we have a structure $(\Gamma', <_1, ..., <_k)$ satisfying the conditions. We will show $\Gamma'^{red}$ has the same finite substructures as the $\Gamma^{red}$ we constructed above; as they are both homogeneous, they will thus be isomorphic.

As all the subquotient orders added to construct $\Gamma^{red}$ were added generically, every finite substructure of $\Gamma'^{red}$ is a substructure of $\Gamma^{red}$. For the converse, we proceed by induction on the size of the substructure. Let $A \cup \set{a}$ be a finite substructure of $\Gamma^{red}$, such that $A$ is a substructure of $\Gamma'^{red}$. We will use Lemma \ref{lemma:1point} to show $p = tp(a/A)$ is realized in $\Gamma'^{red}$.

We may assume $a \not\in A$, otherwise we are done, so condition $(4)$ is satisfied. As (suitable reducts of) $\Gamma^{red}/E$ and $\Gamma'^{red}/E$ both are isomorphic to $\Gamma^*$, and as $a/E$ realizes $p \upharpoonright (<_1'', ..., <_\ell'')$ in the former, there is some $E$-class $C$ realizing it in the latter, so condition $(1)$ is satisfied. For condition $(2)$,  again since the quotient structures are isomorphic, we may pick $C$ such that for each $b \in A$, $C = b/E$ iff $a/E = b/E$. Thus, we are only concerned about $(p \upharpoonright <_i') \upharpoonright (A \cap C)$; but as this restricted type doesn't violate transitivity, it is realized in $C$ since $<_i'$ is dense on $C$. Finally for condition $(3)$, we again have that $p \upharpoonright <_j$ doesn't violate transitivity, and so is realized by some element not in $A$ since $<_j$ is dense on $\Gamma'^{red}$.
\end{proof}



\begin{remark}
Lemma \ref{lemma:uniqueextension} is also true if $(3)$ is relaxed to allow certain restrictions to be the reversals of others. The only case that isn't immediate is if we require $<_i \upharpoonright_E = (<_j \upharpoonright_E)^{opp}$ for $<_i, <_j$ dense. But then $<_i = <_j^{opp}$ by Lemma \ref{lemma:denserestrict}.
\end{remark}

\subsection{The Imprimitive Classification}
We now classify the imprimitive homogeneous structures in a language of 3 linear orders. We present the structures up to definable equivalence, and do so in a language of subquotient orders, each of which is generic, and equivalence relations. Presenting these in the language requires picking a 3-tuple of linear orders interdefinable with the presentation we give, and leads to a mob of examples. The number of $\emptyset$-definable linear orders may be substantial, and the number of suitable subsets quite a bit larger.


We first classify the imprimitive homogeneous 3-dimensional permutation structures $(\Gamma, E, <_1, <_2, <_3)$ in which $\widetilde E = \bbone$, so $\Gamma$ has no non-trivial $\emptyset$-definable congruence. By Corollary \ref{cor:congconv} and Lemmas \ref{lemma:homquotient} and \ref{lemma:uniqueextension}, $\Gamma$ is determined by $(\Gamma/E, <_2)$ and $(\Gamma \upharpoonright_E, <_1, <_2)$, which are themselves primitive homogeneous. There are thus two possibilities.

\begin{enumerate}
\item ($<_1 \upharpoonright_E \neq <_2 \upharpoonright_E$) $\Gamma$ may be presented as $(\Gamma, E, (<_i')_{i=1}^3)$ with $<_1'$ from $\bbzero$ to $\bbone$, $<_2'$ from $\bbzero$ to $E$, and $<_3'$ from $E$ to $\bbone$.
\item ($<_1 \upharpoonright_E = <_2 \upharpoonright_E$) $\Gamma$ may be presented as $(\Gamma, E, (<_i')_{i=1}^2)$ with $<_1'$ from $\bbzero$ to $\bbone$ and $<_2'$ from $E$ to $\bbone$.
\end{enumerate}
$(1)$ is just $\Q^2_{lex}$ with an additional generic order and $(2)$ is the structure described in Example \ref{ex:complex3} in Section 2. 

Also note that when presented in the language of 3 linear orders, $(1)$ uses all 8 2-types, while $(2)$ only uses 6 of them. Thus $(1)$ cannot appear as a factor in a composition, while $(2)$ can.

If $\Gamma$ has a non-trivial $\emptyset$-definable congruence, then it is a composition, whose factors are either primitive or one of the above structures. Below, let $\Gamma^{(g)}_i$ to denote the generic $i$-dimensional permutation structure

If all of the factors are primitive, then each factor is interdefinable with $\Gamma^{(g)}_i$ for $i \in \set{1, 2}$. Each such factor contributes $2^i$ 2-types. As there are at most 8 2-types available, we get at most the following structures.

\begin{enumerate}
\item[(3)] For any multisubset $I \subset \set{1, 2}$ such that $|I|>1$ and $\sum_{i \in I} 2^i \leq 8$, $\Gamma$ is the composition in any order of $\Gamma^{(g)}_i$ for $i \in I$.
\end{enumerate}

Finally, if one of the factors is imprimitive, we noted earlier it must be $(2)$. There are only 2 2-types remaining, so the other factor must be $\Gamma^{(g)}_1$.

\begin{enumerate}
\item[(4)] Let $\Gamma^*$ be the structure from $(2)$. Then $\Gamma = \Gamma^*[\Gamma^{(g)}_1]$ or $\Gamma^{(g)}_1[\Gamma^*]$.
\end{enumerate}

For all of these structures we have only shown that at most 8 2-types are realized, but it is easy to check that each structure can be presented in a language of 3 linear orders by taking restrictions and compositions of the subquotient orders, which concludes our derivation of the catalog.

 This last step prompts the following special case of Question \ref{qu:represent}.

\begin{question}
Let $\Gamma$ be a finite-dimensional permutation structure, with a linear lattice of $\emptyset$-definable equivalence relations. If $\Gamma$ has at most $2^k$ non-trivial 2-types, can $\Gamma$ be presented as a $k$-dimensional permutation structure?
\end{question}

We remark that the linearity hypothesis is necessary, since the full product $\Q^2$ (see Example \ref{ex:fullQ2}) only has 8 non-trivial 2-types, but requires 4 linear orders.

\chapter{Questions around the Catalog}\label{chap:questions}
We collect here the various questions and conjectures that have appeared in the preceding chapters, or that are suggested by them.

The results of this thesis point to the following as a reasonable conjecture.

\begin{conjecture} \label{conj:probclassification}
   Every homogeneous finite-dimensional permutation structure with lattice of $\emptyset$-definable equivalence relations isomorphic to $\Lambda$ is interdefinable with the \fraisse limit of some well-equipped lift of the class of all finite $\Lambda$-ultrametric spaces.
\end{conjecture}

In particular, we state the primitive case of our conjecture separately. 

\begin{conjecture}[Primitivity Conjecture]
Every primitive homogeneous finite dimensional permutation structure can be constructed by the following procedure.
\begin{enumerate}
\item Identify certain orders, up to reversal.
\item Take the \fraisse limit of the resulting amalgamation class, getting a fully generic structure, possibly in a simpler language.
\end{enumerate}
\end{conjecture}

Again, we note that Pierre Simon has confirmed a proof of the Primitivity Conjecture in personal communication just before the submission of this thesis \cite{Simon}. However, we retain the following questions, as it still seems worth considering them via direct amalgamation arguments.

As in Chapter \ref{c:3dim}, we may split the Primitivity Conjecture into two parts, corresponding to Lemmas \ref{lemma:trianglereduce} and \ref{lemma:primtarget}.

\begin{conjecture} \label{conj:generalTriangleReduce}
Let $\Gamma$ be a homogeneous finite-dimensional permutation structure realizing all 3-types. Then $\Gamma$ is fully generic.
\end{conjecture}

\begin{conjecture}
Let $\Gamma$ be a primitive homogeneous finite-dimensional permutation structure. Then all 3-types involving realized 2-types are realized.
\end{conjecture}

Lemma \ref{lemma:4gen} provides a partial result in the direction of Conjecture \ref{conj:generalTriangleReduce}, but the size of the types to be considered goes to infinity with the number of linear orders. This prompts the following question.

\begin{question}
Is there some $n$ such that if any homogeneous finite-dimensional permutation structure contains all $n$-types, then it is fully generic? In particular, can this be proven using some variation of Lachlan's Ramsey argument?
\end{question}

A concrete test of Conjecture \ref{conj:generalTriangleReduce} would be to try to repeat the the proof of Lemma \ref{lemma:trianglereduce} given Lemma \ref{lemma:4gen} in higher dimensions. In particular, up until $k = 14$ linear orders, this requires constructing 5-point configurations as unique amalgams of smaller configurations.

These conjectures are put forward as an answer to the following question of Cameron, which was our starting point.

\begin{problem}[Cameron]
Classify, for each $n$, the homogeneous $n$-dimensional permutation structures.
\end{problem}

A major qualitative consequence of our conjectural classification is the Distributivity Conjecture.

\begin{conjecture} [Distributivity Conjecture]
A finite lattice is isomorphic to the lattice of $\emptyset$-definable equivalence relations in some homogeneous finite-dimensional permutation structure iff it is distributive.
\end{conjecture}

One could try to prove this directly by showing an expansion from the metric language to the full language of linear orders of Figure \ref{fig:6} can be realized in any homogeneous finite-dimensional permutation structure. The amalgamation argument could then be carried out in the full language instead.

In particular, a contradiction to the Distributivity Conjecture would have to arise as follows. Let $\Gamma$ be a homogeneous finite-dimensional permutation structure with a non-distributive lattice of $\emptyset$-definable equivalence relations, and $\Gamma^{met}$ its reduct to the metric language. We know $\Gamma^{met}$ must be non-homogeneous. As $\Lambda$ satisfies the infinite-index property, we may construct each 4-point factor of the diagram from Lemma \ref{lemma:IIPDIS} in $\Gamma^{met}$. Each such factor expands to a substructure of $\Gamma$, but the factors cannot be expanded in such a way that the expansions agree on the base of the diagram.

The most promising counterexample to Conjecture \ref{conj:probclassification} seems to be the configuration appearing in Lemma \ref{lemma:dense}, which there was proven untenable with 3 orders.

\begin{question} \label{question:exotic}
Can a configuration as in Lemma \ref{lemma:dense}, in which $E$-classes define non-trivial Dedekind cuts in other $E$-classes, appear in a homogeneous finite-dimensional permutation structure?
\end{question}

One possible path to ruling out such counterexamples would be to prove the following conjecture. This was proven for the case of 3 linear orders in Lemma \ref{lemma:denserestrict}.

\begin{conjecture}
Let $E$ be a $\emptyset$-definable equivalence relation in a finite-dimensional permutation structure, and let $C_1, C_2$ be $E$-classes. Then $C_2$ remains homogeneous after naming $C_1$. 
\end{conjecture}

We now give the most immediate test case for Question \ref{question:exotic}.

\begin{question}
Assume the Primitivity Conjecture. Let $\Gamma$ be a homogeneous 4-dimensional permutation structure with a linear lattice of $\emptyset$-definable equivalence relations. Can $\Gamma$ realize a configuration as in Question \ref{question:exotic}?
\end{question}

Conjecture \ref{conj:probclassification} presents a conjectural classification of the homogeneous finite-dimensional permutation structures, but not a classification of the homogeneous $k$-dimensional permutation structures for each $k$. This is due to the issues in determining the number of linear orders needed to represent a homogeneous finite-dimensional permutation structure produced by our construction, and sets aside the following combinatorial problem, of significant interest in its own right.

\begin{question} 
Given a homogeneous finite-dimensional permutation structure $\Gamma$ presented in a language of equivalence relations and subquotient orders, what is the minimal $n$ such that $\Gamma$ is quantifier-free interdefinable with an $n$-dimensional permutation structure?
\end{question} 

We have addressed the following related question more directly in Corollary \ref{Corollary:CountLattice}. 

\begin{question}
Given a lattice $\Lambda$, what is the minimal $n$ such that $\Lambda$ is isomorphic to the lattice of $\emptyset$-definable equivalence relations of some homogeneous $n$-dimensional permutation structure? 

In particular, is the bound of Corollary \ref{Corollary:CountLattice} optimal?
\end{question}

Finally, we again mention a general question regarding homogeneous ordered structures, which this thesis may be seen as investigating a simple case of.

\begin{question}
Is every homogeneous ordered structure interdefinable with an expansion of a \emph{homogeneous} proper reduct by a linear order?

Is every \textbf{primitive} homogeneous ordered structure that is the \fraisse limit of a strong amalgamation class interdefinable with an expansion of a homogeneous proper reduct by a \textbf{generic} linear order?
\end{question}






\chapter{Homogeneous Finite-Dimensional Permutation Structures as Ramsey Expansions of $\Lambda$-Ultrametric Spaces}\label{chap:Ramsey}
\section{Introduction}
We now turn to the study of the dynamical properties of the automorphism groups of homogeneous finite-dimensional permutation structures, in the manner of \cite{KPT}, meaning that we will prove a suitable structural Ramsey theorem for the amalgamation classes associated to all the homogeneous finite-dimensional permutation structures in our catalog.

Our Ramsey theorem for well-equipped lifts of $\Lambda$-ultrametric spaces is proved using the results of \cite{HN}, which provide a black box for complicated arguments involving the partite construction. The theorem of \cite{HN} we use is particularly useful for proving a Ramsey theorem for a class whose Fra\"\i ss\'e limit has a non-degenerate algebraic closure operation, in the sense that there exist sets which are not their own algebraic closure. This allows us to make use of a fragment of $M^{eq}$, where the algebraic operations on $\Lambda$ translate into a non-trivial algebraic closure operator.

In particular, we use a combination and generalization of the encoding techniques used in \cite{HN} to prove Ramsey theorems for the free product of Ramsey classes, which involves duplicating the structure and adding unary functions representing a bijection, and for structures that have a chain of definable equivalence relations, which involves adding elements representing classes of $\emptyset$-definable equivalence relations. 
   
A point of considerable technical interest is this: since we are dealing with an arbitrary finite distributive lattice $\Lambda$ of equivalence relations rather than just a chain of such, the algebraic closure operation in the structures we consider is non-unary, that is the algebraic closure of a set is not determined by the algebraic closures of the elements of the set. A non-unary algebraic closure significantly complicates applying the theorems of \cite{HN}, and consequently few classes with a non-unary algebraic closure have been proven to be Ramsey classes, although some new examples were recently given in \cite{HN2}.
   
The main point in the analysis of this closure operator is to show that, in an appropriate category, the closure of a finite set is finite. Rather than analyzing the closure operation explicitly, we derive this from the algebraic closure operation on imaginary elements in the generic $\Lambda$-ultrametric space (as in Lemma \ref{lemma:closure0}).

Another point of considerable technical interest is our use of what we call a \textit{quantifier-free reinterpretation}, a generalization of the argument appearing in \cite{Bod}*{Section 4}, to transfer the Ramsey property between classes. The natural class our arguments would be carried out in has a linear order satisfying many constraints, and the reinterpretation technique allows us to argue in a class where the linear order is more nearly generic, thereby avoiding much bookkeeping.

\section{Multi-Amalgamation Classes and Ramsey Theorems}
We now give an exposition of Theorem 2.2 in \cite{HN}, which provides sufficient conditions for proving a subclass of a known Ramsey class is Ramsey.

\begin{definition}
An $L$-structure $A$ is \textit{irreducible} if for every distinct $x, y \in A$, there is some $R \in L$ and some tuple $\vec t$ containing $x, y$ such that $R(\vec t)$ holds in $A$.

A homomorphism $f: A \to B$ is a \textit{homomorphism-embedding} if $f$ restricted to any irreducible substructure of $A$ is an embedding, i.e. the restriction of $f$ is injective and for any $R \in L$, $R(x_1, ..., x_{r}) \Leftrightarrow R(f(x_1), ..., f(x_{r}))$, where $r$ is the arity of $R$.

Given an $L$-structure $C$ and a class $\KK$ of $L$-structures, we say $C'$ is a \textit{$\KK$-completion} of $C$ if $C' \in \KK$ and there is a homomorphism-embedding $f: C \to C'$.

Given an $L$-structure $C$, an irreducible $B \subset C$, and a class $\KK$ of $L$-structures, we say $C'$ is a \textit{$\KK$-completion of $C$ with respect to copies of $B$} if $C'$ is an irreducible $\KK$-structure and there is a function $f:C \to C'$ such that $f$ restricted to any $\hat B \in \binom C B$ is an embedding.
\end{definition}

\begin{definition}
Given a language $L$, a \textit{closure description} $\UU$ is a set of pairs $(R^U, B)$, where $R^U \in L$ is an $n$-ary relation, and $B$ is a non-empty irreducible $L$-structure on the set $\set{1, ..., m}$ for some $m \leq n$. We call $R^U$ a \textit{closure relation}, and the corresponding structure $B$ its \textit{root}.
\end{definition}

\begin{definition}
Given an $L$-structure $A$, an $n$-ary relation $R \in L$, and $k \leq n$, the \textit{$R$-out-degree} of a $k$-tuple $(x_1, ..., x_k) \in A^k$ is the number of tuples $(x_{k+1}, ..., x_n) \in A^{n-k}$ such that $R(x_1, ..., x_n)$ holds in $A$.  

Given a closure description $\UU$, we say that a structure $A$ is \textit{$\UU$-closed} if, for every $(R^U, B) \in \UU$, the $R^U$-out-degree of any tuple $\vec t$ of elements of $A$ is 1 if $\vec t$ is an embedding of $B$ into $A$, and 0 otherwise.
\end{definition}

Thus, in a $\UU$-closed structure, a closure relation can be thought of as a function assigning additional points to each copy of its root. The strong amalgamation condition in the following definition ensures that, in our applications, the closure relations are such that these functions generate the algebraic closure in the Fra\"\i ss\'e limit.

We are now ready for the main definition and theorem.

\begin{definition}
Let $\RR$ be a Ramsey class of finite irreducible $L$-structures, and let $\UU$ be a closure description in $L$. We say that a subclass $\KK \subset \RR$ is an \textit{$(\RR, \UU)$-multi-amalgamation class} if:
\begin{enumerate}
\item $\KK$ consists of finite $\UU$-closed $L$-structures.
\item $\KK$ is closed under taking $\UU$-closed substructures.
\item $\KK$ has strong amalgamation.
\item \textbf{Locally finite completion property}: Let $B \in \KK$ and $C_0 \in \RR$. Then there exists an $n=n(B, C_0)$ such that for any $\UU$-closed $L$-structure $C$ that satisfies the conditions below, there exists a structure $C'$ that is a $\KK$-completion of $C$ with respect to copies of $B$. The conditions required on $C$ are as follows.
\begin{enumerate}
\item $C_0$ is an $\RR$-completion of $C$.
\item Every substructure of $C$ with at most $n$ elements has a $\KK$-completion.
\end{enumerate}   
\end{enumerate} 
\end{definition}

\begin{theorem} [\cite{HN}*{Theorem 2.2}] \label{theorem:multiamalg}
Let $\RR$ be a Ramsey class. Then every $(\RR, \UU)$-multi-amalgamation class is a Ramsey class.
\end{theorem}

If we wish to prove a class $\KK$ of $L$-structures is Ramsey, the following theorem from \cite{NR} provides, in many cases, a suitable $\RR$ for applying Theorem \ref{theorem:multiamalg}.

\begin{theorem} \label{theorem:NR}
Let $L$ be a finite relational language, such that $<$ is a binary relation in $L$. The class of all finite $L$-structures in which $<$ is interpreted as a linear order is a Ramsey class.
\end{theorem}

\section{The Classes $\KK_0$ and $\KK$}	
Let $\Lambda$ be a finite distributive lattice, ${\AA_\Lambda}$ the class of all finite $\Lambda$-ultrametric spaces, and $\vec \AA_\Lambda$ a well-equipped lift. In order to to prove the locally finite completion property required in Theorem \ref{theorem:multiamalg}, we will need to lift $\vec \AA_\Lambda$ to a linguistically more complex class. The first part of the lift, adding elements representing equivalence classes, is isolated below. It is similar to the lift employed in Lemma 4.28 of \cite{HN} for metric spaces with jumps, and is common in model theory.

A $\KK_0$-structure is meant to be viewed as follows: the elements of $P_{E, 1}$ represent the $E$-classes of a $\Lambda$-ultrametric space, and $U_{E, E'}(x, y)$ holds if $x$ represents an $E$-class and $y$ represents the $E'$-class containing $x$. The metric is not explicitly present in the language of the lift, since it is encoded by the family $\set{U_{E, E'}}$. 

\begin{definition}
Let $\LL_0 = \set{\set{P_{E, 1}}_{E \in \Lambda}, \set{U_{E,E'}}_{E<E' \in \Lambda}}$, be a relational language where the $P_{E, 1}$ are unary and the $U_{E, E'}$ are binary. Let $\UU_U$ be the following closure description: the $U_{E, E'}$ are closure relations, and the root of $U_{E, E'}$ is a single point $x$ such that $P_{E, 1}(x)$.

Let $\KK_0$ consist of all finite $\UU_U$-closed $\LL_0$-structures for which the following hold.
		
		\begin{enumerate}
			\item The family $\set{P_{E, 1}}_{E \in \Lambda}$ forms a partition
			\item If $U_{E,E'}(x,y)$, then $P_{E, 1}(x)$ and $P_{E', 1}(y)$
			\item (Coherence) If $E < E' < E'' \in \Lambda$ and $U_{E, E'}(x, y)$, then $U_{E', E''}(y, z)$ implies $U_{E, E''}(x, z)$
			\item (Downward semi-closure) If $P_{E, 1}(x)$ and $P_{E', 1}(x')$, then there is at most one $y$ such that $P_{E \meet E', 1}(y)$ and $U_{E \meet E', E}(y, x)$, $U_{E \meet E', E'}(y, x')$
		\end{enumerate}
			
\end{definition}

\begin{definition}
Let $\leq_U$ be the relation defined on a $\KK_0$-structure by $x \leq_U y$ if there are $E, E' \in \Lambda$ such that $U_{E, E'}(x, y)$. If $x \leq_U y$ and we wish to specify that $y$ is an $E'$-class, we will write $x/E' = y$.
\end{definition}
		
		\begin{proposition}
			Let $K \in \KK_0$ and let $x, x' \in K$.
			Suppose $x/E_1 = z_1 = x'/E_1$, $x/E_2 = z_2 = x'/E_2$. Then there exists $y \in K$ such that $x/(E_1\meet E_2) = y = x'/(E_1\meet E_2)$.
		\end{proposition} 
		\begin{proof}
			We have $x \in P_{F_1, 1}, y \in P_{F_2,1}$, for some $F_1, F_2 \leq E_1, E_2$. Since $K$ is $\UU_{U}$-closed, there are unique $y = x/(E_1 \meet E_2)$ and $y' = x'/(E_1 \meet E_2)$. By coherence, $y/E_1 = z_1=y'/E_1, y/E_2 = z_2=y'/E_2$. By downward semi-closure, $y=y'$.
		\end{proof}
	
		\begin{definition} For $x, y \in K_0$, define $\delta(x, y)$ to be the least $E$ such that $x/E = y/E$. By the proposition above, this is well-defined.
		\end{definition}
		
		\begin{proposition} \label{proposition:triangle}
			Let $K \in \KK_0$. Then $\delta$  satisfies the triangle inequality in $K$.
		\end{proposition}
		\begin{proof}
			Suppose $\delta(x_1, x_2) = F$, $\delta(x_2, x_3) = F'$. Let $a = x_1/(F \join F')$ and $b = x_3/(F \join F')$. Then $a = x_2/(F \join F') = b$, so $\delta(x_1, x_3) \leq F \join F'$.
		\end{proof}
		
	However, the function $\delta$ is technically not a $\Lambda$-ultrametric, or even a $\Lambda$-pseudoultrametric, since in general $\delta(x, x) \neq \bbzero$. Note that $\delta$ encodes all the information present in the family $\set{U_{E, E'}}$.

\begin{definition}
To any $\Lambda$-ultrametric space $A$, we associate a structure $A^{eq}$, such that if $A \in \AA_\Lambda$, then $A^{eq} \in \KK_0$, as follows.	
\begin{enumerate}
	\item The universe of $A^{eq}$ is $\sqcup_{E \in \Lambda} A/E$.
	\item For each $E \in \Lambda$, label the elements of $A/E$ with the predicate $P_{E, 1}$.
	\item For each $E, E' \in \Lambda$ with $E < E'$, let $U_{E,E'}(x,y)$ hold if $P_{E, 1}(x)$, $P_{E', 1}(y)$, and $y$ represents the $E'$-class containing the $E$-class that $x$ represents.
\end{enumerate}
\end{definition}
Note that this is only a fragment of the full model-theoretic $A^{eq}$, since we are not adding equivalence classes for equivalence relations definable on $A^n$ for $n>1$.

Conversely, to any $K \in \KK_0$, we can associate a structure $A_K \in \AA_\Lambda$. The following construction can be viewed as considering each point in $K$ as representing an equivalence class and picking a generic point in each class, i.e. points that are no closer to each other than necessary.

\begin{definition} \label{def:AK}
 Let $K \in \KK_0$. For each $x \in K$, create a corresponding point $x_A \in A_K$. Then, let $d(x_A, x_A) = \bbzero$, and let distances between distinct points be defined by $d(x_A, y_A) = \delta(x, y)$. By Proposition \ref{proposition:triangle}, the result is a $\Lambda$-ultrametric space.
\end{definition}
		
\begin{proposition} \label{proposition:eqsub}
Let $K \in \KK_0$. Then $K$ embeds into $(A_K)^{eq}$.	
\end{proposition}
\begin{proof}
	For each $x \in K$, if $P_{E, 1}(x)$, we map $x$ to $x_A/E \in (A_K)^{eq}$. 
	
	Suppose, for $x \in K$, that $P_{E, 1}(x)$, and let $y=x$. Then $\delta(x, y) = E = \delta(x_A/E, y_A/E)$. 
	
	Now suppose $x, y \in K$ with $P_{E, 1}(x)$ and $P_{E', 1}(y)$ and $x \neq y$. Let $\delta(x, y) = F$. Then $d(x_A, y_A) = F$. Then in $(A_K)^{eq}$, the least $G \in \Lambda$ such that $x_A/G = y_A/G$ is $G=F$. Since $E, E' < F$, this means $\delta(x_A/E, y_A/E') = F$ as well. Thus, our map preserves the family $\set{P_{E,1}}$ and $\delta$, and so gives an embedding of $K$ into $(A_K)^{eq}$.
\end{proof}

Thus $\KK_0$ is exactly the closure under $\UU_U$-closed substructure of the class obtained by applying the $eq$ operation to $\AA_\Lambda$. We call such structures \textit{upward-closed}. We now consider an additional closure condition.
		
\begin{definition} \label{def: dclosed}
We say $K \in \KK_0$ is \textit{downward closed} if for any $x, y \in K$ such that $P_{E, 1}(x)$, $P_{F, 1}(y)$, and $\delta(x, y) = E \join F$, there is some $z \in K$ such that $P_{E \meet F, 1}(z)$ and $z \leq_U x,y$. 
\end{definition}
		
	\begin{lemma} \label{lemma:closure0}
	Let $K \in \KK_0$. Then there is a finite $\KK_0$-structure $cl_0(K)$ such that
	\begin{enumerate}
	\item $K$ embeds into $cl_0(K)$
	\item $cl_0(K)$ is downward closed
	\end{enumerate}
	\end{lemma}
	\begin{proof}
	By Theorem \ref{theorem:amalg}, $\AA_\Lambda$ is an amalgamation class. Let $\Gamma$ be the Fra\"\i ss\'e limit of $\AA_\Lambda$. Embed $A_K$ into $\Gamma$. Then $(A_K)^{eq}$ is contained in $\Gamma^{eq}$. Let $cl_0(K)$ be the algebraic closure of $(A_K)^{eq}$ in $\Gamma^{eq}$.
	
	Given $x \in cl_0(K)$, with $P_{E, 1}(x)$, for any $E' > E$, $x/E'$ is definable from $x$ by the formula $\phi(y) = U_{E, E'}(x, y)$, and so is in its algebraic closure. Thus $cl_0(K)$ is $\UU_{U}$-closed.
	
	By Proposition \ref{proposition:eqsub}, $(1)$ already holds of $(A_K)^{eq}$.
	
	 We now prove $(2)$. Let $x', y' \in \Gamma$, with $d(x', y') = E \join F$. Since the structure $A = \set{x', y', z'}$, with $d(x', z') = E$, $d(y', z') = F$, $d(x', y') = E \join F$, satisfies the triangle inequality, we have $A \in \AA_\Lambda$. Thus, as $\Gamma$ is the Fra\"\i ss\'e limit of $\AA_\Lambda$, there is some $z' \in \Gamma$ such that $d(x', z') = E$ and $d(y', z') = F$. Given $x, y \in cl_0(K)$ as in Definition \ref{def: dclosed}, there exist $x', y' \in \Gamma^{eq}$ such that $P_{\bbzero, 1}(x')$, $P_{\bbzero, 1}(y')$,  $x'/E = x$, $y'/F = y$, and $\delta(x', y') = E \join F$. Then there is a $z' \in \Gamma^{eq}$ such that $\delta(x', z') = E$, $\delta(y', z') = F$. Thus $z' \leq_U x, y$, and so we may take $z=z'/(E \meet F) \leq_U x, y$. Finally, there can be at most one such $z$, since there is at most one $E \meet F$-class contained in any given $E$-class and $F$-class, and so $z$ is in the algebraic closure of $\set{x, y}$.
	\end{proof}

We now define the full class to which we will lift structures from $\vec \AA_\Lambda$. This will combine adding elements representing equivalence classes with the technique of duplicating the structure and connecting the parts by bijections used in Proposition 4.31 of \cite{HN} for structures with multiple linear orders. Since we are using subquotient orders instead of linear orders, we only need to duplicate part of the structure for each subquotient order.

The reason multiple structures are used to handle multiple linear orders is that Theorem \ref{theorem:NR}, which we plan to use to provide an $\RR$ for Theorem \ref{theorem:multiamalg}, provides a class with only a single linear order. Thus, in \cite{HN}, each linear order is placed on a single copy of the structure, and the copies are ordered one after another to form a single linear order.

The relation $D_{E_1, E_2}(x_1, x_2, y)$ in the definition below is meant to be viewed as stating that $x_1$ and $x_2$ represent an $E_1$ and $E_2$ class, respectively, and $y$ represents their intersection. This intersection of equivalence classes is the reason the algebraic closure operation in the class below will be binary rather than unary.
	
	\begin{definition} \label{def:K}
	For each $E \in \Lambda$ let $N_E \geq 1$, and let $<_{1-types}$ be a linear order on $\set{(E, i) | E \in \Lambda, i \in [N_E]}$. Relative to these parameters, we define $\KK$, a class of structures in the relational language
		$$\LL = \LL_0 \cup \set{\set{P_{E, i}}_{E \in \Lambda, 2 \leq i \leq N_{E}}, \set{B_{E, i, j}}_{E \in \Lambda, i,j \in [N_E]}, \set{D_{E, E'}}_{E \neq E' \in \Lambda}, D^\exists, <}$$ 
	where the relations $P_{E, i}$ are unary, the $B_{E, i, j}$ are binary, the $D_{E, E'}$ are ternary, $D^\exists$ is binary, and $<$ is binary. Let $\KK$ consist of all finite $\LL$-structures for which the following hold.
			
	\begin{enumerate}
		\item The family $\set{P_{E, i}}_{E \in \Lambda, i \in [N_E]}$ forms a partition such that classes that agree in the first index have the same cardinality.
		\item The substructure on the points $x$ such that $P_{E, 1}(x)$ holds for some $E \in \Lambda$ is an $\LL$-expansion of a $\KK_0$-structure.
		\item $<$ is a linear order, which agrees with $<_{1-types}$ between 1-types, i.e. if $P_{E, i}(x), P_{F, j}(y)$, $(E, i) \neq (F, j)$, then $x < y \Rightarrow (E, i) <_{1-types} (F, j)$.
		\item $D^\exists(x, y)$ iff there exists a $z$ such that $U_{E, E \join E'}(x, z)$, $U_{E', E \join E'}(y, z)$.
		\item If $D_{E_1, E_2}(x_1, x_2, y)$, then $P_{E_1, 1}(x_1)$, $P_{ E_2, 1}(x_2)$, $D^\exists(x_1, x_2)$, $P_{E_1 \meet E_2, 1}(y)$, and $y \leq_U x_1, x_2$.
		\item $B_{E, i, j}$ is the graph of a bijection from the points of $P_{E, i}$ to the points of $P_{E, j}$.
		\item If $B_{E, i, j}(x, y)$ and $B_{E, j, k}(y, z)$, then $B_{E, i, k}(x, z)$.
	\end{enumerate}
	\end{definition}
	\begin{definition}		
	We also define a closure description $\UU_{\KK}$ for $\LL$, in which the relations $U_{E, E'}$, $B_{E, i, j}$, and $D_{E, E'}$ are closure relations. The root of $U_{E, E'}$ is a single point $x$ such that $P_{E, 1}(x)$. The root of $B_{E, i, j}$ is a single point $x$ such that $P_{E, i}(x)$. The root of $D_{E, E'}$ is a pair of points $x_1, x_2$ such that $P_{E, 1}(x_1)$, $P_{E', 1}(x_2)$, $D^\exists(x_1, x_2)$, $U_{E, E'}(x_1, x_2)$ if $E < E'$ or $U_{E', E}(x_2, x_1)$ if $E' < E$, and $x_1 < x_2$ if $(E, 1) <_{1-types} (E', 1)$ or $x_2 < x_1$ if $(E', 1) <_{1-types} (E, 1)$. 
	
	\end{definition}
	Although $\KK$ is not closed under taking substructures, the class of $\UU_K$-closed $\KK$-structures is closed under taking $\UU_K$-closed substructures.
	
	\begin{definition}The \textit{metric part} of $K \in \KK$ is the $\KK_0$-structure appearing in condition $(2)$ of Definition \ref{def:K}, with language $\LL_0$. 
	\end{definition}
	
	\begin{remark}
	For $K \in \KK$, we can assign a distance $\delta(x,y)$ between two points in the non-metric part of $K$ as well, by taking the distance between the points $x$ and $y$ are in bijection with in the metric part of $K$. 
	\end{remark}

	\begin{lemma} \label{lemma:Kextension}
	Let $K \in \KK$, let $K_0$ be the metric part of $K$, and let $K_0'$ be a $\KK_0$-structure containing $K_0$. Then there is a $\KK$-structure $K'$ such that $K \subset K'$ and the metric part of $K'$ is $K_0'$.
	
	Furthermore, if $K_0'$ is downward closed, $K'$ can be taken to be $\UU_{\KK}$-closed.
	\end{lemma}
\begin{proof}
For any $x, y \in K_0'$ such that for some $E, F \in \Lambda$, $P_{E, 1}(x)$, $P_{F, 1}(y)$, and $\delta(x, y) = E \join F$, add the relation $D^\exists(x, y)$. Furthermore, if $K_0'$ is downward-closed, there is a $z \leq_U x, y$ such that $P_{E \meet F}(z)$, so add the relation $D_{E, F}(x, y, z)$. 

Then, for every $x_1 \in K_0' \bs K_0$, perform the following. Let $E \in \Lambda$ be such that $P_{E, 1}(x_1)$.
	\begin{enumerate}
		\item for every $2 \leq i \leq N_E$, add a point $x_i$ to $P_{E, i}$
		\item for every $ i, j \in [N_E]$, add the relation $B_{E, i, j}(x_i, x_j)$  
	\end{enumerate}
Finally, complete $<$ arbitrarily to a linear order so that it still agrees with $<_{1-types}$ between 1-types.
\end{proof}
	
	\begin{lemma} \label{lemma:closure1}
	Let $K \in \KK$. Then there is a $\UU_{\KK}$-closed $\KK$-structure $cl(K)$ such that $K$ is a substructure of $cl(K)$.
	\end{lemma}
\begin{proof}
	Let $K_0$ be the metric part of $K$. Let $cl_0(K_0)$ be as in Lemma \ref{lemma:closure0}. Then apply Lemma \ref{lemma:Kextension} to $K$ with $K_0' = cl_0(K_0)$, and let $cl(K)$ be the resulting $K'$.
\end{proof}

\section{Transfer}

In this section, we show that to prove $\vec \AA_\Lambda$ is a Ramsey class, it is sufficient to prove that the class of $\UU_\KK$-closed $\KK$-structures is a Ramsey class.

The first definition describes how we lift an $\vec \AA_\Lambda$-structure to a $\KK$-structure.

\begin{definition} \label{def:prelift}		
Let $\vec A \in \vec{\AA_\Lambda}$, and let $\vec \Gamma$ be the Fra\"\i ss\'e limit of $\vec \AA_\Lambda$. Before we describe how to lift an $\vec \AA_\Lambda$-structure into $\KK$, we first fix the following parameters and notation.
\begin{enumerate}[(a)]
\item For each meet-irreducible $E \in \Lambda$, let $N_E$ be the number of subquotient orders with bottom-relation $E$, and for each meet-reducible $E \in \Lambda$, let $N_E = 1$.
\item Enumerate the subquotient orders with bottom-relation $E$ as $<_{E, i}$ for ${i \in [N_E]}$.
\item For each $E \in \Lambda$, choose $F'_E$ a cover of $E$, and for $E$ meet-reducible choose $F''_E > E$ such that $E = F'_E \meet F''_E$.
\item For each element of $\set{(E, i)|E \in \Lambda, i \in [N_E]}$ with $E$ meet-reducible, use the construction in Lemma \ref{lemma:sqoproduct}, with the above choices of $F'_E$ and $F''_E$, to produce  a quantifier-free formula $\phi_{E, i}$ defining a subquotient order with bottom relation $E$ and top relation $\bbone$ on $\vec \Gamma$.

For $E$ meet-irreducible, use Corollary \ref{corollary:extendsqo} instead. 

Finally, as noted in Remark \ref{rem:choice}, we may assume that whenever the construction has to choose between multiple subquotient orders with a given bottom relation, it chooses the first in our enumeration.
\item Fix a linear order $<_{1-types}$ on $\set{(E, i)| E \in \Lambda, i \in [N_E]}$. 
\item Let $A$ be the metric part of $\vec A$.
\end{enumerate}

Note that, although $\phi_{E, i}$ is defined on elements of $\vec A$, it naturally induces a linear order on the elements of $P_{E, i}$ in $A^{eq}$.
We now associate a $\KK$-structure to $\vec A \in \vec \AA_\Lambda$. Let $L_\KK(\vec A)$ be as follows:
		
	\begin{enumerate}
		\item Construct $A^{eq}$.
		\item For each $E \in \Lambda$, $2 \leq i \leq N_E$, create a copy of the elements of $P_{E, 1}$, and label the elements of that copy with the predicate $P_{E, i}$.
		\item For each $E \in \Lambda$, for each $i, j \in [N_E]$, let $B_{E, i, j}(x, y)$ if $P_{E, i}(x)$, $P_{E, j}(y)$ and $x$ and $y$ represent the same $E$-class.
		\item For each $E \in \Lambda$, $i \in [N_E]$, define $<$ on the elements of $P_{E, i}$ to agree with the order induced by $\phi_{E, i}$ on those elements.
		\item Extend $<$ to a total order by setting $x < y$ if $P_{E, i}(x)$, $P_{F, j}(y)$, and $(E, i) <_{1-types} (F, j)$.
	\end{enumerate}
\end{definition} 

This gives a canonical lifting from $\vec \AA_\Lambda$ to $\KK$. However, we would like to lift elements of $\vec \AA_\Lambda$ to \textit{$\UU_\KK$-closed} $\KK$-structures, which will be done as follows. First, note that we can also apply the $L_\KK$-construction to $\vec \Gamma$.

\begin{definition} \label{def:lift}
Let $\vec A \in \vec{\AA_\Lambda}$, and let $\vec \Gamma$ be the Fra\"\i ss\'e limit of $\vec \AA_\Lambda$. Embed $\vec A$ into $\vec \Gamma$. This induces an embedding of $L_\KK(\vec A)$ into $L_\KK(\vec \Gamma)$, and let $\Lift(\vec A)$ be the algebraic closure of $L_\KK(\vec A)$ in $L_\KK(\vec \Gamma)$.
\end{definition}

By the proof of Lemma \ref{lemma:closure0}, $\Lift(\vec A)$ will be $\UU_\KK$-closed.

\begin{proposition} \label{prop:liftfunc}
Suppose $\vec A, \vec B \in \vec \AA_\Lambda$ and $\vec A$ embeds into $\vec B$. This induces an embedding from $L_\KK(\vec A)$ into $L_\KK(\vec B)$, and there is an embedding of $\Lift(\vec A)$ into $\Lift(\vec B)$ that extends this embedding. In particular, $\Lift(\vec A)$ is well-defined up to isomorphism over $\vec A$.
\end{proposition}
\begin{proof}
The definition of the $\Lift$ operation is based on an embedding of $L_\KK(\vec B)$ into $L_\KK(\vec \Gamma)$. Given such an embedding, it induces a corresponding embedding of $L_\KK(\vec A)$, and relative to these embeddings, $\Lift(\vec A)$ will then be a substructure of $\Lift(\vec B)$.

For the final point, take the embedding of $L_\KK(\vec A)$ to be an isomorphism of $L_\KK(\vec A)$ with $L_\KK(\vec B)$, or its inverse.
\end{proof}

Note that the $\Lift$ operation produces only a subset of the structures in $\KK$, since the order cannot be generic within $1$-types, but is controlled by the $(\phi_{E, i})$ and the $(<_{E, i})$, which remain definable in the lifted structure as appropriate restrictions of $<$. The Ramsey property for $\vec \AA_\Lambda$ corresponds more directly to the Ramsey property for the class of lifted structures, but working with $\UU_\KK$-closed $\KK$-structures reduces much of the bookkeeping. The next definition will provide a way to transfer the Ramsey property from $\UU_\KK$-closed $\KK$-structures to the class of lifted structures (or rather its closure under $\UU_\KK$-closed substructure).

\begin{definition} \label{def:qfr}
Fix a relational language $L$. A \textit{quantifier-free reinterpretation scheme} is a family of quantifier-free $L$-formulas $\Phi = \set{\phi_R : R \in L}$ such that $\phi_R$ has $n_R$ free variables, where $n_R$ is the arity of R.

A quantifier-free reinterpretation scheme $\Phi$ naturally induces a function $f_\Phi$ from $L$-structures to $L$-structures, where $f_\Phi(A)$ is given by reinterpreting each $R \in L$ as $\phi_R$. We call such a function a \textit{quantifier-free reinterpretation}.

Given a class $\RR$ of $L$-structures and a quantifier-free reinterpretation $f_\Phi$, if $f_\Phi$ is a retraction when restricted to $\RR$, we call it a \textit{quantifier-free retraction on $\RR$}. The image of a quantifier-free retraction is a \textit{quantifier-free retract of $\RR$}
\end{definition}

\begin{example}
Let $L$ consist of two binary relations, $<_1, <_2$. Let $\RR$ be the class of all finite $L$-structures where $<_1, <_2$ are linear orders. Then the class $\KK_1 \subset \RR$ consisting of structures where $<_1 = <_2$ is a quantifier-free retract of $\RR$, induced by the quantifier-free reinterpretation scheme $\Phi = \{\phi_{<_1}(x_1, x_2) = x_1 <_1 x_2, \phi_{<_2}(x_1, x_2) = x_1 <_1 x_2 \}$. 

Similarly, the class $\KK_2 \subset \RR$ for which $<_2 = <_1^{opp}$ is a quantifier-free retract, induced by the quantifier-free reinterpretation scheme $\Phi = \{\phi_{<_1}(x_1, x_2) = x_1 <_1 x_2, \phi_{<_2}(x_1, x_2) = x_2 <_1 x_1 \}$.
\end{example}

The first of the above examples essentially appeared in \cite{Bod}, where it was used to argue that if a class of structures with two generic linear orders had the Ramsey property, one could forget one of those linear orders and keep the Ramsey property. We now generalize that argument.

\begin{lemma} \label{lemma:forget}
Let $\KK_1$ be a Ramsey class, and $\KK_2$ a quantifier-free retract of $\KK_1$ such that $\KK_2 \subset \KK_1$. Then $\KK_2$ is a Ramsey class.
\end{lemma}	
\begin{proof}
	Let $A, B \in \KK_2$. Then we also have $A, B \in \KK_1$, and so there is some $C_1 \in \KK_1$ witnessing the Ramsey property for $A, B \in \KK_1$. Let $C_2 \in \KK_2$ be the retract of $C_1$. We claim $C_2$ witnesses the Ramsey property for $A, B \in \KK_2$.
	
	Let $\chi_2$ be a coloring of $\binom{C_2}{A}$. Let $f_\Phi$ be the quantifier-free reinterpretation, as in Definition \ref{def:qfr}, and recall $f_\Phi$ restricts to the identity on copies of $A, B$. Define a coloring $\chi_1$ of $\binom{C_1}{A}$ by $\chi_1(A) = \chi_2(f_\Phi(A))$. Let $\widehat B \subset C_1$ be a monochromatic copy of $B$. Then $f_\Phi(\widehat B) \subset C_2$ is a monochromatic copy of $B$, since $f_\Phi$ is the identity on $\widehat B$. 
\end{proof}

The idea of the retraction we will use is that in any lifted structure, $<$ is determined by the $(\phi_{E, i})$, and certain restrictions of $<$. In a $\KK$-structure, we can take these restrictions of $<$, forget the rest of $<$, and then use the $(\phi_{E, i})$ to define a new order from the restrictions. This is carried out in detail in the next definition.
  
 In the following definition, given $x \in P_{E, i}$ and $F >E$, $x/F$ is the $y \in P_{F, 1}$ such that $y \geq_U x'$, where $B_{E, i, 1}(x, x')$. In order for formulas involving $x/F$ to be quantifier-free, we must consider the relations $B_{E, i, j}$ and $U_{E, E'}$ to be functions. These relations define functions in $\UU_\KK$-closed $\KK$-structures, which is the reason for restricting ourselves to $\UU_\KK$-closed structures below.

\begin{definition} \label{def:K'}
We now define another class $\KK'$ using a quantifier-free retract on $\KK$.

Let $F'_E$, $F''_E$ be as in Definition \ref{def:prelift}.

For each $E \in \Lambda$, we inductively define a formula $\psi_{E, i}$ which gives a linear order on $P_{E, i}$. The case $E = \bbone$ is trivial. Now assume we have defined such $\psi_{F, i}$ for all $F > E$. 

If $E$ is meet-irreducible, let $$\psi_{E, i}(x, y) \Leftrightarrow (\delta(x, y) = F' \meet x < y) \join (\delta(x, y) > F'_E \meet \psi_{F'_E, 1}(x/F'_E, y/F'_E))$$

If $E$ is meet-reducible, let $$\psi_{E, i}(x, y) \Leftrightarrow$$ $$(x \neq y) \meet (((\delta(x, y) = F'_E \meet \psi_{F''_E, 1}(x/F''_E, y/F''_E)) \join (\delta(x, y) \neq F'_E \meet \psi_{F'_E, 1}(x/F'_E, y/F'_E)))$$

Let $\Phi_< = \set{\phi_<}$, where $\phi_<$ is
$$
(\bigvee_{\substack{E \in \Lambda \\ i \in [N_E]}} P_{E, i}(x) \wedge P_{E, i}(y) \wedge \psi_{E, i}(x, y)) \vee (\bigvee_{\substack{E, F \in \Lambda \\ i \in [N_E], j \in [N_F] \\ (E, i) <_{1-types} (F, j)}} P_{E, i}(x) \wedge P_{F, j}(y))
$$

Note that $f_\Phi$ restricts to the identity on structures for which $<$ is suitably encoded by the $(\psi_{F, i})$ and certain restrictions of $<$ on $P_{E, i}$ for meet-irreducible $E$. 

We thus let $\KK'$ be the quantifier-free retract of $\UU_\KK$-closed $\KK$-structures under the above quantifier-free reinterpretation scheme.
\end{definition}

\begin{remark}
Any $\KK'$-structure is $\UU_\KK$-closed.
\end{remark}

\begin{proposition}
Let $\vec A \in \vec \AA_\Lambda$. Then $\Lift(\vec A) \in \KK'$.
\end{proposition}
\begin{proof}
It is clear that $\Lift(\vec A)$ is a $\UU_\KK$-closed $\KK$-structure. We must check that $f_{\Phi_<}$ is the identity on $\Lift(\vec A)$. But $\Phi_<$ was defined so as to make this the case. 
\end{proof}

\begin{proposition} \label{proposition:KtoK'}
Suppose the class of $\UU_{\KK}$-closed $\KK$-structures is a Ramsey class. Then $\KK'$ is a Ramsey class.
\end{proposition}
\begin{proof}
This follows by Lemma \ref{lemma:forget}.
\end{proof}

\begin{proposition} \label{proposition:representation}
Let $K \in \KK'$. Then there is an $\vec A_K \in \vec \AA_\Lambda$ such that $K$ embeds into $\Lift(\vec A_K)$.
\end{proposition}
\begin{proof}
	Let $K_0$ be the metric part of $K$. Taking $A_{K_0}$ as in Definition \ref{def:AK} gives a structure in $\AA_\Lambda$, which needs to be expanded by certain subquotient orders $<_{E, i}$ in order to obtain a structure in $\vec \AA_\Lambda$. 
	
	The $<_{E, i}$ are determined in the following manner: we know that $<_{E, i}$ should have bottom relation $E$, and let $E'$ be its prescribed top-relation. Recall that each $x \in K_0$ gives a point $x_A \in A_{K_0}$. Let $<^*_{E, i}$ be the partial order on points of $A_{K_0}/E$ of the form $x_A/E$ defined as follows. For $x, y \in K_0$,  let $x_A/E <^*_{E, i} y_A/E$ if there are points $x', y' \in P_{E, i}$ such that $B_{E, 1, i}(x, x')$ and $B_{E, 1, i}(y, y')$, we have $\delta(x, y) \leq E'$ and $x' < y'$. 
	
	Then let $<_{E, i}$ be an arbitrary extension of $<^*_{E, i}$ to a subquotient order of $A_{K_0}$ with bottom relation $E$ and top relation $E'$. Then the resulting structure is the desired $\vec A_K$. 
\end{proof}

Thus the class $\KK'$ is exactly the closure under $\UU_\KK$-closed substructure of the class obtained by applying the $\Lift$ operation to $\vec \AA_\Lambda$.
	
\begin{lemma} \label{lemma:K'toA}
Suppose $\KK'$ is a Ramsey class. Then $\vec{\AA_\Lambda}$ is a Ramsey class.
\end{lemma}
\begin{proof}
Let $\vec A, \vec B \in \vec\AA_\Lambda$. Then $\Lift(\vec A), \Lift(\vec B)$ are $\KK'$-structures, and so there is some $C \in \KK'$ witnessing the Ramsey property for $\Lift(\vec A), \Lift(\vec B)$. By possible enlarging $C$, we may assume it has the form $\Lift(\vec C)$ for some $\vec C \in \vec \AA_\Lambda$. We will show that $\vec C$ witnesses the Ramsey property for $\vec A, \vec B$.

Let $\chi$ be a coloring of $\binom{\vec C}{\vec A}$. We wish to lift $\chi$ to a coloring $\widehat \chi$ of $\binom{C}{\Lift(\vec A)}$. 

\begin{claim}
Let $\vec X \in \vec A_\Lambda$, $\widehat X \in \KK'$, and $\widehat X \cong \Lift(\vec X)$. Then there is a unique substructure $\vec X_1$ of $\widehat X$ such that $(\vec X_1, \widehat X) \cong (L_\KK(\vec X), \Lift(\vec X))$.
\end{claim}
\begin{claimproof}
Since $\widehat X$ is equipped with a linear order, it is rigid, and hence there is a unique isomorphism of $\Lift(\vec X)$ with $\widehat X$. The claim follows.
\end{claimproof}

For any $\widehat X \cong \Lift(\vec X)$, we define $ker(\widehat X)$ to be the unique substructure such that $(ker(\widehat X), \widehat X) \cong (L_\KK(\vec X), \Lift(\vec X))$. 

Also, given a structure $X$ of the form $L_\KK(\vec X)$ we define a map $L_X^{opp}$ from $P_{\bbzero, 1} \subset X$  to $\vec A_X$ as defined in Proposition \ref{proposition:representation}, which sends $x$ to the corresponding point $x_A$. Note that if $\vec X \subset \vec Y$, and thus $L_\KK(\vec X) \subset L_\KK(\vec Y)$, then $L_{L_\KK(\vec Y)}^{opp}[L_\KK(\vec X)] \cong \vec X$. Furthermore, if we identify $L_{L_\KK(\vec Y)}^{opp}[L_\KK(\vec Y)]$ with $\vec Y$, then $L_{L_\KK(\vec Y)}^{opp}[L_\KK(\vec X)] = \vec X$.

\begin{claim}
There is a coloring $\widehat \chi$ of $\binom{C}{\Lift(\vec A)}$ such that, for $\widehat A \in \binom{C}{\Lift(\vec A)}$
$$\widehat \chi(\widehat A) = \chi(L_{ker(C)}^{opp}[ker(\widehat A)])$$
\end{claim}
\begin{claimproof}
Because $\widehat A \subset C$, $ker(\widehat A) \subset ker(C)$. Then, since $\widehat A \cong \Lift(\vec A)$, we have $L_{ker(C)}^{opp}[ker(\widehat A)] \in \binom{\vec C}{\vec A}$.
\end{claimproof}

By the Ramsey property for $C$, there is $\widehat B \cong \Lift(\vec B)$ in $C$ which is $\widehat \chi$-monochromatic. We now check $L_{ker(C)}^{opp}[ker(\widehat B)]$ is $\chi$-monochromatic.

If $\vec A_1 \subset \vec B$ with $\vec A_1 \cong \vec A$, then $L_\KK(\vec A_1) \subset L_\KK(\vec B)$, and by Proposition \ref{prop:liftfunc} this extends canonically to an embedding of $\Lift(\vec A_1)$ into $\widehat B$. Thus $\chi(\vec A_1) = \widehat \chi(\widehat A_1)$, with $\widehat A_1$ the image of $\Lift(\vec A_1)$ in $\widehat B$. Thus $\vec B$ is $\chi$-monochromatic.
\end{proof}

Thus, we have the following.
\begin{corollary} \label{lemma:KtoA}
Suppose the class of $\UU_{\KK}$-closed $\KK$-structures is a Ramsey class. Then $\vec{\AA}_\Lambda$ is a Ramsey class.
\end{corollary}
	
\section{Ramsey Theorems}
We now use Theorem \ref{theorem:multiamalg} to prove the class of $\UU_\KK$-closed $\KK$-structures is a Ramsey class.

We first consider the downward-closed $\KK_0$-structures (Definition \ref{def: dclosed}).

\begin{lemma} \label{lemma:K0amalg}
The class of downward closed $\KK_0$-structures is a strong amalgamation class.
\end{lemma}
\begin{proof}
Let the downward closed $\KK_0$-structures $K_1, K_2$ be the factors of an amalgamation problem, and let $K^*$ be their free amalgam. 
\begin{claim*}
$K^* \in \KK_0$. 
\end{claim*}
\begin{claimproof}
Since the $U_{E, E'}$ are unary, and since both factors and the base are $\UU_U$-closed, $K^*$ is $\UU_U$-closed.

 We only check downward semi-closure, since the other constraints follow immediately from the fact that they are satisfied in each factor. Let $x, y \in K^*$ with $P_{E, 1}(x)$, $P_{F, 1}(y)$, and  $\delta(x, y) = E \join F$. 
 
 If $x, y$ are not in the same factor, then there is no $z \leq_U x, y$, since the base is $\UU_U$-closed. If they are in the same factor, then in each factor there is at most one $z \leq_U x, y$ such that $P_{E \meet F, 1}(z)$, since each factor is downward closed. Thus, the only possible contradiction would be if $x, y$ were in the base, and there were $z_1, z_2$ in the first and second factor, respectively, such that $z_i \leq_U x, y$ and $P_{E \meet F, 1}(z_i)$. But this is impossible, since the base is also downward closed.
 \end{claimproof}

Then $cl_0(K^*)$ as provided by Lemma \ref{lemma:closure0} is a downward closed strong amalgam.
\end{proof}

\begin{lemma} \label{lemma:Kamalg} 
The class of $\UU_{\KK}$-closed $\KK$-structures is a strong amalgamation class.
\end{lemma}
\begin{proof}
Let the $\UU_\KK$-closed $\KK$-structures $K_1, K_2$ be the factors of an amalgamation problem. Let $K^*$ be the free amalgam of $K_1, K_2$, and let $K^*_0$ be the metric part of $K^*$. 

\begin{claim*} We can complete $K^*$ to a $\KK$-structure, $K^+$.
\end{claim*}
\begin{claimproof}
 By the claim in Lemma \ref{lemma:K0amalg}, $K^*_0 \in \KK_0$. Because the base is $\UU_\KK$-closed, $\set{B_{E, i, j}}$ are the graphs of the desired bijections in $K^*$. For any $x, y \in K^*$ such that for some $E, F \in \Lambda$, $P_{E, 1}(x)$, $P_{F, 1}(y)$, and $\delta(x, y) = E \join F$, add the relation $D^\exists(x, y)$. Finally, complete $<$ to a linear order that agrees with $<_{1-types}$ between 1-types (this doesn't conflict with any transitivity constraints). The remaining constraints are satisfied since they are satisfied in each factor. Thus, the resulting structure, $K^+$, is in $\KK$.
 \end{claimproof}

Then $cl(K^+)$ as provided by Lemma \ref{lemma:closure1} is a strong amalgam.
\end{proof}

\begin{lemma} \label{lemma:lfc}
	Let $\RR$ be the class of all finite $\LL$-structures where $<$ is a linear order. Then the class of $\UU_\KK$-closed $\KK$-structures has the locally finite completion property with respect to $(\RR, \UU_\KK)$.
\end{lemma}
\begin{proof}
	 Let $B$ be a $\UU_\KK$-closed $\KK$-structure and $C_0 \in \RR$. Set $n(B, C_0) = 0$. Let $C$ be a $\UU_{\KK}$-closed $\LL$-structure with $\RR$-completion $C_0$. We first note that it is sufficient to produce $C'$ a $\KK$-completion of $C$ with respect to copies of $B$, since then $cl(C')$ as provided by Lemma \ref{lemma:closure1} will be a $\UU_{\KK}$-closed $\KK$-completion of $C$ with respect to copies of $B$.
	 
	 Since we only need to produce a $\KK$-completion of $C$ with respect to copies of $B$, rather than a $\KK$-completion, we may assume that $C$ is a union of copies of $B$ and that all relations are between points which lie in a common copy of $B$. However, assuming this means we may only assume that $C$ is \textit{$\UU_\KK$-semi-closed}, meaning that the $R$-out-degree, for any closure relation $R$, of a tuple $\vec t$ is at most 1 if $\vec t$ represents an embedding of the root of $R$ into $C$, and 0 otherwise. (We could actually assume full closure for the relations $U_{E, E'}$ and $B_{E, i, j}$, since they represent unary functions, but it will not be necessary.) We claim this places the following constraints on $C$.
	
	\begin{enumerate}
		\item The family $\set{P_{E, i}}$ forms a partition.
		\item $C$ is $\UU_U$-closed.
		\item The  $U_{E, E'}$ are coherent.
		\item $<$ is an irreflexive, asymmetric, acyclic relation. 
		\item $<$ agrees with $<_{1-types}$ between 1-types.
		\item If $U_{E, E'}(x, y)$, then $P_{E, 1}(x)$, $P_{E', 1}(y)$.
		\item If $D^{\exists}(x,y)$, then there are $E, E' \in \Lambda$ such that $P_{E, 1}(x)$, $P_{E', 1}(y)$, and $\delta(x, y) = E \join E'$.
		\item If $D_{E_1, E_2}(x_1, x_2, y)$, then we have $P_{E_1, 1}(x_1)$, $P_{ E_2, 1}(x_2)$, $D^\exists(x_1, x_2)$, $P_{E_1 \meet E_2, 1}(y)$, and $y \leq_U x_1, x_2$.
		\item If $x, y$ are such that $P_{E, 1}(x)$, $P_{E', 1}(y)$, and $D^\exists(x, y)$, then there is exactly one $z$ such that $D_{E, E'}(x, y, z)$.
		\item $B_{E, i, j}$ is the graph of a bijection from $P_{E, i}$ to $P_{E, j}$.
		\item If $B_{E, i, j}(x, y)$ and $B_{E, j, k}(y, z)$, then $B_{E, i, k}(x, z)$.
		\item $C$ is downward semi-closed.
	\end{enumerate}
	With the exceptions of $(3)$, $(4)$, $(10)$, $(11)$, and $(12)$, the constraints are immediate from the assumption that $C$ is $\UU_\KK$-semi-closed, is a union of copies of $B$, and all relations are between points that lie in the common copy of $B$.
	
	 Constraint $(4)$ follows from the assumption that there is a homomorphism-embedding from $C$ to $C_0$, and $<$ is a linear order on $C_0$.
	
	Before continuing, we observe that if $a \in C$ lies in a given copy of $B$ (perhaps one of several), and $U_{E, E'}(a, b)$ or $B_{E, i, j}(a, b)$, then $b$ lies in that same copy of $B$. For suppose $a \in \widehat B$, but $b \not\in \widehat B$. Then, since $\widehat B$ is $\UU_\KK$-closed, there is some $b' \in \widehat B$ such that we also have $U_{E, E'}(a, b')$ (resp. $B_{E, i, j}(a, b')$). But this is forbidden, since $C$ is $\UU_{\KK}$-semi-closed.
	
	We check constraint $(3)$. Suppose $U_{E, E'}(x, y)$, $U_{E', E''}(y, z)$. By our observation, $x, y, z$ all lie in a single copy of $B$, and so $U_{E, E''}(x, z)$. Constraint $(11)$ follows similarly.
	
	Constraint $(10)$ holds since each $B_{E, i, j}$ is a union of bijections, and $C$ is $\UU_\KK$-semi-closed.
	
	Finally, suppose constraint $(12)$ is violated, so there are $x, y \in C$ such that $P_{E, 1}(x)$, $P_{E', 1}(y)$, and $\delta(x, y) = E \join E'$, and there are distinct $z_1, z_2$ such that $P_{E \meet E', 1}(z_i)$ and $z_i \leq_U x, y$. By our earlier observation, we must have $x, y, z_1$ lying in a single copy of $B$, as well as $x, y, z_2$. Since $B$ is $\UU_\KK$-closed, we then have $D_{E, E'}(x, y, z_1)$ and $D_{E, E'}(x, y, z_2)$. But this is forbidden, since $C$ is $\UU_{\KK}$-semi-closed. 

	We define an equivalence relation $P$ on $C$ whose classes are the family $\set{P_{E, i}}$. Taking the transitive closure of $<$ gives a partial order, for which $P$ is a congruence and which agrees with $<_{1-types}$ between $P$-classes, and so $<$ can be completed to a linear order which is $P$-convex and agrees with $<_{1-types}$ between $P$-classes. Thus, after adding the relation $D^\exists$ where appropriate, we can complete $C$ to a $\KK$-structure $C'$.
\end{proof}

\begin{theorem} \label{theorem:KRamsey}
The class of $\UU_\KK$-closed $\KK$-structures is a Ramsey class.
\end{theorem}
\begin{proof}
By Lemmas \ref{lemma:Kamalg} and \ref{lemma:lfc}, the class of $\UU_\KK$-closed $\KK$-structures is an $(R, \UU_K)$-multi-amalgamation class, where $\RR$ is the class of finite $\LL$-structures where $<$ is interpreted as a linear order. Thus, by Theorems \ref{theorem:multiamalg} and \ref{theorem:NR}, it is a Ramsey class.
\end{proof}

\begin{theorem} \label{theorem:mainthm}
Let $\Lambda$ be a finite distributive lattice, $\AA_\Lambda$ the class of all finite $\Lambda$-ultrametric spaces, and $\vec \AA_\Lambda$ a well-equipped lift of $\AA_\Lambda$. Then $\vec \AA_\Lambda$ is a Ramsey class.
\end{theorem}
\begin{proof}
By Theorem \ref{theorem:KRamsey} and Corollary \ref{lemma:KtoA}, we are done.
\end{proof}

\begin{corollary} \label{corollary:RamseyPermutations}
The amalgamation classes corresponding to all the homogeneous finite-dimensional permutation structures constructed in Theorem \ref{theorem:amalg} are Ramsey.

In particular, for every finite distributive lattice $\Lambda$, there is a Ramsey class such that $\Lambda$ is isomorphic to the lattice of $\emptyset$-definable equivalence relations in the class's \fraisse limit.
\end{corollary}
\begin{proof}
For the first part, all such classes are representable as well-equipped lifts of the class of all $\Lambda$-ultrametric spaces for some finite distributive $\Lambda$.

For the second part, Corollary \ref{theorem:Representation} shows that for every such $\Lambda$, there is a homogeneous finite-dimensional permutation structure with lattice of $\emptyset$-definable equivalence relations isomorphic to $\Lambda$.
\end{proof}

\section{The Expansion Property} \label{sect:expansionproperty}
We now identify a Ramsey lift of $\AA_\Lambda$ with the expansion property as defined below, and use this to compute the universal minimal flow of $Aut(\Gamma)$, where $\Gamma$ is the \fraisse limit of $\AA_\Lambda$. See Chapter \ref{sect:dynamicsBackground} for a review of the necessary background.

We first note a general lemma we will make use of.

\begin{lemma}
Suppose $F^*$ is a precompact expansion of a homogeneous structure $F$, and for every $A^* \in \text{Age}(F^*)$, there is a $B \in \text{Age}(F)$ such that for any expansion $B^* \in \text{Age}(F^*)$ of $B$, $A^*$ embeds into $B^*$. Then $\text{Age}(F^*)$ has the expansion property relative to $\text{Age}(F)$.
\end{lemma}
\begin{proof}
Let $A \in \text{Age}(F)$, and enumerate the expansions of $A$ in $\text{Age}(F^*)$, which there are finitely many of by precompactness, as $(A^*_i)_{i=1}^n$. For each $i \in [n]$, let $B_i \in \text{Age}(F)$ be the structure provided by hypothesis for $A^*_i$. By the joint embedding property, let $B \in \text{Age}(F)$ embed every $B_i$. Then $B$ witnesses the expansion property for $A$.
\end{proof}

\begin{definition}
Let $\Lambda$ be a finite distributive lattice. Enumerate the meet-irreducibles as $E_i$, and for each $i$, let $E_i^+$ cover $E_i$. Then $\vec \AA_\Lambda^{\min}$ is the class of all expansions of elements of $\AA_\Lambda$ by a single subquotient order for each $i$, with bottom relation $E_i$ and top relation $E_i^+$.
\end{definition}

The first step is to reduce the desired expansion property to one which will be easier to prove.

\begin{definition}
 Let $\KK'^{\min}$ be the closure under $\UU_\KK$-closed substructure of structures of the form $\Lift(\vec A)$ for some $\vec A \in \vec \AA_\Lambda^{\min}$. 
 
 Let $\KK'^{\min}_r$ be the reduct of $\KK'^{\min}$ forgetting the order.
 \end{definition}

\begin{lemma} \label{lemma:expansiontransfer}
Suppose $\KK'^{min}$ has the expansion property relative to $\KK'^{min}_r$. Then $\vec \AA_\Lambda^{\min}$ has the expansion property relative to $\AA_\Lambda$.
\end{lemma}
\begin{proof}




Let $\vec A \in \vec \AA_\Lambda^{\min}$. By assumption, there is a $B \in \KK'^{min}_r$ witnessing the expansion property for $\Lift(\vec A)$. Let $B_0$ be the metric part of $B$, and consider $A_{B_0} \in \AA_\Lambda$. We claim $A_{B_0}$ witnesses the expansion property for $\vec A$.

Let $\vec A_{B_0} \in \vec \AA_\Lambda^{\min}$ be an arbitrary expansion of $A_{B_0}$. This induces an expansion $\vec B \in \KK'^{min}$ of $B$, such that $\vec A_{\vec B}$, as defined in Proposition \ref{proposition:representation}, is isomorphic to $\vec A_{B_0}$. 

Then $\vec B$ embeds $\Lift(\vec A)$, and therefore embeds $L_\KK(\vec A)$. Composing the natural injection of $\vec A$ into $L_\KK(\vec A)$ with the embedding of $L_\KK(\vec A)$ into $\vec B$, and composing the result with the map from $\vec B$ to $\vec A_{B_0}$ sending $x$ to $x_A$, gives an embedding of $\vec A$ into $\vec A_{B_0}$.
\end{proof}

We now use a standard argument (see \cite{digraphs}*{Theorem 8.6} for example, although it appeared earlier) to prove the expansion property for the lifted classes

\begin{lemma} \label{lemma:expansion}
$\KK'^{min}$ has the expansion property relative to $\KK'_r$.
\end{lemma}

\begin{proof}
Let $\vec A \in \KK'^{\min}$. Enumerate the meet-irreducibles in $\Lambda$ as $(E_i)_{i=1}^n$, and let $E_i^+$ cover $E_i$. For each $i$, consider the structure on $\set{x_1, x_2, x_3}$, such that
\begin{enumerate}
 \item $P_{E_i, 1}(x_1)$, $P_{E_i, 1}(x_2)$, $P_{E_i^+, 1}(x_3)$ \item $x_1, x_2 \leq_U x_3$
 \item $x_1 < x_2$
 \item The remaining order information is determined by $<_{1-types}$.
 \end{enumerate}
  Let $\vec p_i \in \KK'^{\min}$ be the minimal $\UU_\KK$-closed substructure of the \fraisse limit of $\KK'^{min}$ containing the above structure.
 
 By possibly enlarging $\vec A$, we may assume it has the form $\Lift(\vec B)$ for some $\vec B \in \vec \AA_\Lambda^{min}$. For each $i \in [n]$, let $<_{E_i}$ be the unique subquotient order on $\vec B$ from $E_i$ to $E_i^+$. For each $I \subset [n]$, let $\vec B_I$ be $\vec B$, but with $<_{E_i}$ reversed for every $i \in I$, and let $\vec A_I \cong \Lift(\vec B_I)$.

 Let $\vec C_0 \in \KK'^{min}$ embed $\vec A_I$ for every $I \subset [n]$, and let $\vec C_{i+1} \rightarrow (\vec C_i)^{\vec p_{i+1}}_2$ for every $i \in [n]$. Let $C_n \in \KK'_r$ be the reduct of $\vec C_n$. We will show $C_n$ witnesses the expansion property for $\vec A$.

Let $(C_n, \prec) \in \KK'^{min}$ be an expansion of $C_n$. For each $i \in [n]$, let $\prec_i$ be the restriction of $\prec$ to pairs $x, y \in P_{E_i, 1}$ such that $\delta(x, y) = E_i^+$, let $<_{E_i, 1}$ be the corresponding restriction of $<$, and let $\chi_i$ be a coloring of $\binom {\vec C_n}{\vec p_i}$ defined by
$$\chi_i(\vec p_i) =  \begin{cases} 
      0 & x_2 \prec_i x_1 \\
      1 & x_1 \prec_i x_2 \\
   \end{cases}$$
Then iterated applications of the Ramsey property yield $\vec C_0^* \subset \vec C_n$, a copy of $\vec C_0$ in which, for every $i \in [n]$, either $\prec_i = <_{E_i, 1}$, or $\prec_i = <_{E_i, 1}^{opp}$. Let $I \subset [n]$ be such that $\prec_i = <_{E_i, 1}^{opp}$ iff $i \in I$. We have that $\vec C_0^*$ contains a copy $\vec A_I^*$ of $\vec A_I$. Letting ${A_I}^*$ be the reduct of $\vec A_I^*$ to $\KK'_r$, we have $\vec A \cong \left(A_I^*, \set{\prec_i}_{i = 1}^n \right)$.
\end{proof}

\begin{theorem} \label{theorem:UMF}
 Let $\Lambda$ be a finite distributive lattice, $\AA_\Lambda$ the class of all finite $\Lambda$-ultrametric spaces, $\Gamma$  the \fraisse limit of $\AA_\Lambda$, and $\vec \Gamma^{min} = {(\Gamma, (<_{E_i})_{i=1}^n)}$ the \fraisse limit of $\vec \AA_\Lambda^{min}$. Then 
 \begin{enumerate}
 \item $\vec \AA_\Lambda^{\min}$ is a Ramsey class and has the expansion property relative to $\AA_\Lambda$.
 \item The logic action of $\text{Aut}(\Gamma)$ on $\overline{\text{Aut}(\Gamma) \cdot (<_{E_i})_{i=1}^n}$ is the universal minimal flow of $\text{Aut}(\Gamma)$.
 \end{enumerate}
\end{theorem}
\begin{proof}
By Theorem \ref{theorem:mainthm}, $\vec \AA_\Lambda^{\min}$ is a Ramsey class. By Lemmas \ref{lemma:expansiontransfer} and \ref{lemma:expansion}, $\vec \AA_\Lambda^{\min}$ has the expansion property relative to $\AA_\Lambda$. The second part then follows by Theorem \ref{theorem:NVT}.
\end{proof}
\chapter{The Decision Problem for Joint-Embedding}\label{chap:JEP}
\section{Introduction} \label{sec:jepIntro}
\setcounter{figure}{0}
 We now take up a different type of issue relating to the intersection of model theory and combinatorics in the vicinity of permutation avoidance classes.
 When studying finitely-presented structures, it is natural to ask whether various properties are decidable, and such questions have been considered for finitely-constrained permutation avoidance classes for some time, as evidenced by Ru\v{s}kuc's talk \cite{Rusk}, where, among several other questions, the decidability of \textit{atomicity} was raised. A permutation avoidance class is called atomic if it cannot be expressed as a union of two proper subclasses. One of the primary uses of atomicity is the following lemma, which reduces the calculation of growth rates for a wide class of permutation avoidance classes to the calculation for atomic classes (see \cite{Vatter} for a reference). 
 
 \begin{lemma}
 Suppose $\KK$ is a permutation avoidance class, with no infinite antichain in the containment order. Then $\KK$ can be expressed as a finite union of atomic subclasses. Furthermore, the upper growth rate of $\KK$ is equal to the maximum upper growth rate among its atomic subclasses. 
 \end{lemma}
 
 Atomicity is easily proven to be equivalent to the joint-embedding property (see \cite{Vatter}), and so we may rephrase Ru\v{s}kuc's question as follows.
 
  \begin{question}
  Is there an algorithm that, given finite set of forbidden permutations, decides whether the corresponding permutation avoidance class has the joint-embedding property?
  \end{question}
  
  This problem is known to be decidable in certain restricted classes of permutations, such as \textit{grid classes} \cite{Waton}. 
  Also, whether a permutation avoidance class is a \textit{natural class}, which is a stronger condition, is decidable \cite{Murphy}.
  
  We believe there is a strong possibility this decision problem is undecidable in general. We are not aware of many undecidability results in the permutation avoidance classes literature, although \cite{Pak}, using methods that seem quite different from ours,  proves an undecidability results about comparing the parity of the number of permutations of size $n$ in two permutation classes. As a first approximation to Ru\v{s}kuc's problem, we examine the corresponding problem in the category of graphs. This problem seems to be easier due to the fact that the edge relation in the infinite random graph has the independence property, and so the base theory does not constrain how we can add edges between the factors during the joint-embedding procedure. As the generic permutation is NIP, we cannot hope for such behavior there. For graphs, we prove the following theorem, via a reduction to the tiling problem.
  
  \begin{theorem} \label{theorem:inducedJEP}
 There is no algorithm that, given a finite set of forbidden induced subgraphs, decides whether the corresponding hereditary graph class has the JEP.
  \end{theorem}
  
  This is first proven for graphs enriched by a sufficient supply of unary predicates, and then a formal reduction to the pure graph language is given. A very rough sketch of the proof is as follows. The first two steps ensure that the tiling problem is equivalent to whether we can joint-embed two particular graphs, and the third step ensures that joint-embedding for the class is equivalent to joint-embedding for those two graphs.

\begin{enumerate}
\item Construct two graphs $A^*$, representing a grid, and $B^*$ representing a suitable collection of tiles.
\item Choose a finite set of constraints to ensure that successfully joint-embedding $A^*$ and $B^*$ requires producing a valid tiling of the grid points in $A^*$ with the tiles from $B^*$
\item Show that if the tiling problem admits a solution, then the chosen class admits a joint-embedding procedure.
\end{enumerate}

Finally, we consider a variation on the JEP called the joint homomorphism property, which is of interest for infinite-domain constraint satisfaction problems. Modifying our proof of Theorem \ref{theorem:inducedJEP} yields the following theorem.

  \begin{theorem}
 There is no algorithm that, given a finite set of forbidden induced subgraphs, decides whether the corresponding hereditary graph class has the joint homomorphism property.
  \end{theorem}

\section{Graphs with Unary Predicates and the Canonical Models}
We will work in the following language.
\begin{enumerate}
\item $E$: a binary relation, to represent edges
\item $O^i, P'^i, G^i, T^1$ for $i \in \set{0, 1}$: unary predicates, which will denote origin vertices, non-origin path vertices, grid vertices, and tile vertices 
\item $C_i$ for $1 \leq i \leq 4$: additional coding vertices
\end{enumerate}

We also define a unary predicate $P^i = O^i \cup P'^i$, which will denote path vertices.

As is evident from the presence of the $C_i$, we are already doing some coding; the most natural language would instead have an additional two colored edge relations and a directed edge relation.

\subsection{The Canonical Models} \label{sec:canonical}
We here further flesh out steps $(1)$ and $(2)$ from the proof sketch.

$A^*$ will consist of a 1-way infinite directed path, with vertices in $P^0$, and a marked origin in $O^0$. Directed edges will be simulated using the edge relation $E$ and coding vertices of type $C_1$ and $C_2$. To every pair of points in this path, we attach a $G^0$-vertex, representing a grid point with coordinates taken from the related path points. Because we must distinguish between $x$ and $y$-coordinates, we use coding vertices of type $C_3$ and $C_4$ to simulate two additional types of colored edge, and use these edges to attach the grid point to its coordinates. 

$B^*$ will look like a copy of $A^*$, using 1-superscripted predicates instead, but with a path of length $T$ (where $T$ is the number of tile types in the given tiling problem) $T^1$-vertices attached to each $G^1$-point. These represent a full tile-set available at each coordinate, with the different tile-types being distinguished by their distance from the corresponding $G^1$-point.

When we try to joint-embed $A^*$ and $B^*$, we wish our constraints to force the following: for every $G^0$-point in $A$, with coordinates $(x, y)$, we must add an edge to exactly one tile-point attached to the $G^1$-point in $B$ with the same coordinates. This is interpreted as tiling the point $(x, y)$ by the corresponding tile-type, and our constraints should further enforce the local tiling rules.

For the particular classes of structures we are dealing with here, namely graphs with forbidden induced subgraphs, our choice of $B^*$ is rather baroque. We could have simply chosen $B^*$ to be a collection of $n$ tile points, with some further coding to distinguish the different tile-types. However, the construction presented here is more flexible and better adapted to handling more complex classes of structures.

In particular, the relevant fact seems to be that in the \fraisse limit of the class of all graphs, the edge relation, which is the relation we are using to connect a tile to a given grid point, has a model-theoretic property called the \textit{independence property}, which corresponds to the presence of arbitrary induced bipartite subgraphs. Thus, the base theory puts no constraints on how we may connect a set of tile points to a set of grid points; we may build an arbitrary bipartite graph between them. Such constraints only arise from the induced subgraphs we have chosen to forbid. In cases where the relation we are using to connect tiles to grid points does not have the independence property, the base theory adds additional constraints that must be handled. In permutations this manifests in the requirement that the orders be transitive, and such constraints require the flexibility of the two-grid construction. 

\section{Constraints}
In addition to the constraints forcing a valid tiling to be attempted when joint-embedding the canonical models, we have several constraints which ensure that the origin, path, and grid points encode something grid-like. We would ideally be able to choose further constraints which ensure that every structure in our class $\GG_\TT$ looks like $A^*$ or $B^*$. We would like every grid point to have coordinates from the path, or every $G^1$-point to have a complete tile-set. However, as we cannot enforce such ``totality'' conditions using forbidden structures, we must allow for partial structures.

In the previous section, we noted that we would wish our constraints to force a $G^0$-point to be tiled using a tile from a $G^1$-point with the same coordinates. However, as we are forbidding a \textit{finite} number of finite structures, our constraints must have a \textit{local} character; as figuring out the coordinates of a grid point requires walking back to the origin, and thus looking at an unbounded number of vertices, we cannot use our constraints as desired. Instead, we will start the tiling at the origin, and then propagate it by local constraints. 

We will now precisely state our constraints, but first will establish some notation.

\begin{definition} \label{def:graphRelations}

We first define the ``special'' edges our construction will use.
\begin{enumerate}
\item $x \rightarrow^i y$ if $x, y \in P^i$ and there exist $a \in C_1, b \in C_2$ such that $xEaEbEy$. In this case, we say $x$ is the predecessor of $y$, and $y$ the successor of $x$.
\item $\Pi^i_1(v, w)$ if $v \in G^i$, $w \in P^i$, and there exists $a \in C_3$ such that $vEaEw$. In this case, we say $w$ is an $x$-projection of $v$.
\item $\Pi^i_2(v, w)$ if $v \in G^i$, $w \in P^i$, and there exists $a \in C_4$ such that $vEaEw$. In this case, we say $w$ is an $y$-projection of $v$.
\end{enumerate}

We say $g$ is a \textit{$G^i$-origin}, or sometimes a \textit{grid origin}, if there is an $x \in O^i$ such that $\Pi^i_1(g, x)$ and $\Pi^i_2(g, x)$.

Our constraints will force $x$ and $y$-projections to be unique. Given $g, g'$ in $G^i$, we say $g'$ is a horizontal successor of $g$ if they have the same $y$-projection, and the $x$-projection of $g'$ is the $\rightarrow$-successor of the $x$-projection of $g$. Similarly for vertical successor, but with $x$ and $y$ switched.

We now define binary relations related to the tiles.
\begin{enumerate}
\item For $i \in [T]$, we say $\tau_i(x, y)$ if $x \in G^1$, $y \in T^1$, and there exist $v_1, ..., v_i \in T^1$ such that $v_i = y$ and $xEv_1E...Ev_i$. In this case, we say $y$ is a \textit{tile of type $i$ associated to $x$}.
\item $\tau(x, y)$ if $x \in G^0$, $y \in G^1$, and there exist a $t$ such that $\tau_i(y, t)$ for some $i \in [T]$. In this case, we say \textit{$x$ is tiled by $y$} or that \textit{$x$ is tiled by a tile of type $i$}.
\end{enumerate}

Finally, we say $x$ has a $\textit{full set of tiles}$ if there exist $t_i$ for $i \in [T]$ such that for all $i$, $\tau_i(x, t_i)$.
\end{definition}

Given a tiling problem $\TT$, we now define $\GG_\TT$ as the class of all finite graphs with the following constraints.
\begin{enumerate}
\item The unary predicates in the language are disjoint.
\item A path vertex has at most 1 $\rightarrow$-predecessor
\item An origin vertex has no $\rightarrow$-predecessor.
\item A grid vertex has at most 1 $x$-projection and 1 $y$-projection.
\item Tile vertices are associated to at most one grid point, i.e. given $t \in T^1$, there do not exist distinct $g, h \in G^1$ such that $\tau_i(g, t)$ and $\tau_j(h, t)$.
\item Tile vertices have a unique type, i.e. if $\tau_i(g, t)$ and $\tau_j(g, t)$ then $i=j$.
\item The tiling rules of $\TT$ are respected.

Suppose $\TT$ forbids a tile of type $j$ to the right of (respectively, above) a tile of type $i$. Then we forbid the following as a non-induced subgraph.

 Let $g, g' \in G^0$ with $g'$ a horizontal (resp. vertical) successor of $g$. Let $h, h' \in G^1$ with $h'$ a horizontal (resp. vertical) successor of $h$. Finally, let $\tau(g, t_{h,i}), \tau(g', t_{j', j})$ where $t_{h, i}$ is a tile of type $i$ associated to $h$ and $t_{h', j}$ is a tile of type $j$ associated to $h'$.
 \item If $g \in G^0$ and $h \in G^1$ are grid-origins, and $h$ has a full set of tiles, then $g$ must be tiled by a tile associated to $h$.
 \item Suppose $g, g' \in G^0$ with $g'$ a horizontal (resp. vertical) successor of $g$, and $h, h' \in G^1$ with $h'$ a horizontal (resp. vertical) successor of $h$. Suppose $\tau(g, t)$ where $t$ is a tile vertex associated with $h$. If $h'$ has a full tileset, then $g'$ must be tiled by a tile vertex associated with $h'$.
\end{enumerate}

We note that only the last two constraints require the presence of edges, and so are the only ones that require forbidding induced subgraphs.

\section{The Proof for Graphs with Unary Predicates}

We wish to prove the following.

\begin{proposition} \label{prop:unaryundecidable}
Let $\TT$ be a tiling problem, and $\GG_\TT$ be the hereditary graph class defined above. Then $\GG_\TT$ has the JEP iff $\TT$ has a solution.
\end{proposition}

\subsection{An Informal Proof}
\begin{proof}
\textit{The easy direction: from the JEP to a tiling} \\
Suppose $\GG_\TT$ has the JEP. Note that $A^*, B^*$ as described above are in $\GG_\TT$, so we may joint-embed them. By (8), the origin of the grid in $A^*$ must be tiled, and by (9) this must propagate to a tiling of the whole grid, while respecting the tiling rules by (7). We may then read a solution to the tiling problem off the resulting graph.

\textit{The delicate direction: from a tiling to the JEP} \\
Here, we are a bit sketchier. We first fix a solution $\theta: \N^2 \to [T]$ to the tiling problem $\TT$. Given $A, B \in \GG_\TT$, we initially take the disjoint union $C = A \sqcup B$.

As only constraints (8) and (9) require the presence of edges, these are the only constraints that may be violated at this point, and in fact only (8) may be. We thus use $\theta(0, 0)$ to tile all $G^0$-origins in one factor from all full tilesets attached to $G^1$-origins in the other factor. However, now there may be violations of constraint (9). We continue using $\theta$ to appropriately tile our grids. The key point here is constraints (2-4) ensure that every grid point we must work with has well-defined coordinates, so we have a definite input to give to $\theta$.
\end{proof}

\subsection{From the JEP to a Tiling} \label{sec:JEPtoTiling}
For this direction, we may largely repeat the informal version.

Let $\Pi^0 = \set{p^0_i | i \in \N}$, and let $\Gamma^0 = P^2$, whose elements we denote $g^0_{i, j}$ rather than $(p_i, p_j)$. Let $A^*$ have vertex set $\Pi^0 \cup \Gamma^0$, with $p_0 \in O^0$, $Pi^0\bs\set{p_0} \subset P'^0$, and $\Gamma^0 \subset G^0$. Also, add  coding vertices in $C_i$ and the associated edges needed to encode the relations $p_i \rightarrow^0 p_{i+1}$ for each $p_i \in \Pi^0$, and $\Pi^0_1(g_{i,j}, p_i)$ and $\Pi^0_2(g_{i,j}, p_j)$.

 Let $B^*$ be constructed as $A^*$, but using 1-superscripted points, sets, and predicates in place of 0-superscripted ones. Let $\Theta^1 = \Gamma^1 \times [T]$, and denote its elements as $t^1_{g, i}$ rather than $(g, i)$, and add these vertices to $B^*$. Finally, for each $g \in \Gamma^1$, add edges so that $gEt_{g,1}E...Et_{g,T}$.

By inspection, $A^*, B^* \in \GG_\TT$. Let $C \in \GG_\TT$ joint-embed $A^*$ and $B^*$. We claim $C$ encodes a solution to $\TT$.

By (1), no points in $A^*$ and $B^*$ got identified in $C$, except perhaps coding vertices. By (8) and (9), for every $(i,j) \in \N^2$ there is some $k \in [T]$ such that $\tau(g^0_{i,j}, t^1_{g^1_{i,j}, k})$. Define the function $\theta: \N^2 \to [T]$ by picking one such $k$ for each $(i,j)$. By (7), $\theta$ yields a solution to $\TT$.

\subsection{From a Tiling to the JEP}

For this section, we fix a solution $\theta: \N^2 \to [T]$ to $\TT$.

We begin by establishing some effects of constraints (2-4), which will allow us to assign coordinates to grid points. We note that, although it would add little additional overhead, it is not necessary to constrain the number of $\rightarrow$-successors, and so constraints (2) and (3) actually allow the path vertices to form a forest.

\begin{definition}
Let $\rightarrow^i_n$ be the $n$-fold composition of $\rightarrow^i$.

Given $p \in P^i$ and $o \in O^i$, we say \textit{$p$ is on a path with origin $o$} if there is some $n \in \N$ so that $o \rightarrow^i_n p$. In this case, we say \textit{$p$ is at distance $n$ from $o$}.

Let $G^i_*$ be the set of all $g \in G^i$ such that there exist $o \in O^i$ and $x, y \in P^i$ with $\Pi^i_1(g, x), \Pi^i_2(g, y)$ and $x$ and $y$ are on paths with origin $o$. In this case, if $x$ is at distance $n$ from $o$, and $y$ at distance $m$, we say \textit{$g$ has coordinates $(n,m)$}.
\end{definition} 

Constraints (2) and (3) ensure that if $p$ is on a path with origin $o$ and a path with origin $o'$, then $o= o'$. They also ensure that the distance of $p$ from $o$ is unique. This, together with constraint (4), ensures that the coordinates of a grid point are unique. 

\begin{definition}
Let $\theta_*: G^0_* \to [T]$ be defined by $\theta_*(g) = i$ iff $g$ has coordinates $(n, m)$ and $\theta(n,m) = i$.
\end{definition}

We are now ready to state our joint-embedding procedure. Let $A, B \in \GG_\TT$. Let $C_0$ be the disjoint union $A \sqcup B$. We construct an extension $C$ of $C_0$ by adding edges of the form $(g, t)$ when the following conditions are met.
\begin{enumerate}
\item $(g, t) \in A \times B \cup B \times A$
\item $g \in G^0$
\item $g$ has coordinates $(n,m)$ for some $n,m \in \N$
\item There is $h \in G^1$ with coordinates $(n,m)$ such that $\tau_{\theta_*(g)}(h, t)$
\end{enumerate}

\begin{remark}
This procedure may add many more tiling-relations than would be required to satisfy the constraints. For example, we tile any grid point with coordinates, even if preceding grid points are missing that block propagation from the origin, and we may tile using tiles from incomplete tilesets.
\end{remark}

We now wish to show that $C \in \GG_\TT$ by showing it satisfies each constraint.

As constraint (1) only involves unary predicates, and these remain unchanged by taking the disjoint union and adding edges, it remains satisfied in $C$. 

\begin{lemma}
$C$ satisfies constraints (2-6).
\end{lemma}
\begin{proof}
For all these constraints, the forbidden configuration is connected, and thus they are satisfied in $C_0$. However, our procedure then only adds edges from $G^0$ to $T^1$-vertices, which by constraint (1) are not of any other type. As none of the forbidden configurations involve both $G^0$ and $T^1$-vertices, such edges cannot cause them to be violated, and so they continue to be satisfied in $C$. 
\end{proof}

\begin{lemma} \label{lemma:Con7}
$C$ satisfies constraint (7).
\end{lemma}
\begin{proof}
Again, our constraint is connected, and so satisfied in $C_0$. Fix a violation of (7), say of the horizontal rule, with vertices as in the constraint description. As we only add edges from $G^0$-vertices to $T^1$-vertices, we must have added either the edge $(g, t_{h,i})$ or $(g', t_{h', i})$. However, if we have only added one such edge, the configuration without that edge would be connected and would have been present in $C_0$, and so be entirely contained in one factor. This is a contradiction, as we only add edges between points in distinct factors. Thus we must have added both these edges.

Thus $t_{h, i}$ is a tile of type $\theta_*(g)$ and $t_{h', j}$ is a tile of type $\theta_*(g')$, and by constraints (5) and (6) these types are unique. Suppose $g$ has coordinates $(n, m)$; as $g'$ is a horizontal successor of $g$, it must have coordinates $(n+1, m)$. But then $t_{h, i}$ is of type $\theta(n, m)$ and $t_{h, j}$ is of type $\theta(n+1, m)$, so they cannot violate (7).
\end{proof}

\begin{lemma} \label{lemma:Con8}
$C$ satisfies constraint (8).
\end{lemma}
\begin{proof}
Let $X=\set{g, c, d, o}, Y=\set{g', c', d', o', t'_1, ..., t'_T}$, and $X \cup Y$ witness a violation of constraint (8), with $o \in O^0$, $g \in G^0$ with $x$ and $y$-projections equal to $o$, and $c$ and $d$ the requisite coding vertices; let $g', c', d', o'$ be a corresponding configuration using 1-superscripted predicates, and let $t_i \in T^1$ for each $i \in T$ with $g'Et'_1E...Et'_T$.

As $X$ and $Y$ are each connected, they must each lie in a single factor, and these factors must be distinct. Thus in $C_0$, $g$ and $g'$ both have coordinates $(0, 0)$ and $g'$ has a full tileset, so our procedure adds an edge from $g$ to $t'_{\theta(0,0)}$, which satisfies the constraint. 
\end{proof}

\begin{lemma}
$C$ satisfies constraint (9).
\end{lemma}
\begin{proof}
Consider a violation of constraint (9), with labels as in the constraint description (although the violation also requires suitable path and coding vertices). As in Lemma \ref{lemma:Con7}, we must have added the edge from $g$ to $t$. As in Lemma \ref{lemma:Con8}, the violation then splits into two connected components, one in each factor; one component contains $g, g'$, and their associated path and coding vertices while the other contains $h, h'$, and their associated tilesets and path and coding vertices.  

As we added an edge from $g$ to $t$, $g$ and $h$ must have had coordinates in $C_0$. Thus $g'$ and $h'$ also have coordinates in $C_0$. As $h'$ has a full tileset in $C_0$, our procedure adds an edge from $g'$ to a tile in this tileset, which satisfies the constraint. 
\end{proof}

\section{Moving to the Language of Graphs} \label{sec:movingToGraphs}
Given a finitely-constrained hereditary class $\GG_\TT$ in the language with unary predicates, we wish to produce a finitely-constrained hereditary graph class that has the JEP iff $\GG_\TT$ does. For this, we need some means of interpreting the unary predicates in the pure graph language. Our plan is to associate the $i^{th}$ unary predicate to some graph $G_i$, and to represent ``$v$ is in the $i^{th}$ predicate'' by freely joining a copy of $G_i$ over $v$. In order for this coding to be unambiguous, the graphs we choose must form an antichain under embeddings.

 We remark that we do not actually require an infinite antichain in the following definition, merely one with as many graphs as we have unary predicates. For our argument, the minimum size will be 13.

\begin{definition} \label{def:graphAntichain}
We now fix an infinite collection of 2-connected graphs with basepoints $(G_i, a_i)_{i \in \N}$, such that $\set{G_i}_{i \in \N}$ is an antichain under embeddability.

As we do not require that there are no automorphisms of any $G_i$ moving the basepoint, we will refer to any image of $a_i$ under an automorphism as a \textit{possible basepoint} of $G_i$.

\end{definition}

\begin{definition}
Let $\CC_k$ be the class of finite graphs with $k$ unary predicates, which we will refer to as colors $\set{1, ..., k}$. Let $\CC^*_k \subset \CC_k$ be the subclass in which the colors partition the vertices, and in which any (colored) copy of the $G_i$ are forbidden.
\end{definition}

\begin{definition} \label{def:wedge}
Define $\wedge: \CC^*_k \to \set{graphs}$ as follows: for each vertex of the graph, if it has color $i$, freely attach a copy of $G_i$ over it as the basepoint. These copies of $G_i$ will be called \textit{attached copies}.

The image of $A \in C^*_k$ will be denoted by $\widehat A$. We will also denote the pointwise image of $\GG \subset C^*_k$ as $\widehat \GG$.
\end{definition}

\begin{lemma} \label{lemma:attachedcopies}
Let $G \in \CC^*_k$. Any copy of $G_i$ in $\widehat G$ is an attached copy.
\end{lemma}
\begin{proof}
As $G_i$ is 2-connected, any copy must be contained in a single block of $\widehat G$. As the copies of $G_i$ are freely attached, the blocks of $\widehat G$ are those of $G$ as well as the attached $G_j$ for various $j$. Thus, any copy of $G_i$ must be contained in one of the attached $G_j$. As $\set{G_i}$ is an antichain, it must be one of the attached copies of $G_i$.
\end{proof}

\begin{lemma} \label{lemma:hatinverse}
Let  $\vee: \set{graphs} \to \CC_k$ be given by taking a graph, and for each copy of $G_i$ free over its basepoint (picking one such basepoint if there are several), retaining the basepoint and giving it color $i$, and forgetting the remaining vertices. Then $\vee(\widehat G) \cong G$.

In particular, $\wedge$ is injective. 
\end{lemma}
\begin{proof}
This is nearly immediate from Lemma \ref{lemma:attachedcopies}. The only subtlety is that we have not required that each $G_i$ be rigid, and so if there are automorphisms moving the basepoint, there will be multiple possible basepoints to choose from for a given copy of $G_i$. However, we claim this does not matter.

There are two cases to consider. If the copy of $G_i$ has no external edges incident upon any of its vertices, then $\vee$ will send it to a single isolated vertex with color $i$, no matter which possible basepoint we pick. If there is an external edge incident upon a vertex, that vertex must be the basepoint, as $G_i$ is freely attached over its basepoint.
\end{proof}

\begin{lemma} \label{lemma:imageconditions}
A graph is in the image of $\wedge$ iff it satisfies the following properties.
\begin{enumerate}
\item For each $i$, every copy of $G_i$ is free over its basepoint.
\item If $v$ is the basepoint of a copy $H_1$ of $G_i$ and $H_2$ of $G_j$, then $H_1 = H_2$.
\item Every vertex is, for some $i$, part of a copy of $G_i$.
\end{enumerate}
\end{lemma}
\begin{remark}
As in Lemma \ref{lemma:hatinverse}, it will not matter which basepoint of $G_i$ we choose for (1). Also, (2) implicitly uses that $\set{G_i}$ is an antichain.
\end{remark}
\begin{proof}
Suppose we start with $G \in \CC^*_k$. Then $\widehat G$ is produced by making each vertex the basepoint of a copy of $G_i$, for the appropriate $i$. Thus (3) is satisfied. Conditions (1) and (2) are satisfied by Lemma \ref{lemma:attachedcopies}.

Now suppose we are given a graph $G$ of this form. By conditions (1) and (2), the vertex set of $\vee(G)$ consists of the basepoints of copies of $G_i$, each given color $i$, and with edges between them induced by $G$. Then, using condition (3), we have $G = \widehat{ \vee(G)}$.
\end{proof}

\begin{lemma}
$\wedge$ preserves embeddings, i.e. there exists an embedding $A \into B$ iff there exists an embedding $\widehat A \into \widehat B$
\end{lemma}
\begin{proof}
The forward direction should be clear. 

For the other direction, suppose $\widehat A \into \widehat B$. Then for each copy of $G_i \subset \widehat A$, the basepoint (picking one if there are several) must be mapped to such a basepoint in $\widehat B$. By Lemma \ref{lemma:imageconditions}, each of these basepoints in $\widehat B$ has a free copy of $G_i$ over it, and so can be identified with a vertex in $\vee(\widehat B)$. Furthermore, it will receive the same color as the corresponding point in $\vee(\widehat A)$. Finally, $\vee$ preserves the induced graph on the points it retains, so $\vee(\widehat A) \into \vee(\widehat B)$, and so by Lemma \ref{lemma:hatinverse} we are finished.
\end{proof}

As $\wedge$ preserves embeddings, the class $\GG_\TT$ in the language with unary predicates will have the JEP iff its image under $\wedge$ does. However, this image is not a hereditary graph class, and it is not clear that its downward closure will be finitely-constrained. So our goal now is to find some finitely-constrained hereditary graph class such that every member can be completed to an element in the image of $\GG_\TT$ under $\wedge$, which must satisfy the conditions of Lemma \ref{lemma:imageconditions}.

The following constraints are meant to enforce conditions (1) and (2) of Lemma \ref{lemma:imageconditions}.

\begin{definition} \label{def:HHConstraints}
Let $\HH_1$ be the set of graphs consisting of, for each $i$: 
\begin{enumerate}
\item a copy of $G_i$ and an additional vertex adjacent to a point that is not a possible basepoint of $G_i$
\item a copy of $G_i$ and additional vertices $v_i, v_j$, possibly with $v_i = v_j$, adjacent to two distinct possible basepoints of $G_i$
\end{enumerate}

Let $\HH_2$ be the set of graphs consisting of a copy of $G_i$ and $G_j$ freely joined over their basepoints, for each $i, j$, allowing $i=j$.
\end{definition}

\begin{definition}
Given a set $\GG$ of graphs, we define $\neg \GG$ to be the corresponding hereditary graph class forbidding the graphs in $\GG$.
\end{definition}

Keeping in mind condition (3) of Lemma \ref{lemma:imageconditions}, the plan for our completion algorithm is to freely attach a copy of $G_i$ for some $i$ over every vertex that is not already in some copy of one of the $\set{G_i}$. However, randomly assigning colors may produce a forbidden structure. Thus, we make sure we have a ``dummy'' color, which is not in any non-trivial constraint, available and only use its associated $G_i$ for our completion.

\begin{lemma}
Let $\GG \subset \CC_k$, such that $\neg \GG \subset \CC^*_k$. Further suppose that the only graphs in $\GG$ containing a $k$-colored vertex are multicolored single vertices and colored copies of the $\set{G_i}$.

 Then every graph in $\neg(\widehat \GG \cup \HH_1 \cup \HH_2)$ embeds into one in $\widehat {\neg \GG}$.
\end{lemma}
\begin{proof}
Let $G \in \neg (\widehat \GG \cup \HH_1 \cup \HH_2)$. Since $G \in \neg \HH_1$, it satisfies (1) from Lemma \ref{lemma:imageconditions}. Since $\GG$ contains all multicolored single vertices, then since $G \in \neg (\widehat \GG \cup \HH_2)$, it also satisfies (2) from Lemma \ref{lemma:imageconditions}.

For every vertex $v$ for which there is no $i$ such that $v$ is in copy of $G_i$ free over its basepoint, we freely attach to $v$ a copy of $G_k$, identifying $v$ with the basepoint. Call the resulting graph $G^+$, and note it satisfies (3) from Lemma \ref{lemma:imageconditions}. 

Using the 2-connectedness of the $\set{G_i}$ as in Lemma \ref{lemma:attachedcopies}, $G^+$ still satisfies (1) and (2) from Lemma \ref{lemma:imageconditions}.

We claim it is also still in $\neg \widehat \GG$, as we have only added copies of $G_k$. Suppose $\widehat H \in \widehat \GG$ embeds into $G^+$. Then $H \in \GG$ embeds into $\vee(\GG^+)$. 

As $G^+$ satisfies (2) from Lemma \ref{lemma:imageconditions}, $H$ cannot be a multicolored vertex. As $G^+ \in \neg \HH_1$, $H$ cannot be a colored copy of any of the $\set{G_i}$. Thus $H$ does not contain any $k$-colored vertices.

Consider the subgraph $A \subset G^+$ induced by all vertices which are not the basepoint of a freely-attached copy of $G_k$. Then $H$ must embed into $\vee(A)$. But then $\widehat H$ embeds into $A$ and thus into $G$.
\end{proof}


\begin{lemma} \label{lemma:transfer}
Let $\GG \subset \CC^k$, such that $\neg \GG \subset \CC^*_k$. Further suppose that the only graphs in $\GG$ containing a $k$-colored vertex are multicolored single vertices and colored copies of the $\set{G_i}$. Then $\neg (\widehat \GG \cup \HH_1 \cup \HH_2)$ has the JEP iff $\neg \GG$ has the JEP.
\end{lemma}
\begin{proof}
Suppose $\neg \GG$ has the JEP. Let $A, B \in \neg (\widehat \GG \cup \HH_1 \cup \HH_2)$. Extend them to $A^+, B^+ \in \widehat {\neg \GG}$. Then, there is some $C \in \neg \GG$ embedding $\vee(A^+), \vee(B^+)$. Thus $\widehat C$ embeds $A^+, B^+$, and so $A, B$ as well.

Now suppose $\neg (\widehat \GG \cup \HH_1 \cup \HH_2)$ has the JEP. Let $A, B \in \neg \GG$. Then there is some $C \in \neg (\widehat \GG \cup \HH_1 \cup \HH_2)$ embedding $\widehat A, \widehat B$. Extend $C$ to $C^+ \in \widehat {\neg \GG}$. Then $\vee(C^+)$ embeds $A, B$.
\end{proof}

In order to finally prove our main theorem, we must choose a suitable set $\set{(G_i, a_i)}$. The graphs must be 2-connected and form an antichain under embedding. Finally, in order to have $\GG_\TT \subset \CC^*_k$, no colored version of them may embed into our canonical models $A^*, B^*$, and they must not be produced by our joint-embedding process for the graphs with unary predicates.

\begin{definition}
Given $n \in \N$, a \textit{necklace of $n$ triangles} is the graph obtained from $n$ triangles $\set{T_1, ..., T_n}$ by identifying a single point of $T_i$ with a point from $T_{i+1 \text{ (mod $n$)}}$ for each $i$.  
\end{definition}

\begin{notation}
For the remainder of this section, will let $G_i$ consist of a necklace of $i+2$ triangles, and allow any points for the basepoints.
\end{notation}

\begin{theorem}
  There is no algorithm that, given a finite set of forbidden induced subgraphs, decides whether the corresponding hereditary graph class has the JEP.
\end{theorem}
\begin{proof}
By Proposition \ref{prop:unaryundecidable}, it is undecidable whether $\GG_\TT$ has the JEP, as $\TT$ varies. We may modify $\GG_\TT$ to $\GG_\TT^*$ by introducing an extra color and forbidding all uncolored vertices. We claim we may also add constraints forbidding $\set{G_i}$, as well as constraints forbidding any two grid vertices from being connected to each other, or any two tile vertices from being connected to each other.

 Note that our canonical models contain no triangles (as the coding vertices break up edges), and thus no copies of the $\set{G_i}$, and they also satisfy the other new constraints. As our joint-embedding procedure only adds edges from grid vertices to tile vertices, by the second additional constraint this will produce no triangles, and thus no copies of any of the $\set{G_i}$. Again, because our joint-embedding procedure only adds edges from grid vertices to tile vertices, it will also preserve the other new constraints.

 We may thus to apply Lemma \ref{lemma:transfer} to $\GG^*_\TT$ to produce a family of finitely-constrained hereditary graph classes for which the JEP is undecidable as $\TT$ varies.
\end{proof}

\section{The Joint Homomorphism Property}
A class of structures has the \textit{joint homomorphism property (JHP)} if, given any two structures in the class, there is a third that admits homomorphisms from both. This notion naturally arises in infinite-domain constraint satisfaction problems. For example, the constraint satisfaction problem for a theory can be realized as the constraint satisfaction problem for a particular model iff the models of the theory have the JHP \cite{Bodthes}. The following question was posed by Bodirsky in January 2018 (personal communication).

\begin{question}
  Is there an algorithm that, given a finite set of forbidden induced subgraphs, decides whether the corresponding hereditary graph class has the JHP?
\end{question}

In this section, our main result is a negative answer to this question, obtained by modifying our construction for the JEP.

\begin{theorem} \label{theorem:jhpUndecidable}
 There is no algorithm that, given a finite set of forbidden induced subgraphs, decides whether the corresponding hereditary graph class has the JHP.
\end{theorem}

Theorem \ref{theorem:jhpUndecidable} will be proven by modifying our proof of Theorem \ref{theorem:inducedJEP}. The reader should be familiar with the brief sketch of the proof of Theorem \ref{theorem:inducedJEP} appearing at the end of Section \ref{sec:jepIntro} and the discussion at the beginning of Section \ref{sec:movingToGraphs} about removing the unary predicates; relevant results and definitions will be recalled or referenced as needed.

Unlike the JEP, the JHP is sensitive to changing between interdefinable languages. For example, we get the following as a corollary to Theorem \ref{theorem:inducedJEP}, but will later have to work much more without the non-edge relation present.

\begin{proposition} \label{prop:nonEdgeJHP}
Work in a language with relations for edges and non-edges. Then there is no algorithm that, given a finite set of forbidden induced subgraphs, decides whether the corresponding hereditary graph class has the JHP.
\end{proposition}
\begin{proof}
Our goal is to alter our canonical models (the graphs $A^*, B^*$ from the proof sketch at the end of Section \ref{sec:jepIntro}, although here we really want their interpretations in the pure graph language) so that any homomorphism is actually an embedding. 

Suppose we are given finite a set $\CC^{red}$ of forbidden induced subgraphs in the language with just the edge relation. Let $\CC$ be the set of graphs, in the enriched language, with the non-edge relation added between any non-adjacent points. Let $\CC^+$ be the union of $\CC$ with the graphs on two points in which either both relations or neither relation is present, ensuring the relations act as edges and non-edges. Then $\neg\CC^{red}$ has the JEP iff $\neg \CC^+$ has the JHP.
\end{proof}

As in Proposition \ref{prop:nonEdgeJHP}, the plan for proving Theorem \ref{theorem:jhpUndecidable} will be to modify our canonical models so that any $C$ witnessing the JHP also witnesses the JEP. In Proposition \ref{prop:nonEdgeJHP}, we did this by adding the non-edge relation between any two non-adjacent vertices to make our structures clique-like. Here we do the following.
\begin{enumerate}
\item Forbid $K_4$.
\item Find some graph $G$ such that any non-identity homomorphic image of $G$ contains a copy of $K_4$.
\item Over any two non-adjacent basepoints of $(G_i, a_i)$ (the graphs we are using to code unary predicates, see Definition \ref{def:graphAntichain}) in our canonical models, freely join a copy of $G$, while keeping the vertices non-adjacent. 
\end{enumerate}

The procedure above ensures that homomorphisms cannot identify the basepoints of the $(G_i, a_i)$ in our new canonical models, nor add edges between them, just as adding the non-edge relation did in Proposition \ref{prop:nonEdgeJHP}. The constraint set $\HH_1$ from Definition \ref{def:HHConstraints} ensures that we cannot add edges between non-basepoints and any points outside the copy of $G_i$ they lie in, nor can we identify such points. Thus the only possible issue is if the homomorphisms of our new canonical models fail to be embeddings within a single copy of some antichain element $G_i$.

This last possibility will be removed by forbidding all non-identity homomorphic images of each $G_i$ that we use from our antichain. However, these forbidden homomorphic images of $G_i$ might embed into $G_i$; to handle this point, we replace the $\set{G_i}$ with their \textit{cores}, for which this problem disappears. 

\begin{definition}
A \textit{retract} of a graph $G$ is a subgraph $H \subset G$ such that there exists a homomorphism $\phi:G \to H$ with $\phi \upharpoonright_H = id$. In particular, $H$ is an induced subgraph of $G$. 

Given a finite graph $G$, the \textit{core of $G$} is its unique (up to isomorphism) minimal retract.
\end{definition}  

The next lemma gives the key fact about cores. The proof may be found in standard references, e.g. \cite{HT}.

\begin{lemma}
Let $C$ be the core of $G$. Then any endomorphism of $C$ is an automorphism.
\end{lemma}

\begin{notation}
Let $W_5$ be the \textit{5-wheel}, i.e. a 5-cycle with an additional point adjacent to all others. 
\end{notation}

\begin{remark}
Every non-identity homomorphic image of $W_5$ contains a copy of $K_4$. Also, $W_5$ is 2-connected and every edge is contained in a triangle.
\end{remark}

when considering the JEP, we  let $G_i$ consist of a necklace of $i+2$ triangles, and allowed any points for the basepoints. However, we will need a different choice of $G_i$ for this section, and so now $N_i$ will refer to the necklace of $i+2$ triangles, with any points allowed for the basepoints.

We will now work toward constructing this section's choice of 2-connected graphs $\set{(G_i, a_i)}$ forming an antichain under embeddings, and prove some preparatory lemmas about them.

\begin{definition}
Given a graph $G$, we construct an \textit{augmented copy of $G$}, denoted $G^+$, as follows. First, we start with a copy of $G$. Then over every non-adjacent pair of vertices, we freely join a copy of $W_5$, identifying that pair of vertices with a pair of non-adjacent vertices in $W_5$.
\end{definition}

\begin{lemma} \label{lemma:k4Appearance}
Consider a graph $G^+$, and let $H \subset G^+$ be a copy of $G$ in which any 2 non-adjacent vertices have a copy of $W_5$ freely joined over them in $G^+$. Then a copy of $K_4$ embeds into the image of any homomorphism of $G^+$ that is not an embedding when restricted to $H$.

Furthermore, if $u, v \in H$ are non-adjacent and $\phi$ is a homomorphism of $G^+$ such that $\phi(u) \neq \phi(v)$, then $\phi(u)$ and $\phi(v)$ have 2 common neighbors in $\phi(G^+)$.
\end{lemma}
\begin{proof}
The first part is immediate from the definition of $H$ and the fact that any proper homomorphic image of $W_5$ embeds a copy of $K_4$.

For the second part, note that any 2 non-adjacent points of $H$ have 2 common neighbors in $G^+$, lying in the copy of $W_5$ freely joined over them, and these common neighbors are adjacent to each other.
\end{proof}

\begin{lemma} \label{lemma:G+antichain}
The $\set{N^+_i}$ are 2-connected, have every edge contained in a triangle, contain no copies of $K_4$, and form an antichain under homomorphisms.
\end{lemma}
\begin{proof}
The first two points are clear by inspection. For the third point, first notice that all the new edges in $N^+_i$ contain at least one new vertex; as $N_i$ contains no copies of $K_4$, any copy of $K_4$ in $N^+_i$ must contain at least one new vertex. However, any new vertex is only connected to other vertices in the same copy of $W_5$, and thus cannot be contained in any copy of $K_4$.

For the last point, let $\phi: N^+_i \to N^+_j$ be a homomorphism, with $i \neq j$. Let $H_i \subset N^+_i$ be the copy of $N_i$ such that any 2 non-adjacent vertices have a copy of $W_5$ freely joined over them in $N^+_i$, and define $H_j \subset N^+_j$ similarly. As $N^+_j$ contains no $K_4$, Lemma \ref{lemma:k4Appearance} implies $\phi$ must embed $H_i$ into $N^+_j$. 

This embedding of $H_i$ cannot be done solely using triangles from $H_j$. As every triangle in $N^+_j$ is either contained in $H_j$ or in some copy of $W_5$, the image of $H_i$ must contain some triangle $T$ from some copy $W$ of $W_5$. As any 2 triangles of $W$ that intersect in a single point are connected by an edge, no other triangle from $W$ can be in the image of $H_i$. However, any triangle from $W$ shares a point with at most one triangle outside $W$. Thus $T$ cannot be contained in a necklace of triangles. 
\end{proof}

\begin{lemma} \label{lemma:2ConOrK4}
Any homomorphic image of $N^+_i$ is either 2-connected or contains a copy of $K_4$.
\end{lemma}
\begin{proof}
Let $G$ be the homomorphic image of $N^+_i$ via the homomorphism $\phi:N^+_i \to G$. We will assume that no copy of $K_4$ appears in $G$, and show it is 2-connected. 

Let $H \subset N^+_i$ be the copy of $N_i$ such that any 2 non-adjacent vertices have a copy of $W_5$ freely joined over them in $N^+_i$. By Lemma \ref{lemma:k4Appearance}, we may assume the homomorphism is an embedding when restricted to $H$.

Let $B$ be the block of $G$ containing $\phi[H]$. We wish to show $B = G$. Otherwise, there is some copy of $W$ of $W_5$ in $N^+_i$ such that $\phi[W]$ is not contained in $B$. Note that $\phi$ must be an embedding when restricted to $W$, since otherwise $G$ would contain a copy of $K_4$. But then $\phi[W]$ is 2-connected and meets $B$ in two vertices, and so is contained in $B$.
\end{proof}

\begin{notation}
For each $i$, we will let $(G_i, a_i)$ be the core of $(N^+_i, a_i)$.
\end{notation}

\begin{lemma} \label{lemma:G'antichain}
The $\set{G_i}$ are 2-connected, have every edge contained in a triangle, and form an antichain under homomorphisms.
\end{lemma}
\begin{proof}
For each $i$, $G_i$ is an induced subgraph and homomorphic image of $N^+_i$. Thus $G_i$ contains no copy of $K_4$, and so is 2-connected by Lemma \ref{lemma:2ConOrK4}.

Let $u,v$ be adjacent in $G_i$. We wish to show $u,v$ have a common neighbor. Let $\phi:N^+_i \to G_i$ be a retract. We know $u, v$ have a common neighbor $w$ in $N^+_i$. But then $\phi(w)$ is a common neighbor of $u,v$ in $G_i$.

 For the last point, let $\phi: G_i \to G_j$ be a homomorphism, for $i \neq j$. Then as there exist homomorphisms from $N^+_i \to G_i$ and from $G_j \to N^+_j$, composition gives a homomorphism from $N^+_i \to N^+_j$, contradicting Lemma \ref{lemma:G+antichain}. 
\end{proof}

We will use $A^+, B^+$ to denoted the augmented copies of our canonical models $A^*, B^*$, where the new vertices we have added are marked with a new unary predicate $C_5$. Note that $A^+, B^+$ contain $A^*, B^*$ as induced subgraphs.

\begin{lemma} \label{lemma:unaryCanonical}
$A^+, B^+$ contain no homomorphic images, including the identity, of any of the $\set{G_i}$, nor any copies of $K_4$.
\end{lemma}
\begin{proof}
As neither $A^*$ nor $B^*$ contain triangles, the only triangles in $A^+, B^+$ are in the copies of $W_5$ we have added. Thus there are no copies of $K_4$.

Now suppose there is some homomorphic image $H$ of $G_i$ in $A^+$ or $B^+$. As $H$ does not contain a copy of $K_4$, by Lemma \ref{lemma:k4Appearance} it must contain a necklace of triangles and additional vertices such that any two non-adjacent vertices of the necklace have at least 2 common neighbors. As the only triangles are in copies of $W_5$, the necklace of triangles must contain triangles in 2 distinct copies of $W_5$. However, any points in distinct copies of $W_5$ will not have 2 common neighbors.
\end{proof}

We now shift from the language with unary predicates to the pure graph language. Given the choice of $(G_i, a_i)$ to encode unary predicates, for any choice of tiling problem $\TT$ we get a hereditary graph class $\HH_\TT$, which has the JEP iff $\TT$ has a solution. We wish to add extra constraints to this graph class. In particular we wish to forbid $K_4$ and non-identity homomorphic images of the $\set{G_i}$, for $i \leq 14$. (We choose $i = 14$ because our original construction in a language with unary predicates used 12 unary predicates. We have added another predicate $C_5$ in this section, and require a ``dummy'' predicate for the translation to the pure graph language.) We will call the resulting hereditary graph class $\HH_\TT^+$.

\begin{remark}
Recall the function $\wedge$ from Definition \ref{def:wedge}. The definition depends on the choice of $\set{G_i, a_i}$, which is differs between this section on the previous one. As before, we will use $\widehat{G}$ to denote $\wedge(G)$.
\end{remark}

\begin{lemma} \label{lemma:augmentedCanonical}
Let $\widehat{A^+}, \widehat{B^+}$ be the canonical models in the pure graph language, obtained by applying the function $\wedge$ to $A^+,B^+$. Then $\widehat{A^+}, \widehat{B^+}$ do not contain copies of $K_4$ or any non-identity homomorphic images of the $\set{G_i}$, and so are in $\HH_\TT^+$.
\end{lemma}
\begin{proof}
As $K_4$ and any non-identity homomorphic images of the $\set{G_i}$ not containing $K_4$ are 2-connected, if one of them is contained in $\widehat{A^+}$ or $\widehat{B^+}$ then it must be contained in a single block. We know they are not contained in any of the copies of $\set{G_i}$ attached by $\wedge$ as the $\set{G_i}$ are cores and form an antichain under homomorphisms, so they must have been present in $A^+,B^+$. But by Lemma \ref{lemma:unaryCanonical}, we know this is not the case.
\end{proof}

As we already know $\HH_\TT$ has the JEP when $\TT$ has a solution, to check that $\HH_\TT^+$ has the JEP, it suffices to check that our joint-embedding procedure for $\HH_\TT$ does not create any new copies of $K_4$ or homomorphic images of $\set{G_i}$.

Recall the two steps of our joint-embedding procedure in the pure graph language. First, for every vertex $v$ such that there is no $i$ such that $v$ is in a copy of $G_i$ free over its basepoint, we attach a copy of $G_k$ freely over $v$, which gets identified with the basepoint, where $G_k$ represents a unary predicate specially reserved for this completion process (in our case, $k = 14$). We may then interpret the resulting graph in the language with unary predicates, and in the next step we add edges as we would have done there.

\begin{lemma} \label{lemma:augmentedJEP}
Let $\TT$ be a tiling problem with a solution, and suppose $A, B \in \HH_\TT^+$. Then applying our joint-embedding procedure to $A,B$ creates no homomorphic images of any of the $\set{G_i}_{i \leq 14}$ except for copies of $G_{14}$, nor any copies of $K_4$, and so produces a graph in $\HH_\TT^+$.
\end{lemma}
\begin{proof}
In the first step of our joint-embedding procedure, we add copies of $G_{14}$ freely over various vertices. As $K_4$ and homomorphic images of the $\set{G_i}$ not containing $K_4$ are 2-connected, any new copies of these graphs must appear in the attached copies of $G_{14}$. First, $K_4$ does not embed into $G_{14}$. Then, as $G_{14}$ is a core and the $\set{G_i}$ form an antichain under homomorphisms, the only homomorphic image of any of the $\set{G_i}$ embedding in $G_{14}$ is $G_{14}$ itself. 

Let $A'$ and $B'$ be the graphs obtained from $A$ and $B$ as a result of this first step. As the graphs $\set{G_i}$ and $K_4$ are connected, no copies of $K_4$ or the $\set{G_i}$ are created by passing to the disjoint union $A' \sqcup B'$. We now continue on to the second step of our joint-embedding procedure, in which edges between the factors are added to $A' \sqcup B'$. The key point in this step is that no edge we add is contained in a triangle. This immediately rules out creating any copies of $K_4$.

Now suppose our joint embedding-procedure creates some graph $H$, a homomorphic image of one of the $\set{G_i}$. Let $\phi:G_i \to H$ be a homomorphism. We divide the edges of $H$ into two classes. An edge $(u,v)$ of $H$ will be \textit{old} if $G_i$ contains an edge between some element of $\phi^{-1}(u)$ and some element of $\phi^{-1}(v)$, and otherwise the edge will be \textit{new}.

First, note that as all the edges of $G_i$ are contained in a triangle, the same is true for all the old edges of $H$. Thus our joint-embedding procedure cannot add any old edges. 

Let $H'$ be the graph $H$ with all the new edges removed. Then $H'$ is still a homomorphic image of $G_i$, and must be contained in the disjoint union $A' \sqcup B'$. As $H'$ is connected, it must be contained in one of the factors. As $A', B' \in \HH^+$, this is only possible if $H'$ is a copy of $G_i$. Our joint-embedding procedure will not add edges to any non-basepoint of a copy of $G_i$, so any graph produced by adding edges from $H'$ to other vertices will not be 2-connected and won't contain a copy of $K_4$, and so by Lemma \ref{lemma:2ConOrK4} cannot be a homomorphic image of $G_i$.
\end{proof}

\begin{theorem}
There is no algorithm that, given a finite set of forbidden induced subgraphs, decides whether the corresponding hereditary graph class has the JHP.
  
In particular, given a tiling problem $\TT$, $\HH_\TT^+$ has the JHP iff $\TT$ has a solution.
\end{theorem}
\begin{proof}
First, suppose $\TT$ has a solution. Then by Lemma \ref{lemma:augmentedJEP}, $\HH_\TT^+$ has the JEP, and thus the JHP.

Now suppose $\HH_\TT^+$ has the JHP. Then there is some $C \in \HH_\TT^+$ that $\widehat{A^+}, \widehat{B^+}$ both have homomorphisms into. We now wish to argue any homomorphism of $\widehat{A^+}$ into $C$ must be an embedding, and similarly for $\widehat{B^+}$.

Consider taking a homomorphism of $\widehat{A^+}$ whose image must be in $\HH_\TT^+$.  We cannot identify or add edges between any two basepoints of any of the $\set{G_i}$, as they are either already adjacent or have a copy of $W_5$ freely joined over them, so the identification or new edge would create a copy of $K_4$. We cannot identify any non-basepoint of a copy of one of the $\set{G_i}$ with any point outside of that copy of $G_i$ as that would create an edge incident to the non-basepoint, forbidden by $\HH_1$ (Definition \ref{def:HHConstraints}), unless we identified the entire copy of $G_i$ with another copy of $G_i$; however the latter is forbidden as the basepoints cannot be identified. We also cannot add an edge to a non-basepoint from outside the copy of $G_i$ it is in. Finally, we cannot add edges or identify points within a given copy of one of the $\set{G_i}$, since all non-identity homomorphic images of the $\set{G_i}$ are forbidden. Thus the homomorphisms must be the identity.

Thus $\widehat{A^+}, \widehat{B^+}$ actually joint-embed in $C$, and as in Section \ref{sec:JEPtoTiling} this must encode a solution to $\TT$.
\end{proof}

\chapter{\text{ Questions} on the Model Theory of Permutation Avoidance Classes}\label{chap:PermutationQuestions}
As noted in the thesis introduction, the model-companion of the theory of a permutation avoidance class provides a notion of limit structure to which model-theoretic methods can be applied. In particular, if the model companion is $\omega$-categorical, i.e. has a unique countable model, then there is a canonical countable limit structure generalizing the \fraisse limit. Such a structure would also have much stronger universality properties than those provided by the notions studied in the permutation avoidance literature.

The permutation avoidance literature includes atomicity, equivalent to the JEP, which is the minimal condition for being able to represent a permutation avoidance class in terms of a single structure. It also includes the stronger notion of \textit{natural classes}, in which both the linear orders on the finitely universal limit structure provided by the JEP must have order type $\N$ \cite{natural}, as well as the intermediate notion of \textit{supernatural classes} in which only one of the linear orders need have order type $\N$ \cite{supernatural}.

Rather than just specifying the class of finite substructures of the limit structure, one may ask for a countable limit structure into which every other possible countable limit structure can be embedded. The existence of an $\omega$-categorical model companion is stronger yet. Although not considered in the permutation avoidance literature, the analogous problems for graphs classes specified by forbidden (non-induced) subgraphs have been intensely studied (see \cite{WQOU} for a survey) and lead back to finitary problems. One may hope for similar developments in permutation avoidance classes.

We state most of our questions in terms of the decidability, given finitely many forbidden permutations, of various properties of a permutation avoidance class. One could more broadly ask for reasonable necessary and/or sufficient conditions for these properties, which would also remove the restriction of considering only finitely-constrained classes.

\begin{question}
Is there an algorithm that, given finite set of forbidden permutations, decides whether the theory of the corresponding permutation avoidance class has a model companion?
\end{question}

\begin{question} \label{question:OmegaDecidable}
Is there an algorithm that, given finite set of forbidden permutations, decides whether the theory of the corresponding permutation avoidance class has an $\omega$-categorical model companion?
\end{question}

These two questions have been considered in \cite{LH} in the cases where the forbidden permutations have size at most 3. Based on the evidence there, the following conjecture was put forward. While extravagant, it is consistent with what is currently known.

\begin{conjecture}[\cite{LH}]
If a permutation pattern avoidance class has the joint embedding property and is not finite,
then its first-order theory has an $\omega$-categorical model companion.
\end{conjecture}

If this should be true, then we are led again to Question \ref{question:permJEP} about the decidability of the JEP.

When studying the analog of Question \ref{question:OmegaDecidable} in the case of graph classes forbidding finitely many non-induced subgraphs, the key tool in the analysis is the model-theoretic algebraic closure operator, understood as being taken in an existentially-complete structure.

\begin{definitionC} \label{def:acl}
Let $\GG$ be a class of graphs specified by forbidden subgraphs. Let $G \in G$ and $A \subset G$. Then $v \in acl(A)$ if for any $G^* \in C$ such that $A \subset G^*$, the set of images of $v$ under embeddings of $G$ into $G^*$ over $A$ is finite.
\end{definitionC}

In particular, the class has an $\omega$-categorical model companion iff the algebraic closure operator is locally finite, i.e. $acl(A)$ is finite for finite $A$. 

However, as noted in \cite{WQOU}, considering forbidden induced graphs allows points to be replaced by infinite sets of indiscernibles, trivializing the algebraic closure operator without affecting the existence of an $\omega$-categorical model companion. As we are considering forbidding induced subpermutations, we next consider whether the algebraic closure operator may still be of use.

\begin{question}
Can the notion of algebraic closure from Definition \ref{def:acl}, or some variation of it, be used to determine whether an $\omega$-categorical model companion exists, perhaps in certain restricted permutation avoidance classes?
\end{question}


More broadly than $\omega$-categorical model companions, we may consider the existence of countable universal structures for a class, i.e. countable models into which all other countable models embed. These still provide a notion of limit structure, and if the model companion has a saturated countable model, that would provide a canonical notion. 

General model theory provides conditions for when such a saturated countable model of the model companion exists, see e.g. \cite[Theorem 2]{CSS}.

\begin{theorem}
Let $T$ be a complete theory, with model companion $T^*$. Then the following are equivalent.
\begin{enumerate}
\item The models of $T$ contain a countable universal structure under embedding.
\item The models of $T^*$ contain a countably universal structure under elementary embeddings.
\item There exist countable saturated models of $T^*$.
\item $T^*$ is small, i.e. has countably many $n$-types for each $n$.
\end{enumerate}
\end{theorem}

Based on the study of graph classes, we might expect the questions of the existence of $\omega$-categorical model companions and of countable universal structures to be intimately related.

\begin{question}
Is there an algorithm that, given a finite set of forbidden permutations, decides whether the corresponding permutation avoidance class has a countable universal structure?
\end{question}

Once a limit theory has been provided in the form of a model companion, model theory provides further tools for separating tame behavior from wild. As the structures we consider are ordered, the most appropriate such dividing line is whether the theory is NIP, although the more restrictive notion of distality may also be interesting.

\begin{question}
Is there an algorithm that, given finite set of forbidden permutations such that the theory of the corresponding permutation avoidance class has a model companion, decides whether that model companion is NIP? distal?
\end{question}

\begin{question}
Do the various properties we are trying to decide algorithmically have combinatorial significance?
\end{question}



\end{document}